\documentclass  [12pt]{article}

\newtheorem{theorem}{Theorem}[section]
\newtheorem{lemma}{Lemma}[section]

\newtheorem{conjecture}{Conjecture}[section]

\newtheorem{remark}{Remark}[section]
\newcommand{\eqnsection}{
   \renewcommand{\theequation}{\thesection.\arabic{equation}}
   \makeatletter
   \csname @addtoreset\endcsname{equation}{section}
   \makeatother}
 
\def \ov{\overline}

\def \be{\begin{equation}}
\def \ee{\end{equation}}
\def \bt{\begin{theorem}} 
\def \et{\end{theorem}}
\def \bl{\begin{lemma}} 
\def \el{\end{lemma}}
\def \bea{\begin{eqnarray}}
\def \eea{\end{eqnarray}}
\def \bas{\begin{eqnarray*}}
\def \eas{\end{eqnarray*}}

\def \al{\alpha}

\def \de{\delta}
\def \De{\Delta}
\def \ep{\epsilon}

\def \la{\lambda}

\def \si{\sigma}

\def \th{\theta}

\def \ze{\zeta}

\def \ff{\infty}
\def \wh{\widehat}
\def \wt{\widetilde}

\def \rar{\rightarrow}

\def \cd{\,\cdot\,}

\def \RR{{\cal R}}

\def \E{E}

\def \({\left(}
\def \){\right)}

\def \lc{\left\{}
\def \rc{\right\}}

\def \nn{\nonumber}

\def \bc{\begin{center} }
\def \ec{\end{center} }
\def \bs{\begin{slide} }
\def \es{\end{slide} }

\def\square{{\vcenter{\vbox{\hrule height.3pt
        \hbox{\vrule width.3pt height5pt \kern5pt
           \vrule width.3pt}
        \hrule height.3pt}}}}
\def\qed{{\hfill $\square$ \bigskip}}

\eqnsection
\begin{document}

\def\wh{\widehat}
\def\ol{\overline}

\title{A CLT  for the third integrated moment \\of  Brownian local time increments}

\author{  Jay Rosen\thanks
     {  This research   was   supported, in part, by grants from PSC-CUNY and the National Science
Foundation.}}


\maketitle

\bibliographystyle{amsplain}

\begin{abstract} 
 Let $\{L^{ x }_{ t}\,;\,(x,t)\in R^{ 1}\times  R^{  1}_{ +}\}$ denote the local time of  Brownian motion. Our main result is to show that  for each fixed $t$
\bea
&& { \int ( L^{ x+h}_{t}- L^{ x}_{ t})^{ 3}\,dx-12h\int ( L^{ x+h}_{t }- L^{ x}_{ t })L^{ x}_{ t}\,dx-24h^{2}t\over h^{ 2}}\nonumber\\
&&\hspace{1 in}
\stackrel{\mathcal{L}}{\Longrightarrow}\sqrt{ 192}\( \int ( L^{ x}_{ t})^{ 3}\,dx\)^{1/2}\,\,\eta\nn
\eea
 as $h\rar 0$, where
$\eta$ is a normal random variable with mean zero and variance one that is independent of $L^{ x }_{ t}$. This generalizes our previous result for the second moment. We also explain  
 why our approach will not work for higher moments.

\end{abstract}



\section{Introduction}\label{sec-intro}

   Let $\{L^{ x }_{ t}\,;\,(x,t)\in R^{ 1}\times  R^{  1}_{ +}\}$ denote the local time of  Brownian motion.   
 Let
\begin{equation}
\al_{p, t}=\int ( L^{ x}_{ t})^{ p}\,dx\label{5.1w}
\end{equation} 
   (an integral sign without limits is to be read as $\int_{-\ff}^{\ff}$,) and let $\eta=N(0,1)$ be  
independent of $\al_{p, t}$. The main result of  \cite{CLMR} is the following weak limit theorem.

\begin{theorem}\label{theo-clt2} For each fixed $t$
\begin{equation} { \int ( L^{ x+h}_{t}- L^{ x}_{ t})^{ 2}\,dx- 4ht\over h^{ 3/2}}
\stackrel{\mathcal{L}}{\rightarrow}c\sqrt{\al_{2,t}}\,\,\eta\label{5.0weak}
\end{equation} 
 as $h\rar 0$,  where $c=\( 64 / 3 \)^{ 1/2}$. 
 
 Equivalently
\begin{equation} { \int ( L^{ x+1}_{t}- L^{ x}_{ t})^{ 2}\,dx- 4t\over t^{ 3/4}}
\stackrel{\mathcal{L}}{\rightarrow}c\sqrt{\al_{2,1}}\,\,\eta\label{5.0tweak}
\end{equation}
as $t\rar\ff$.
\end{theorem}

\medskip	
In this paper we provide  the analogous CLT for the third power.
\begin{theorem}\label{theo-trip}
 For each fixed $t$
\begin{equation} { \int ( L^{ x+h}_{t}- L^{ x}_{ t})^{ 3}\,dx-12h\int ( L^{ x+h}_{t }- L^{ x}_{ t })L^{ x}_{ t}\,dx-24h^{2}t\over h^{ 2}}
\stackrel{\mathcal{L}}{\rightarrow}c\sqrt{\al_{3, t}}\,\,\eta\label{1.1}
\end{equation} 
 as $h\rar 0$,  where $c=\sqrt{ 192}$. 
 
 Equivalently
\begin{equation} { \int ( L^{ x+1}_{t}- L^{ x}_{ t})^{ 3}\,dx -12\int ( L^{ x+1}_{t }- L^{ x}_{ t })L^{ x}_{ t}\,dx-24t^{2}\over t}
\stackrel{\mathcal{L}}{\rightarrow}c\sqrt{\al_{3, 1}}\,\,\eta\label{1.2}
\end{equation}
as $t\rar\ff$.	
\end{theorem}

We explain below why the approach we use will not work for moments larger than three.

The equivalence of (\ref{1.1}) and (\ref{1.2}) 
 follows from the scaling   relationship 
  \begin{equation}
\{ L^{ x}_{ h^{-2}t}\,;\,( x,t)\in R^{ 1}\times R^{ 1}_{ +}\}
\stackrel{\mathcal{L}}{=}\{h^{ -1} L^{h x}_{ t}\,;\,( x,t)\in R^{ 1}\times
R^{ 1}_{ +}\},\label{scale}
\end{equation}
see e.g. \cite[Lemma 10.5.2]{book},
which implies that 
\begin{equation}
\int ( L^{ x+h}_{t}- L^{ x}_{ t})^{ 3}\,dx\stackrel{\mathcal{L}}{=}h^{4} \int ( L^{ x+1}_{t/h^{2}}- L^{ x}_{ t/h^{2}})^{ 3}\,dx,\label{scl1}
\end{equation}
and 
\begin{equation}
\int ( L^{ x+h}_{t}- L^{ x}_{ t})L^{ x}_{ t} \,dx\stackrel{\mathcal{L}}{=}h^{3} \int ( L^{ x+1}_{t/h^{2}}- L^{ x}_{ t/h^{2}})L^{ x}_{ t/h^{2}} \,dx.\label{scl2}
\end{equation}
  Using this, and (\ref{1.1}) with $t=1$, and then setting $h^{2}=1/t$ gives (\ref{1.2}).

\medskip	

Theorem \ref{theo-trip} is derived using the method of moments.  
Note that the right hand side of (\ref{1.1}) is $c\sqrt{\al_{3, t}}\,\,\eta$ .
Unfortunately, we can only show that  $\sqrt{\al_{p, t} }\,\,\eta$ is determined by its moments   if $p=2$ or $3$,  so we cannot use our approach to prove an analog of Theorem \ref{theo-trip} for moments larger than three.

In Section \ref{sec-est} we give  some estimates on the potential densities and transition densities  of Brownian motion which are used throughout this paper. Their proof is deferred until Section \ref{sec-Prooflemvprop}.  In Section \ref{sec-BM} we show how Theorem \ref{theo-clt2} will  follow from a result, Lemma \ref{lem-2weak}, on the moments of an analogous expression where   $t$ is replaced by an independent exponential time. This
 Lemma is proven in Section \ref{sec-expmom}. Other lemmas are that  used in the proof of  Theorem \ref{theo-clt2} are derived in Sections \ref{sec-6.1}-\ref{sec-var}.
 
  This paper extends the basic approach used in \cite{CLMR}. The main novelty in this paper is the need to subtract a non-random term in (\ref{1.1}) in order to get a Central Limit Theorem. Dealing with this non-random  subtraction term, and in particular the need to keep track of delicate cancellations, makes this paper considerably more difficult than \cite{CLMR}. Although, as mentioned, the approach of the present paper will not work for higher moments, Theorem \ref{theo-clt2} does suggest what a Central Limit Theorem for higher moments should look like. Here is our conjecture for the fourth integrated moment.

 \begin{conjecture} For each fixed $t$ 
 \bea
&&{ \int ( \De^{h}L^{ x}_{ t})^{ 4}\,dx-24h\int  ( \De^{h}L^{ x}_{ t})^{ 2}L^{ x}_{ t}\,dx+48h^{2}\int  (L^{ x}_{ t})^{ 2}- ( \De^{h}L^{ x}_{ t})L^{ x}_{ t}  \,dx\over h^{ 5/2}}\nonumber\\
&&\hspace{3 in}
\stackrel{\mathcal{L}}{\rightarrow}c_{4}\sqrt{\al_{4,t}}\,\,\eta\label{p2.4}
\eea
  as $h\rar 0$,where $c_{q}=\sqrt{ {2^{2q+1} q!\over q+1}}$ and $\De^{h}L^{ x}_{ t}=L^{ x+h}_{t}- L^{ x}_{ t}$.

\end{conjecture}

\section{Estimates for  the   potential density of \newline Brownian motion}\label{sec-est}

 Let $p_{t}(x)$ denote the   probability density of Brownian motion. 
The $\al$-potential density of Brownian motion
 \begin{equation}
u^{\al}(x)=\int_{0}^{\ff}e^{-\al t}p_{t}(x)\,dt={e^{-\sqrt{2\al}|x|} \over \sqrt{2\al}}\label{pot.1w}.
\end{equation}
Let $\la_{\al}$ be an independent exponential random variable with mean $1/\al$.

 Kac's moment formula, \cite[Theorem 3.10.1]{book}, states that
\begin{equation} E^{ x_{ 0}}\(\prod_{ j=1}^{ n}L^{ x_{ j}}_{ \la_{\al}} \)=\sum_{
\pi}\prod_{ j=1}^{ n}u^{\al}( x_{\pi( j)}-x_{\pi( j-1)})\label{1.2w}
\end{equation} where the sum runs over all permutations $\pi$ of $\{ 1,\ldots,
n\}$ and  
$\pi(0)=0.$

Let $\De_{ x}^{ h}$  denote
the finite difference operator on the variable $x$, i.e.
\begin{equation}
\De_{x}^{ h}\,f(x)=f(x+h)-f(x).\label{pot.3w}
\end{equation}
We write $\De^{ h}$ for $\De_{x}^{ h}$ when the variable $x$ is clear.

\medskip	The next lemma collects some facts about $u^{\al}(x)$ that are used in this paper.

\begin{lemma}\label{lem-vprop}Fix $\al >0$. For $0<h\leq 1$,
\bea
\De_{ x}^{ h}\De_{ y}^{ h} u^{\al}(x-y)\Bigg\vert_{ y=x}&= &2\({1-e^{-\sqrt{2\al}\,h} \over \sqrt{2\al}}\)=2h+O( h^{ 2}),\label{1.8}
\\\nn\\
v^{\al}(x)=:|\De ^{ h}\,u^{\al}(x)|&\leq &Ch\, u^{\al}( x),\label{1.3x}
\\\nn\\
w^{\al}(x)=:|\De^{ h}\De^{ -h} u^{\al}(x )|&\leq&  \left\{\begin{array}{ll}
C h \, u^{\al}( x ),   \\\\
C h^{2}\, u^{\al}( x ), \hspace{.2 in}\forall \,|  x|\geq h. 
\end{array}
\right. \label{1.3y}
\eea
We have
\begin{equation}
\int  \(w^{\al}(x)\)^{q}\,dx=O( h^{ q+1}) \label{li.13}
\end{equation}
and 
\begin{equation}
\int_{|x|\geq h} \(w^{\al}(x)\)^{q}\,dx=O( h^{ 2q}).\label{1.30gb}
\end{equation}
In addition, for any $q\geq 2$
\bea
&&
\int  \(\De^{ h}\De^{ -h}\,u^{\al}(x)\)^{q} \,dx=( 2^{q+1}/(q+1)+O( h))h^{ q+1},\label{1.30g}
\eea
In all these statements the constants $C$ and the terms $O( h^{ \cd})$ may depend on $\al$.
\end{lemma} 

The proof is provided in Section \ref{sec-Prooflemvprop}.

 We note that the same proof shows that for any $\al_{1}, \ldots, \al_{q}>0$
 \begin{equation}
 \int  \prod_{i=1}^{q} \(\De^{ h}\De^{ -h}\,u^{\al_{i}}(x)\)\,dx=( 2^{q+1}/(q+1)+O( h))h^{ q+1}.\label{mult.u}
 \end{equation}
 
  \begin{remark}
 {\rm   In Lemma \ref{lem-vprop} we have taken $h$ positive. 
  Using the fact that $u^{\al}(x)$ is an even function of $x$ it is easy to check that we obtain the analog of (\ref{1.3x}) for all $|h|\leq 1$  if on the right hand side   we replace $h$ by $|h|$.}
  \end{remark}

  The following estimates, which will be used in the proof of Lemma \ref{lem-3.1j}, are also proven in Section \ref{sec-Prooflemvprop}.
  
   \begin{lemma}\label{lem-vpropt}Let $0<h\leq 1$ and $0<T<\ff$. Then for some $C_{T}<\ff$
 \begin{equation}
u_{T}( x)=:\int_{0}^{T} \,p_{t}(x)\,dt\leq C_{T}\, e^{-|x|},\label{9.300}
\end{equation} 
\begin{equation}
v_{T}( x)=:\int_{0}^{T} |\De ^{ h}\,p_{t}(x)|\,dt\leq C_{T}h\, e^{-|x|},\label{9.3x}
\end{equation} 
and
 \begin{equation}
w_{T}(x)=:\int_{0}^{T} |\De^{ h}\De^{ -h} p_{t}(x )|\,dt\le C_{T}h^{2}\frac{e^{-x^{2}/32T}}{|x|},\hspace{.2 in} 
|x|\geq 2h.\label{9.3w}
\end{equation}
Also
\begin{equation} 
\int w_{T}(x) \,dx\leq C_{T}h^{ 2}|\log h|,\label{9.13t}
\end{equation}
and for any $q\geq 2$
\begin{equation} 
\int w_{T}^{q}(x)\,dx\leq C_{T} h^{ q+1},\label{9.30g}
\end{equation}
and
\begin{equation} 
\int_{|x|\geq \sqrt{h}} \,w_{T}^{q}(x) \,dx\leq C_{T}h^{3q/2+1/2}.\label{9.30gb}
\end{equation}
\end{lemma} 

   \begin{lemma}\label{lem-vpropd}Let $0<h\leq 1$ and $0<\de<T<\ff$. Then for some $C_{\de,T}<\ff$
 \begin{equation}
u_{T}( x)=:\sup_{\de\leq t\leq T} \,p_{t}(x) \leq C_{\de,T}\, e^{-x^{2}/2T},\label{d9.300}
\end{equation} 
\begin{equation}
v_{T}( x)=:\sup_{\de\leq t\leq T} |\De ^{ h}\,p_{t}(x)| \leq C_{\de,T}h\, e^{-x^{2}/2T}.\label{d9.3x}
\end{equation} 
and
 \begin{equation}
w_{T}(x)=:\sup_{\de\leq t\leq T} |\De^{ h}\De^{ -h} p_{t}(x )|  \le C_{\de,T}h^{2}e^{-x^{2}/2T}.\label{d9.3w}
\end{equation}
\end{lemma}

 \begin{lemma}\label{lem-big}Let $0<h\leq 1$. 
For $q\geq 2$
\begin{equation} 
\int  \(\int_{0}^{\ff} \De^{ h}\De^{ -h} p_{t}(x )\,dt\)^{q}\,dx=( 2^{q+1}/(q+1)+O( h^{1/2}))h^{ q+1},\label{big.1}
\end{equation}
and
\begin{equation} 
\int  \(\int_{0}^{h} \De^{ h}\De^{ -h} p_{t}(x )\,dt\)^{q}\,dx=( 2^{q+1}/(q+1)+O( h^{1/2}))h^{ q+1}.\label{big.2}
\end{equation}
\end{lemma}

\section{Proof of Theorem \ref{theo-trip}}\label{sec-BM}

 Theorem \ref{theo-trip} will follow from the next lemma.

\bl\label{lem-6.1} For each integer $m\ge 0$ and $t\in R_{+}$
\begin{eqnarray} &&
\lim_{ h\rar 0}E\(\({ \int ( L^{ x+h}_{t}- L^{ x}_{ t})^{
3}\,dx- 12h\int L^{ x}_{ t}( L^{ x+h}_{t}- L^{ x}_{ t})\,dx -24h^{2}t\over h^{ 2}}\)^{m}\)\nn\\ 
&&\hspace{ .5in}  =\left\{\begin{array}{ll}
\displaystyle{  ( 2n)!\over 2^{ n}n!}\( 192\)^{ n} E\lc\(\int (L^{ x}_{ t})^{ 3}\,dx\)^{
n}\rc &\mbox{ if }m=2n\\\\
0&\mbox{ otherwise.}
\end{array}
\right.
\label{7.53}
\end{eqnarray}
\el

\noindent{\bf  Proof of Theorem \ref{theo-trip} }   It follows from  \cite[(6.12)]{CLR}   that   for any $q$
 \begin{equation}
 E\lc\(\int (L^{ x}_{ t})^{ q}\,dx\)^{ n}\rc\leq C_{t}^{ n}( n!)^{ (q-1)/2}.\label{rb.1}
 \end{equation}
 Therefore, since $ \sqrt{( 2n)!}\le 2^{ n}n! $,
the right hand side of (\ref{7.53}), which is the $2n-$th  moment of  $\wt C\sqrt{\int (L^{ x}_{ t})^{ 3}\,dx}\,\,\eta$ is bounded above by $  \wt C^{2n}C^{ n}(2n)! $. 
This implies that   $\wt C\sqrt{\int (L^{ x}_{ t})^{ 3}\,dx}\,\,\eta$ is determined by its moments; (see \cite[p. 227-228]{Feller}).
Lemma \ref{lem-6.1} together with the  method of moments, \cite[Theorem 30.2]{B}, then  gives us (\ref{1.1}).\qed

\noindent{\bf Proof of Lemma \ref{lem-6.1} } 
Let $\la_{\ze}$ be an exponential random variable with mean $1/\ze$.
It follows from Lemma \ref{lem-2weak} below that for each integer $m\ge 0$,
\begin{eqnarray} &&
\lim_{ h\rar 0}E\(\({ \int (\De_{x}^{h}L^{ x}_{ \la_{\ze}})^{
3}\,dx- 12h \int L^{ x}_{\la_{\ze}}( \De_{x}^{h}L^{ x}_{ \la_{\ze}})\,dx -24h^{2}\la_{\ze}\over h^{ 2}}\)^{m}\)\nn\\ 
&&\hspace{ .5in}  =\left\{\begin{array}{ll}
  \displaystyle  {( 2n)!\over 2^{ n}n!}\(  192\)^{ n} E\lc\(\int (L^{ x}_{ \la_{\ze}})^{ 3}\,dx\)^{
n}\rc &\mbox{ if }m=2n\\\\
0&\mbox{ otherwise.}
\end{array}
\right.
\label{7.54a}
\end{eqnarray}

We  write (\ref{7.54a}) as 
\begin{eqnarray} &&
\int_{0}^{\ff}e^{- \ze s } E\(\({ \int (\De_{x}^{h}L^{ x}_{s})^{
3}\,dx- 12h \int L^{ x}_{s}( \De_{x}^{h}L^{ x}_{s})\,dx-24h^{2}s \over h^{ 2}}\)^{m}\) \,ds 
\nn \\
&&\qquad\longrightarrow\int_{0}^{\ff} e^{- \ze s }E\lc  \eta^{m}\( 192\int  (L^{ x}_{s})^{3}   \,dx\,\,\)^{
m/2}\rc \,ds\label{77.1d}
\end{eqnarray}
as $h\rar 0$. 
For $h>0 $ let  
\be 
\wh F_{m,h}(s):=  E\(\({ \int (\De_{x}^{h}L^{ x}_{s})^{
3}\,dx- 12h \int L^{ x}_{s}( \De_{x}^{h}L^{ x}_{s})\,dx -24h^{2}s\over h^{ 2}}\)^{m}\), \label{77.1e}
\ee
and set
\begin{eqnarray} &&
\wh F_{m,0}(s):= E\lc  \eta^{m}\(  192\int  (L^{ x}_{s})^{3}   \,dx\,\,\)^{
m/2}\rc. \label{77.1e}
\end{eqnarray}
Then (\ref{77.1d}) can be written as 
\begin{equation}
\lim_{ h\rar 0}\int_{0}^{\ff}e^{- \ze s } \wh {F}_{m,h}(s) \,ds =\int_{0}^{\ff}e^{- \ze s } \wh {F}_{m,0}(s) \,ds.\label{77.4}
\end{equation}

We consider first the case when $m$ is even and write   $m=2n$. In this case  $\wh {F}_{2n,h}(s)\geq 0$ and the extended continuity theorem 
\cite[XIII.1, Theorem 2a]{Feller} applied to  (\ref{77.4}) implies that
\begin{equation}
\lim_{ h\rar 0}\int_{0}^{t}  \wh {F}_{2n,h}(s) \,ds =\int_{0}^{t}  \wh {F}_{2n,0}(s) \,ds\label{77.5}
\end{equation}
for all $t$. In particular,
\begin{equation}
\lim_{ h\rar 0}\int_{t}^{t+\de}  \wh {F}_{2n,h}(s) \,ds =\int_{t}^{t+\de}  \wh {F}_{2n,0}(s) \,ds.\label{77.6}
\end{equation}
It follows from the fact that $L^{ x}_{s}$ is almost surely continuous and increasing in $s$ that $ \wh {F}_{2n,0}(s)$ is continuous in $s$. 
(We saw in (\ref{rb.1}) that it is bounded.) Consequently,
\begin{equation}
\lim_{ \de\rar 0}\lim_{ h\rar 0}{1 \over \de}\int_{t}^{t+\de}  \wh {F}_{2n,h}(s) \,ds = \wh {F}_{2n,0}(t). \label{77.6}
\end{equation}
When $t=0$ we get
\begin{equation}
\lim_{\de\to 0}
\lim_{h\to 0}{1\over \de}\int_0^{\de}\wh {F}_{2n,h}(s)ds=0.\label{77.13}
\end{equation}

To obtain (\ref{7.54a}) when $m$ is even we must show that
\begin{equation}
\lim_{ h\rar 0}   \wh {F}_{2n,h}(t)  = \wh {F}_{2n,0}(t) \label{77.15}.
\end{equation}
This follows from  (\ref{77.6}) once we  show that
\begin{equation}
\lim_{ \de\rar 0}\lim_{ h\rar 0}{1 \over \de}\int_{t}^{t+\de}  \wh {F}_{2n,h}(s) \,ds = \lim_{ h\rar 0}   \wh {F}_{2n,h}(t)  . \label{77.6j}
\end{equation}
We proceed to obtain (\ref{77.6j}).

For $s\geq t$ we  write 
\begin{eqnarray}
&&\int (\De_{x}^{h}L^{ x}_{s})^{
3}\,dx=\int \(\De_{x}^{h}L^{ x}_{t}+\De_{x}^{h}(L^{ x}_{s}-L^{ x}_{t})\)^{
3}\,dx
\label{fu.1}\\
&&=\int (\De_{x}^{h}L^{ x}_{t})^{
3}\,dx +3\int (\De_{x}^{h}L^{ x}_{t})^{
2}\De_{x}^{h}(L^{ x}_{s}-L^{ x}_{t})\,dx  \nonumber\\
&& +3\int \De_{x}^{h}L^{ x}_{t}\(\De_{x}^{h}(L^{ x}_{s}-L^{ x}_{t})\)^{
2}\,dx  +\int  \(\De_{x}^{h}(L^{ x}_{s}-L^{ x}_{t})\)^{
3}\,dx\nonumber
\end{eqnarray}
and
\begin{eqnarray}
&&\int L^{ x}_{s}( \De_{x}^{h}L^{ x}_{s})\,dx =\int (L^{ x}_{t}+(L^{ x}_{s}-L^{ x}_{t}))\(\De_{x}^{h}L^{ x}_{t}+\De_{x}^{h}(L^{ x}_{s}-L^{ x}_{t})\)\,dx 
\nn\\
&& =\int L^{ x}_{t}( \De_{x}^{h}L^{ x}_{t})\,dx +\int L^{ x}_{t} \De_{x}^{h}(L^{ x}_{s}-L^{ x}_{t}) \,dx   \label{fu.2}\\
&& +  \int (L^{ x}_{s}-L^{ x}_{t})( \De_{x}^{h}L^{ x}_{t})\,dx +\int (L^{ x}_{s}-L^{ x}_{t}) \De_{x}^{h}(L^{ x}_{s}-L^{ x}_{t}) \,dx \nonumber
\end{eqnarray}
so that
\begin{eqnarray}
&&\int (\De_{x}^{h}L^{ x}_{s})^{
3}\,dx-12h\int L^{ x}_{s}( \De_{x}^{h}L^{ x}_{s})\,dx-24h^{2}s
\label{fu.3}\\
&&=\int (\De_{x}^{h}L^{ x}_{t})^{
3}\,dx-12h \int L^{ x}_{t}( \De_{x}^{h}L^{ x}_{t})\,dx-24h^{2}t\nonumber\\
&&+3\int \lc  (\De_{x}^{h}L^{ x}_{t})^{
2}-4h  L^{ x}_{t}        \rc\De_{x}^{h}(L^{ x}_{s}-L^{ x}_{t})\,dx  \nonumber\\
&& +3\int \De_{x}^{h}L^{ x}_{t}\lc   \(\De_{x}^{h}(L^{ x}_{s}-L^{ x}_{t})\)^{
2} -4h   (L^{ x}_{s}-L^{ x}_{t})    \rc\,dx \nonumber\\
&& +\int  \(\De_{x}^{h}(L^{ x}_{s}-L^{ x}_{t})\)^{
3}\,dx-12h\int (L^{ x}_{s}-L^{ x}_{t}) \De_{x}^{h}(L^{ x}_{s}-L^{ x}_{t}) \,dx-24h^{2}(t-s).\nonumber
\end{eqnarray}

Note that, using $\wt B_{t}$ to denote an independent Brownian motion, and then using  translation invariance 
\begin{eqnarray}
&&\int  \(\De_{x}^{h}(L^{ x}_{s}-L^{ x}_{t})\)^{
3}\,dx-12h\int (L^{ x}_{s}-L^{ x}_{t}) \De_{x}^{h}(L^{ x}_{s}-L^{ x}_{t}) \,dx-24h^{2}(t-s)
\nn\\
&&=\int  \(\De_{x}^{h}L^{ x}_{s-t} \)^{
3}\circ\th_{t}\,\,dx-12h\int \lc L^{ x}_{s-t}\,\,  \(\De_{x}^{h}L^{ x}_{s-t} \)\rc  \circ\th_{t}\,\,dx -24h^{2}(t-s) \nonumber\\
&&\stackrel{law}{=} \int  \(\De_{x}^{h}L^{x -\wt B_{t}}_{s-t} \)^{
3}\,dx-12h\int L^{ x-\wt B_{t}}_{s-t} \,\, \(\De_{x}^{h}L^{ x-\wt B_{t}}_{s-t}\)  \,dx-24h^{2}(t-s)  \nonumber\\
&&=\int (\De_{x}^{h}L^{ x}_{s-t})^{
3}\,dx-12h\int L^{ x}_{s-t}( \De_{x}^{h}L^{ x}_{s-t})\,dx-24h^{2}(t-s).\label{ti.1}
\end{eqnarray}
Also, using $\wt L^{x}_{t}$ to denote an independent copy of Brownian local time
\begin{eqnarray}
&&\int \lc  (\De_{x}^{h}L^{ x}_{t})^{
2}-4h  L^{ x}_{t}        \rc\De_{x}^{h}(L^{ x}_{s}-L^{ x}_{t})\,dx 
\label{ti.2}\\
&& =\int \lc  (\De_{x}^{h}L^{ x}_{t})^{
2}-4h  L^{ x}_{t}        \rc \(\De_{x}^{h}L^{ x}_{s-t}\circ\th_{t}\) \,dx   \nonumber\\
&& \stackrel{law}{=} \int \lc  (\De_{x}^{h}L^{ x}_{t})^{
2}-4h  L^{ x}_{t}        \rc  \De_{x}^{h}\wt L^{ x-B_{t}}_{s-t}  \,dx   \nonumber\\
&& = \int \lc  (\De_{x}^{h}L^{ x+B_{t}}_{t})^{
2}-4h  L^{ x+B_{t}}_{t}        \rc  \De_{x}^{h}\wt L^{ x}_{s-t}  \,dx   \nonumber\\
&&\stackrel{law}{=}\int \lc  (\De_{x}^{h}L^{ x}_{t})^{
2}-4h  L^{ x}_{t}        \rc  \De_{x}^{h}\wt L^{ x}_{s-t}  \,dx   \nonumber
\end{eqnarray}
where we have used the fact that $\{L^{ x+B_{t}}_{t}\,,\,x\in R^{1}\} \stackrel{law}{=}
\{L^{ x}_{t}\,,\,x\in R^{1}\}$ which follows from time reversal.
Similarly,
\begin{eqnarray}
&&\int \De_{x}^{h}L^{ x}_{t}\lc   \(\De_{x}^{h}(L^{ x}_{s}-L^{ x}_{t})\)^{
2} -4h   (L^{ x}_{s}-L^{ x}_{t})    \rc\,dx
\label{ti.3}\\
&& \stackrel{law}{=}  \int \De_{x}^{h}\wt L^{ x}_{t}\lc   \(\De_{x}^{h}(L^{ x}_{s-t} \)^{
2} -4h   L^{ x}_{s-t}     \rc\,dx. \nonumber
\end{eqnarray}

Let
\bea
\wh {G}_{m, h}(t,r)&=&:E\(h^{-2}\int \lc  (\De_{x}^{h}L^{ x}_{t})^{
2}-4h  L^{ x}_{t}       \rc\De_{x}^{h}\wt L^{ x}_{r} \,dx\)^{m}\label{fu.4}
\eea
and set
\begin{eqnarray} &&
\wh {G}_{m, 0}(t,r):= E\lc  \eta^{m}\(  64\int  (L^{ x}_{t})^{2}\wt L^{ x}_{r}   \,dx\,\,\)^{
m/2}\rc. \label{77.1eg}
\end{eqnarray}
We then use the triangle inequality with respect to the norm $\|\,\cdot\,\|_{2n}$ together with (\ref{fu.3})-(\ref{ti.3}) to see that
\begin{eqnarray} 
\wh {F}_{2n, h}^{1/(2n)}(s)&\leq & \wh {F}_{2n, h}^{1/(2n)}(t)+\wh {F}_{2n, h}^{1/(2n)}(s-t)
 \label{fu.7aa}\\
&+ &   
3\wh {G}^{1/(2n)}_{2n, h}(t,s-t)  +3\wh {G}^{1/(2n)}_{2n, h}(s-t,t). \nn 
\end{eqnarray}
Similarly we have
\begin{eqnarray} 
\wh {F}_{2n, h}^{1/(2n)}(s)&\geq & \wh {F}_{2n, h}^{1/(2n)}(t)-\wh {F}_{2n, h}^{1/(2n)}(s-t)
 \label{fu.7ab}\\
&- &
3\wh {G}^{1/(2n)}_{2n, h}(t,s-t) (t)-3\wh {G}^{1/(2n)}_{2n, h}(s-t,t). \nn 
\end{eqnarray}
We now use the triangle inequality   with respect to the norm in $L^{2n}([t,t+\de],   \de^{-1}\,ds)$ to see that
\begin{eqnarray} 
&&
 \lc {1 \over \de}\int_{t}^{t+\de}\wh {F}_{2n, h}(s)\,ds\rc^{1/(2n)}
 \label{fu.7ac}\\   && \leq  \wh {F}_{2n, h}^{1/(2n)}(t)+ \lc{1 \over \de}\int_{t}^{t+\de}\wh {F}_{2n, h}(s-t)\,ds\rc^{1/(2n)}      \nonumber\\
&& +
3 \lc{1 \over \de}\int_{t}^{t+\de}\wh {G}_{2n, h}(t,s-t) \,ds\rc^{1/(2n)}  +3 \lc{1 \over \de}\int_{t}^{t+\de}\wh {G}_{2n, h}(s-t,t)  \,ds\rc^{1/(2n)} \nn 
\end{eqnarray}
and
\begin{eqnarray} 
&&
 \lc {1 \over \de}\int_{t}^{t+\de}\wh {F}_{2n, h}(s)\,ds\rc^{1/(2n)}
 \label{fu.7ad}\\   && \geq  \wh {F}_{2n, h}^{1/(2n)}(t)- \lc{1 \over \de}\int_{t}^{t+\de}\wh {F}_{2n, h}(s-t)\,ds\rc^{1/(2n)}      \nonumber\\
&& -
3 \lc{1 \over \de}\int_{t}^{t+\de}\wh {G}_{2n, h}(t,s-t) \,ds\rc^{1/(2n)}  -3 \lc{1 \over \de}\int_{t}^{t+\de}\wh {G}_{2n, h}(s-t,t)  \,ds\rc^{1/(2n)} \nn 
\end{eqnarray}
Hence, in light of (\ref{77.13}), to prove (\ref{77.6j}) we need only show that for each $t$
\begin{equation}
\lim_{ \de\rar 0}\lim_{ h\rar 0}{1 \over \de}\int_{0}^{\de} \wh {G}_{2n, h}(t,s ) \,ds = 0 \label{fu.8}
\end{equation}
and
\begin{equation}
\lim_{ \de\rar 0}\lim_{ h\rar 0}{1 \over \de}\int_{0}^{\de} \wh {G}_{2n, h}(s,t)  \,ds = 0. \label{fu.9}
\end{equation}

We use $E^{y,z}(\cdot)$ to denote expectation with respect to the independent Brownian motions $B_{t}$ starting at $y$ and $\wt B_{t}$ starting at $z$.Let $\la_{\ze}, \la_{\ze'}$ be independent  exponential random variables with mean $1/\ze, 1/\ze'$.
We show in  Lemma \ref{lem-6.3a} below  that for each integer $n\ge 0$,
\begin{eqnarray} &&
\lim_{ h\rar 0}E^{y,z}\(\({ \int \lc  (\De_{x}^{h}L^{ x}_{\la_{\ze}})^{
2}-4h  L^{ x}_{\la_{\ze}}       \rc\De_{x}^{h}\wt L^{ x}_{\la_{\ze'}}\,dx \over h^{ 2}}\)^{2n}\)\label{6.3lem}\\ 
&&\hspace{ .5in}  =  {( 2n)!\over 2^{ n}n!}\(   \displaystyle   64  \)^{ n} E^{y,z}\lc\(\int (L^{ x}_{ \la_{\ze}})^{ 2}\wt L^{ x}_{ \la_{\ze'}}\,dx\)^{
n}\rc 
\nn
\end{eqnarray}
 uniformly in $y,z.$

Just as (\ref{77.4}) implied (\ref{77.5}), it follows from (\ref{6.3lem}) that 
\begin{equation}
\lim_{ h\rar 0}\int_{0}^{t}\int_{0}^{q}  \wh {G}_{2n,h}(s,r) \,dr\,ds= \int_{0}^{t}\int_{0}^{q}  \wh {G}_{2n,0}(s,r) \,dr\,ds\label{fu.10}
\end{equation}
for all $t$. In particular,
\begin{equation}
\lim_{ h\rar 0}\int_{t}^{t+\de}\int_{q}^{q+\de'}  \wh {G}_{2n,h}(s,r)\,dr\,ds =\int_{t}^{t+\de}\int_{q}^{q+\de'}   \wh {G}_{2n,0}(s,r) \,dr\,ds.\label{fu.11}
\end{equation}
It follows as with $ \wh {F}_{2n,0}(s)$ that $ \wh {G}_{2n,0}(s,r)$ is continuous in $s,r$. Consequently,
\begin{equation}
\lim_{ \de,\de'\rar 0}\lim_{ h\rar 0}{1 \over \de\de'} \int_{t}^{t+\de}\int_{q}^{q+\de'}  \wh {G}_{2n,h}(s,r) \,dr\,ds=\wh {G}_{2n,0}(t,q). \label{fu.12}
\end{equation}
When $t=0$ we get
\begin{equation}
\lim_{ \de,\de'\rar 0}\lim_{ h\rar 0}{1 \over \de\de'} \int_{0}^{ \de}\int_{q}^{q+\de'}  \wh {G}_{2n,h}(s,r)\,dr\,ds=0.\label{fu.13}
\end{equation}
Similarly we have
\begin{equation}
\lim_{ \de,\de'\rar 0}\lim_{ h\rar 0}{1 \over \de\de'} \int_{t}^{t+ \de}\int_{0}^{ \de'}  \wh {G}_{2n,h}(s,r) \,dr\,ds=0.\label{fu.14}
\end{equation}

For $s\geq t$ we  write  
\bea
&& \int \lc  (\De_{x}^{h}L^{ x}_{s})^{2}-4hL^{ x}_{s}      \rc \De_{x}^{h}\wt L^{ x}_{r} \,dx \label{fu.4}\\
&& =  \int \lc  \(\De_{x}^{h}L^{ x}_{t}+\De_{x}^{h}(L^{ x}_{s}-L^{ x}_{t})\)^{2}-4h
\( L^{ x}_{t}+ (L^{ x}_{s}-L^{ x}_{t})\)  \rc\De_{x}^{h}\wt L^{ x}_{r} \,dx      \nonumber\\
&& =  \int \lc  \(\De_{x}^{h}L^{ x}_{t}\)^{2}-4h L^{ x}_{t}  \rc\De_{x}^{h}\wt L^{ x}_{r} \,dx      \nonumber\\
&& +\int \lc  \(\De_{x}^{h}(L^{ x}_{s}-L^{ x}_{t})\)^{2}-4h\(L^{ x}_{s}-L^{ x}_{t})\)  \rc \De_{x}^{h}\wt L^{ x}_{r} \,dx\\
&& +2\int   \De_{x}^{h}L^{ x}_{t}  \,\De_{x}^{h}(L^{ x}_{s}-L^{ x}_{t}) \, \De_{x}^{h}\wt L^{ x}_{r} \,dx      \nonumber
\eea
Hence as before we obtain 
\begin{eqnarray}
&&\lc  {1 \over \de\de'}\int_{t}^{t+ \de}\int_{0}^{ \de'}\wh {G}_{2n, h}(s,r)\,dr\,ds\rc^{1/(2n)} \geq  \lc{1 \over  \de'} \int_{0}^{ \de'} \wh {G}_{2n, h}(t,r)  \,dr\rc^{1/(2n)} \nonumber\\
&&\hspace{1 in}-
\lc{1 \over \de\de'} \int_{t}^{t+ \de}\int_{0}^{ \de'}\bar {G}_{2n, h}(s-t,r)\,dr\,ds\rc^{1/(2n)} 
\label{fu.5}\\
&&\hspace{1 in} -\lc {2 \over \de\de'} \int_{t}^{t+ \de}\int_{0}^{ \de'}\wh {H}_{2n, h}(t,s-t,r)\,dr\,ds\rc^{1/(2n)} \nonumber
\end{eqnarray}
where
\bea
&&
\bar {G}_{m, h}(s-t,r)\label{fu.4d}\\
&&=:E\(h^{-2}\int \lc  (\De_{x}^{h}L^{ x}_{s-t})^{
2}-4h  L^{ x}_{s-t}       \rc\circ\th_{t}\De_{x}^{h}\wt L^{ x}_{r} \,dx\)^{m}\nn\\
&&=\int E^{y,0}\lc \(h^{-2}\int \lc  (\De_{x}^{h}L^{ x}_{s-t})^{
2}-4h  L^{ x}_{s-t}\rc       \De_{x}^{h}\wt L^{ x}_{r} \,dx\)^{m}\rc p_{t}(y)\,dy\nn
\eea
and 
\begin{equation}
\wh {H}_{m, h}(t,s,r)=E\(h^{-2}    \int   \De_{x}^{h}L^{ x}_{t}  \,\(\De_{x}^{h}  L^{ x}_{s }\circ\th_{t}\) \, \De_{x}^{h}\wt L^{ x}_{r} \,dx    \)^{m}.\label{fu.6}
\end{equation}
We show in Lemma \ref{lem-3.1j} below  that for each integer $n\ge 0$,
\begin{eqnarray} &&
\lim_{ h\rar 0}E\(\({ \int   \De_{x}^{h}L^{ x}_{t}  \,\(\De_{x}^{h}  L^{ x}_{s }\circ\th_{t}\) \, \De_{x}^{h}\wt L^{ x}_{r} \,dx \over h^{ 2}}\)^{2n}\)\label{6.3lev}\\ 
&&\hspace{ .5in}  =  {( 2n)!\over 2^{ n}n!}\(   \displaystyle   64  \)^{ n} E\lc\(\int L^{ x}_{t} 
\(  L^{ x}_{s }\circ\th_{t}\)     \wt L^{ x}_{ r}\,dx\)^{
n}\rc,
\nn
\end{eqnarray}
locally uniformly in $r,s,t$ on $t>0$.
(\ref{fu.8}) then follows  by arguing as we did to obtain (\ref{fu.14}).

We can also write
\bea
&& \int \lc  (\De_{x}^{h}L^{ x}_{s})^{2}-4hL^{ x}_{s}      \rc \De_{x}^{h}\wt L^{ x}_{r} \,dx \label{fu.4}\\
&& = \int \lc  (\De_{x}^{h}L^{ x}_{s})^{2}-4hL^{ x}_{s}      \rc \De_{x}^{h}\wt L^{ x}_{q} \,dx     \nonumber\\
&& +\int \lc  (\De_{x}^{h}L^{ x}_{s})^{2}-4hL^{ x}_{s}      \rc \De_{x}^{h}(\wt L^{ x}_{r}-   \wt L^{ x}_{q}   )  \,dx    \nonumber \\
&& = \int \lc  (\De_{x}^{h}L^{ x}_{s})^{2}-4hL^{ x}_{s}      \rc \De_{x}^{h}\wt L^{ x}_{q} \,dx     \nonumber\\
&& +\int \lc  (\De_{x}^{h}L^{ x}_{s})^{2}-4hL^{ x}_{s}      \rc \De_{x}^{h}\wt L^{ x}_{r-q}\circ\wt \th_{q}   \,dx    \nonumber     
\eea
and hence as before this leads to 
\begin{eqnarray}
&& \lc{1 \over \de\de'}\int_{0}^{\de}\int_{q}^{q+ \de'}\wh {G}_{2n, h}(s,r)\,dr\,ds \rc^{1/(2n)} \geq   \lc{1 \over  \de} \int_{0}^{ \de} \wh {G}_{2n, h}(s,q)  \,ds\rc^{1/(2n)}\nonumber\\
&&\hspace{1 in}-
 \lc{1 \over \de\de'} \int_{0}^{\de}\int_{q}^{q+ \de'}\wt {G}_{2n, h}(s,r-q)\,dr\,ds\rc^{1/(2n)}
\label{fu.5a}
\end{eqnarray}
where
\bea
&&
\wt {G}_{m, h}(s,r-q)\label{fu.4e}\\
&&=:E\(h^{-2}\int \lc  (\De_{x}^{h}L^{ x}_{s })^{
2}-4h  L^{ x}_{s }       \rc \De_{x}^{h}\wt L^{ x}_{r-q}\circ\wt \th_{q} \,dx\)^{m}\nn\\
&&=\int E^{0,z}\lc \(h^{-2}\int \lc  (\De_{x}^{h}L^{ x}_{s })^{
2}-4h  L^{ x}_{s }\rc       \De_{x}^{h}\wt L^{ x}_{r-q} \,dx\)^{m}\rc p_{q}(z)\,dz\nn
\eea
and then (\ref{fu.9}) follows by arguing as we did to obtain   (\ref{fu.13}).

 Thus we obtain (\ref{77.15}) and hence  (\ref{7.53}) when $m$ is even.

\medskip	
In order to obtain (\ref{77.15}) when  $m$ is odd we first show that 
\begin{equation}
\sup_{h>0}\wh F_{2n,h}(t)\leq Ct^{2n }.\label{eq.1n}
\end{equation}
To this end, it clearly suffices to show that
\begin{equation}
\sup_{h>0}\wh F^{(0)}_{2n,h}(t)\leq Ct^{2n }\label{eq.1}
\end{equation}
where
\be 
\wh F^{(0)}_{m,h}(t):=  E\(\({ \int (\De_{x}^{h}L^{ x}_{t})^{
3}\,dx- 12h \int L^{ x}_{t}( \De_{x}^{h}L^{ x}_{t})\,dx\over h^{ 2}}\)^{m}\). \label{77.1en}
\ee
To see this we  observe that by first changing variables and then using  the scaling relationship (\ref{scale}) with $h=\sqrt{t}$,
we have 
\begin{eqnarray}
\int ( L^{ x+h}_{t}- L^{ x}_{ t})^{3}\,dx&=&\sqrt{t}\int ( L^{\sqrt{t}( x+ht^{-1/2})}_{t}- L^{\sqrt{t} x}_{ t})^{3}\,dx\label{eq.1a}\\
&=& t^{2}\int ( L^{  x+ht^{-1/2}}_{1}- L^{  x}_{ 1})^{3}\,dx\nonumber
\end{eqnarray}
and
\begin{eqnarray}
\int L^{ x}_{t}( L^{ x+h}_{t}- L^{ x}_{t})\,dx&=&\sqrt{t}\int L^{\sqrt{t} x}_{t}( L^{\sqrt{t}( x+ht^{-1/2})}_{t}- L^{\sqrt{t} x}_{ t}) \,dx\label{eq.1a}\\
&=& t^{3/2}\int L^{ x}_{1}( L^{ x+ht^{-1/2}}_{1}- L^{ x}_{1})\,dx.\nonumber
\end{eqnarray}
Therefore
\begin{eqnarray}
&&
{\int ( L^{ x+h}_{t}- L^{ x}_{ t})^{3}\,dx-12h \int L^{ x}_{t}( L^{ x+h}_{t}- L^{ x}_{t})\,dx
 \over h^{2}}  \label{eq.2}\\
 && \stackrel{\mathcal{L}}{=}
  {t^{ 2} \int ( L^{  x+ht^{-1/2}}_{1}- L^{  x}_{ 1})^{3}\,dx-12ht^{3/2}\int L^{ x}_{1}( L^{ x+ht^{-1/2}}_{1}- L^{ x}_{1})\,dx
 \over h^{ 2}} \nn \\
&&= t \,\,{\(\int ( L^{  x+ht^{-1/2}}_{1}- L^{  x}_{ 1})^{3}\,dx-12(ht^{-1/2}) \int L^{ x}_{1}( L^{ x+ht^{-1/2}}_{1}- L^{ x}_{1})\,dx     \) 
 \over (ht^{-1/2})^{ 2}}\nn
\end{eqnarray}
so that for any integer $m$
\begin{equation}
\wh F^{(0)}_{m,h}(t)=t^{ m }\wh F^{(0)}_{m,ht^{-1/2}}(1).\label{eq.3}
\end{equation}
Therefore to prove (\ref{eq.1}) we  need only show that
\begin{equation}
\sup_{t}\sup_{ h>0 }\wh F^{(0)}_{2n,ht^{-1/2}}(1)\leq C.\label{eq.4}
\end{equation}

It follows from (\ref{77.15}) that for some $\de>0$
\begin{equation}
\sup_{\{t,h\,|\,ht^{-1/2}\leq \de\}} \wh{F}^{(0)}_{2n,ht^{-1/2}}(1)\leq C.\label{eq.5}
\end{equation}
 On the other hand, for $ht^{-1/2}\geq \de$
\begin{eqnarray}
&&\Bigg|{\(\int ( L^{  x+ht^{-1/2}}_{1}- L^{  x}_{ 1})^{3}\,dx-12(ht^{-1/2}) \int L^{ x}_{1}( L^{ x+ht^{-1/2}}_{1}- L^{ x}_{1})\,dx     \) 
 \over (ht^{-1/2})^{ 2}}\Bigg|
\nn\\
&& \qquad\leq \de^{- 2}\int ( L^{  x+ht^{-1/2}}_{1}- L^{  x}_{ 1})^{3}\,dx +12\de^{-1}
|\int L^{ x}_{1}( L^{ x+ht^{-1/2}}_{1}- L^{ x}_{1})\,dx| \nonumber\\
&&\qquad \leq 4\de^{- 2}\int (L^{  x}_{ 1})^{3}\,dx +24\de^{-1}\int (L^{  x}_{ 1})^{2}\,dx \label{eq.6}
\end{eqnarray}
which has finite moments since 
  $\int (L^{  x}_{ 1})^{2}\,dx$ and $\int (L^{  x}_{ 1})^{3}\,dx$ have finite moments, see  (\ref{rb.1}). Using this  and (\ref{eq.5}) we get (\ref{eq.4}) and hence (\ref{eq.1}). As already noted, this implies (\ref{eq.1n}). It
then follows from   the Cauchy-Schwarz inequality that
\begin{equation}
\sup_{h>0}|\wh{F}_{m,h}(t)|\leq Ct^{ m }\label{eq.7}
\end{equation}
 for all integers $m$.
 
We next show that for any integer  $m$, the family of functions $\{\wh{F}_{m,h}(t);\,h\}$ is  equicontinuous in $t$, that is, for each $t$ and $\ep>0$ we can find a $\de>0$ such that  
\begin{equation}
\sup_{\{s\,|\,|s-t|\leq \de \}}\sup_{ h>0 }|\wh{F}_{m,h}(t)-\wh{F}_{m,h}(s)  |\leq \ep.\label{eq.8}
\end{equation}
Let
\begin{equation}
\Phi_{h}(t):={\int ( L^{ x+h}_{t}- L^{ x}_{ t})^{3}\,dx-12h \int L^{ x}_{t}( L^{ x+h}_{t}- L^{ x}_{t})\,dx
  -24h^{2}t\over h^{2}}.\label{eq.9}
\end{equation}
Applying the identity
$A^{m}-B^{m}=\sum_{j=0}^{m-1}A^{j}(A-B)B^{m-j-1}$  with $A=\Phi_{h}(t),\,B=\Phi_{h}(s)$ gives
\begin{equation}
   \wh{F}_{m,h}(t)-\wh{F}_{m,h}(s)=\sum_{j=0}^{m-1}\Phi_{h}(t)^{j}(\Phi_{h}(t)-\Phi_{h}(s))\Phi_{h}(s)^{m-j-1}
   \end{equation}
Consquently by using  the Cauchy--Schwarz inequality twice and  (\ref{eq.7}),  we see that
\begin{equation}
\sup_{\{s\,|\,|s-t|\leq \de \}}\sup_{ h>0 }|\wh{F}_{m,h}(t)-\wh{F}_{m,h}(s)  |\leq C_{t,m}\sup_{\{s\,|\,|s-t|\leq \de \}}\sup_{ h>0 }\| \Phi_{h}(t)-\Phi_{h}(s)  \|_{2}.\label{eq.10}
\end{equation}
Using (\ref{fu.3})-(\ref{ti.3}), we see that to obtain (\ref{eq.8}) it  suffices to show that for any $\ep>0$ we can find some  $\de>0$ such that 
\begin{equation}
\sup_{\{s\,|\,s\leq \de \}}\sup_{ h>0 }\wh{F}_{2,h}(s)  \leq \ep\label{eq.11}
\end{equation}
and  for any $T<\ff$
\begin{equation}
\sup_{\{\,t\leq T \}}\sup_{\{ \,s\leq \de \}}\sup_{ h>0 }\E\bigg[{1\over h^{2}}
\int \lc  (\De_{x}^{h}L^{ x}_{t})^{2}-4hL^{ x}_{t}     \rc\De_{x}^{h}\wt L^{ x}_{s} \,dx\bigg]^{2}\leq \ep\label{eq.12a}
\end{equation}
and
\begin{equation}
\sup_{\{\,s\leq T \}}\sup_{\{ \,t\leq \de \}}\sup_{ h>0 }\E\bigg[{1\over h^{2}}
\int \lc  (\De_{x}^{h}L^{ x}_{t})^{2}-4hL^{ x}_{t}     \rc\De_{x}^{h}\wt L^{ x}_{s} \,dx\bigg]^{2}\leq \ep.\label{eq.12b}
\end{equation}

By (\ref{eq.1})
\begin{equation}
\sup_{h>0}F_{2,h}(s)\leq Cs^{2}\label{eq.13}
\end{equation}
which immediately gives (\ref{eq.11}), while (\ref{eq.12a}) and (\ref{eq.12b}) follow from Lemma 
\ref{lem-var} below. This establishes (\ref{eq.8}).

 We now obtain (\ref{7.53}) when $m$ is odd. By equicontinuity, for any sequence $h_{n}\rar 0$, we can find a 
 subsequence $h_{n_{j}}\rar 0$,  such that 
 \begin{equation}
\lim_{j\rar\ff}\wh F_{m,h_{n_{j}}}(t)\label{eq.16}
 \end{equation}
converges   to a continuous function which we denote by $\ov F_{m}(t)$. It remains to show that
\begin{equation}
\ov F_{m}(t)\equiv 0.\label{eq.17}
\end{equation}
Let 
\begin{equation}
G_{m,h }(t):=e^{-t}\wh F_{m,h }(t) \hspace{.2 in}\mbox{and}\qquad\ov G_{m}(t):=e^{-t}\ov F_{m}(t).\label{eq.18}
\end{equation}
By (\ref{eq.7})
\begin{equation}
\sup_{h>0}\sup_{t}|G_{m,h}(t)|\leq C\quad\mbox{ and}\quad\lim_{t\rar \ff}\sup_{h>0}G_{m,h}(t)=0. \label{eq.19}
\end{equation}
 It then  follows from (\ref{77.4}) and the dominated convergence theorem that for all $\ze>0$
\begin{equation}
\int_{0}^{\ff}e^{- \ze s } \ov G_{m}(s) \,ds =0.\label{eq.20}
\end{equation}
We obtain  (\ref{eq.17}) by showing that $\ov G_{m}(s)\equiv 0$.

It   follows from (\ref{eq.19}) that  $\ov G_{m}(t) $ is a continuous bounded function on $R_{+}$ that goes to zero as $t\to\ff$. By the Stone--Weierstrass Theorem; (see \cite[Lemma 5.4]{K}), for any $\ep>0$, we can find a finite linear combination of the form $\sum_{i=1}^{n}c_{i}e^{- \ze_{i} s }$ such that
\begin{equation}
\sup_{s}|\ov G_{m}(s)-\sum_{i=1}^{n}c_{i}e^{- \ze_{i} s }|\leq \ep.\label{eq.21}
\end{equation}
Therefore, by  (\ref{eq.20})
\begin{eqnarray} 
\int_{0}^{\ff}e^{-  s } \ov G^{2}_{m}(s) \,ds\label{eq.22}&
 =&\int_{0}^{\ff}e^{-  s }\(\sum_{i=1}^{n}c_{i}e^{- \ze_{i} s }\) \ov G_{m}(s) \,ds\\
 &&\quad+\int_{0}^{\ff}e^{-  s } \(\ov G_{m}(s)-\sum_{i=1}^{n}c_{i}e^{- \ze_{i} s }\)\ov G_{m}(s) \,ds\nn\\
 &=&\int_{0}^{\ff}e^{-  s } \(\ov G_{m}(s)-\sum_{i=1}^{n}c_{i}e^{- \ze_{i} s }\)\ov G_{m}(s) \,ds\nn\\
 &\le &2\ep\(\int_{0}^{\ff}e^{-  s } \ov G^{2}_{m}(s) \,ds\)^{1/2}
\end{eqnarray}
by the Cauchy--Schwarz inequality and  (\ref{eq.21}).   Thus $\int_{0}^{\ff}e^{-  s } \ov G^{2}_{m}(s) \,ds=0$ which implies that $\ov G_{m}(s)\equiv 0$.

\qed

\section{Moments at exponential times}\label{sec-expmom}

We often write $u_{h,-h}^{\ze}(0)$ for $\De^{h}\De^{-h}u^{\ze}(0)=2\(u^{\ze}(0)-u^{\ze}(h)\)$.
\begin{lemma}\label{lem-2weak}
For each $m$, as $h\rar 0$
\begin{eqnarray} &&
E\(\({ \int \lc (\De_{x}^{h}L^{ x}_{ \la_{\ze}})^{
3} - 6u_{h,-h}^{\ze}(0)   L^{ x}_{\la_{\ze}}\De_{x}^{h}L^{ x}_{ \la_{\ze}} 
- 6\(u_{h,-h}^{\ze}\)^{2}  L^{ x}_{\la_{\ze}}\rc\,dx\over h^{ 2}}\)^{m}\)\nn\\ 
&&\hspace{ .8in} \Longrightarrow \left\{\begin{array}{ll}
{( 2n)!\over 2^{ n}n!}\( 192\)^{ n} E\lc\(\int (L^{ x}_{ \la_{\ze}})^{ 3}\,dx\)^{
n}\rc &\mbox{ if }m=2n\\
0&\mbox{ otherwise.}
\end{array}
\right.
\label{7.54}
\end{eqnarray}
\end{lemma}

{\bf  Remark: } Of course, $\int  L^{ x}_{ \la_{\ze}} \,dx=\la_{\ze}$. Using (\ref{1.8}) and the continuity of local time, Lemma \ref{lem-2weak} implies (\ref{7.54a}).

{\bf  Proof of Lemma \ref{lem-2weak}}: 
For any
integer
$m$ we have
\bea &&  E\(\( \int \lc (\De_{x}^{h}L^{ x}_{ \la_{\ze}})^{
3} - 6u_{h,-h}^{\ze}(0)   L^{ x}_{\la_{\ze}}\De_{x}^{h}L^{ x}_{ \la_{\ze}} 
- 6\(u_{h,-h}^{\ze}\)^{2}  L^{ x}_{\la_{\ze}}\rc\,dx\)^{ m}
\)\nn\\ &&=E\(\prod_{ i=1}^{ m}\(  \int \lc (\De_{x_{ i}}^{h}L^{ x_{ i}}_{ \la_{\ze}})^{
3} - 6u_{h,-h}^{\ze}(0)   L^{ x_{ i}}_{\la_{\ze}}\De_{x_{ i}}^{h}L^{ x_{ i}}_{ \la_{\ze}} 
- 6\(u_{h,-h}^{\ze}\)^{2}  L^{ x_{ i}}_{\la_{\ze}}\rc\,dx_{ i}\)
\)\nn\\ &&=\sum_{A,B,C} ( -1)^{m- |B|- |C|}(6u_{h,-h}^{\ze}(0))^{|B|}\(6\(u_{h,-h}^{\ze}\)^{2}\)^{|C|}\label{1.16g}\\ &&
\hspace{ .1in}E\(
\(\prod_{ i\in A}  \int (\De_{x_{ i}}^{h}L^{ x_{ i}}_{ \la_{\ze}})^{
3}\,dx_{
i}\)\(\prod_{ j\in B}\int L^{ x_{ j}}_{
\la_{\ze}}\De_{x_{ j}}^{h}L^{ x_{ j}}_{ \la_{\ze}}\,dx_{ j}\)\(\prod_{ k\in B}\int L^{ x_{ k}}_{
\la_{\ze}}\,dx_{ k}\) \),\nn
\eea
where the sum runs over all partitions of $[1,m]$ into three parts, $A,B,C$.

We initially calculate 
\begin{equation} E\(\prod_{ i\in A} \De_{x_{ i}}^{h}L^{ x_{ i}}_{ \la_{\ze}}\De_{y_{ i}}^{h}L^{ y_{ i}}_{ \la_{\ze}}\De_{z_{ i}}^{h}L^{ z_{ i}}_{ \la_{\ze}}
\prod_{j\in B}L^{ x_{ j}}_{\la_{\ze}}\De_{y_{ j}}^{h}L^{ y_{ j}}_{ \la_{\ze}}\prod_{ k\in C}L^{ x_{ k}}_{\la_{\ze}}\)\label{1.18g}
\end{equation}
 and eventually we set $y_{l}=x_{ l}=z_{l}$ for all $l$. 
 Using (\ref{1.2w}) we have 
\begin{eqnarray} && E\(\prod_{ i\in A} \De_{x_{ i}}^{h}L^{ x_{ i}}_{ \la_{\ze}}\De_{y_{ i}}^{h}L^{ y_{ i}}_{ \la_{\ze}}\De_{z_{ i}}^{h}L^{ z_{ i}}_{ \la_{\ze}}
\prod_{j\in B}L^{ x_{ j}}_{\la_{\ze}}\De_{y_{ j}}^{h}L^{ y_{ j}}_{ \la_{\ze}}\prod_{ k\in C}L^{ x_{ k}}_{\la_{\ze}}\)\label{1.19g}\\ && =\(\prod_{ i\in A}\De_{ x_{ i}}^{ h}\De_{ y_{ i}}^{
h}\De_{ z_{ i}}^{ h}\prod_{ j\in B} \De_{y_{ j}}^{h} \)E\(\prod_{ i\in A}L^{ x_{ i}}_{
\la_{\ze}} L^{ y_{ i}}_{ \la_{\ze}}L^{ z_{ i}}_{ \la_{\ze}}\prod_{ j\in B} L^{ x_{ j}}_{\la_{\ze}}L^{ y_{ j}}_{ \la_{\ze}}\prod_{ k\in C}L^{ x_{ k}}_{\la_{\ze}}\)
\nonumber\\ && =\(\prod_{ i\in A}\De_{ x_{ i}}^{ h}\De_{ y_{ i}}^{
h}\De_{ z_{ i}}^{ h}\prod_{ j\in B} \De_{y_{ j}}^{h} \)\,\,
\sum_{ \si}\prod_{ j=1}^{m+2|A|+|B|}u^{\ze}( \si( j)-\si( j-1)) \nonumber
\end{eqnarray} where the sum runs over all
bijections
\bas
&&\hspace{-.5 in}
\si:\,[1,2,\ldots,m+2|A|+|B|]\mapsto \\
&&\hspace{1 in}
\{ x_{ i}, y_{ i},z_{i},  i\in A\}\cup \{ x_{ j},y_{ j}, j\in B\}\cup \{ x_{ k},  k\in C\}.
\eas
We then use the product
rule 
\begin{equation}\De_{ x}^{ h}\{ f( x)g( x)\}=
\{ \De_{ x}^{ h} f( x)\}g( x+h)+f( x)\{  \De_{ x}^{ h}g( x)\}\label{pr1}
\end{equation}
to expand the right hand side of (\ref{1.19g}) into a sum of many terms, over all
$\si$ and all ways to allocate each $\De_{ x_{ i}}^{ h},\De_{ y_{ i}}^{ h}, \De_{ z_{ i}}^{ h}$ or $\De_{ y_{ k}}^{ h}$ to
a single $u^{\ze}$ factor.

Consider first the case where $A=\{ 1,\ldots,m\}$. For a given term in the above
expansion, we will say that $x_{ i}$ is $3$-bound if $x_{ i},y_{ i},z_{i} $ are adjacent, (in other words,  for some $j$ we have $ (\si(j),\si(j+1),\si(j+2))=(x_{ i},y_{ i},z_{i})$ or one of its $6$ possible permutations), and  $\De_{ x_{ i}}^{ h},\De_{ y_{ i}}^{ h},\De_{ z_{ i}}^{ h}$   are all attached    to the  $u^{\ze}$ factors which connect $x_{ i},y_{ i},z_{i}$. Thus if $ (\si(j),\si(j+1),\si(j+2))=(x_{ i},y_{ i},z_{i})$, the  $\De_{ x_{ i}}^{ h},\De_{ y_{ i}}^{ h},\De_{ z_{ i}}^{ h}$   are all attached to $u^{\ze}( y_{ i}-x_{ i})u^{\ze}( z_{ i}-y_{ i})$. We return shortly to analyze this case.

If $x_{ i}$ is not $3$-bound, we will say that
it is $2$-bound if   any two of the three elements $x_{ i},y_{ i},z_{i} $ are adjacent, for example if $x_{ i},y_{ i} $ are adjacent, (in other words either $(x_{ i},y_{ i})=(\si(j),\si(j+1))$ or $(y_{ i},x_{ i})=(\si(j),\si(j+1))$ for some $j$), and both $\De_{ x_{ i}}^{ h}$ and
$\De_{ y_{ i}}^{ h}$ are attached to the factor $u^{\ze}( x_{ i}-y_{ i})$.   In applying (\ref{pr1}) we are free to choose which function plays the role of $f$, and which the role of $g$. In case $x_{ i}$ is $2$-bound, taking our example with $(x_{ i},y_{ i})=(\si(j),\si(j+1))$, when using  (\ref{pr1})  to expand  $\De_{ x_{ i}}^{ h}$  we take  $g$ to be $u^{\ze}( x_{ i}-y_{ i})$ and
similarly when using  (\ref{pr1})  to expand  $\De_{ y_{ i}}^{ h}$. In this way we guarantee that we have not added   $\pm h$ to the arguments of any other factors.
By (\ref{1.8}),  setting $ x_{ i}=y_{ i}$   turns the factor $\De_{ x_{ i}}^{ h}\De_{ y_{ i}}^{ h}u^{\ze}( x_{ i}-y_{ i})$ into $\De^{h}\De^{-h}u^{\ze}(0)$, and since for every such $\si$ there is precisely one other bijection which differs from  $\si$ only in that it 
permutes 
$x_{ i},y_{ i} $ we obtain a factor of $2\De^{h}\De^{-h}u^{\ze}(0)$. This is precisely what we
would have obtained if instead of $\De_{ x_{ i}}^{ h}\De_{ y_{ i}}^{ h}L^{ x_{
i}}_{\la_{\ze}} L^{ y_{ i}}_{ \la_{\ze}}$ in (\ref{1.19g}) we had $2\De^{h}\De^{-h}u^{\ze}(0)L^{ x_{ i}}_{\la_{\ze}} $. There are ${3 \choose 2 }=3$ ways to pick two letters from among 
$\{x_{i},y_{i},z_{i}\}$. By considering all such cases,  we obtain precisely what we
would have obtained if instead of $\De_{ x_{ i}}^{ h}\De_{ y_{ i}}^{ h}\De_{ z_{ i}}^{ h}L^{ x_{
i}}_{\la_{\ze}} L^{ y_{ i}}_{ \la_{\ze}}L^{ z_{ i}}_{ \la_{\ze}}$ in (\ref{1.19g}) we had $6\De^{h}\De^{-h}u^{\ze}(0)L^{ x_{ i}}_{\la_{\ze}} \De_{ z_{ i}}^{ h}L^{ z_{
i}}_{\la_{\ze}}$. Consider then a term 
in the expansion of  (\ref{1.19g}) with $A=\{ 1,\ldots,m\}$ and $J=\{i\,|\,x_{ i} \mbox{ is $2$-bound}\} $  non-empty, but $ \{i\,|\,x_{ i} \mbox{ is $3$-bound}\}= \emptyset$.
By  (\ref{1.19g}) there will be an identical contribution from the last line of  (\ref{1.16g}) from any other   $A, B, C$ with $B\subseteq J$ and $C=\emptyset$.
Since 
$\sum_{B\subseteq J} ( -1)^{ |B|}=0$, we see that in the
expansion of (\ref{1.16g}) there will not be any contributions from $2$-bound $x$'s.

We emphasize that if $x_{ i}$ is $2$-bound, and, for example, $ (\si(j),\si(j+1),\si(j+2))=(x_{ i},y_{ i},z_{i})$, with both   $\De_{ x_{ i}}^{ h},\De_{ y_{ i}}^{ h}$   attached to $u^{\ze}( y_{ i}-x_{ i})$,   then   $\De_{ z_{ i}}^{ h}$ cannot be assigned to    $u^{\ze}( z_{ i}-y_{ i})$. Otherwise, $x_{ i}$ would be $3$-bound.

We now return to analyze the case where $x_{ i}$ is $3$-bound. Consider the case that $ (\si(j),\si(j+1),\si(j+2))=(x_{ i},y_{ i},z_{i})$. We first apply   the $\De_{ x_{ i}}^{ h}$ and $\De_{ z_{ i}}^{ h}$ operators  to 
$u^{\ze}( y_{ i}-x_{ i})u^{\ze}( z_{ i}-y_{ i})$ to obtain  $\De_{ x_{ i}}^{ h}u^{\ze}( y_{ i}-x_{ i})\,\,\De_{ z_{ i}}^{ h}u^{\ze}( z_{ i}-y_{ i})$.
Then by (\ref{pr1}) we have 
\begin{eqnarray}
&& \De_{ y_{ i}}^{ h}\(\De_{ x_{ i}}^{ h}u^{\ze}( y_{ i}-x_{ i})\,\,\De_{ z_{ i}}^{ h}u^{\ze}( z_{ i}-y_{ i})\)
\label{nu.1}\\
&& =  \(\De_{ y_{ i}}^{ h}\De_{ x_{ i}}^{ h}u^{\ze}( y_{ i}-x_{ i})\)\,\,\De_{ z_{ i}}^{ h}u^{\ze}( z_{ i}-y_{ i}-h)\nonumber\\
&& + \De_{ x_{ i}}^{ h}u^{\ze}( y_{ i}-x_{ i})\,\,\(\De_{ y_{ i}}^{ h}\De_{ z_{ i}}^{ h}u^{\ze}( z_{ i}-y_{ i})\).\nonumber
\end{eqnarray}
If we now set $y_{ i}=x_{ i}=z_{i}$ we obtain
\begin{equation}
\De^{h}\De^{-h}u^{\ze}(0)\(u^{\ze}(0)-u^{\ze}(h)\)+\(u^{\ze}(h)-u^{\ze}(0)\)\De^{h}\De^{-h}u^{\ze}(0)=0.\label{nu.2}
\end{equation}
Thus,  $3$-bound variables such as $x_{ i}$ make no contribution to (\ref{1.19g}). However, there will be an analogous contribution from $6\De^{h}\De^{-h}u^{\ze}(0)L^{ x_{ i}}_{\la_{\ze}} \De_{ z_{ i}}^{ h}L^{ z_{i}}_{\la_{\ze}}$ which is not cancelled by $2$-bound variables. This is the case where 
$x_{ i},z_{ i} $ are adjacent,  say $(x_{ i},z_{ i})=(\si(j),\si(j+1))$, and  $\De_{ z_{ i}}^{ h}$ is attached to the factor $u^{\ze}( x_{ i}-z_{ i})$. As before we may do this without adding an $h$ to the arguments in any other factors. After setting $x_{ i}=z_{ i} $ we obtain $6\De^{h}\De^{-h}u^{\ze}(0)\(u^{\ze}(h)-u^{\ze}(0)\)=-3\(\De^{h}\De^{-h}u^{\ze}(0)\)^{2}$. Since we can also interchange 
$x_{ i},z_{ i}$, such $x_{i}$ contribute $-6\(\De^{h}\De^{-h}u^{\ze}(0)\)^{2}$, which will be exactly canceled by the term $-6\(\De^{h}\De^{-h}u^{\ze}(0)\)^{2}L^{ x_{i}}_{\la_{\ze}}$.
Furthermore, this completely exhausts the contribution to (\ref{1.16g}) of all  
$A\neq \{ 1,\ldots,m\}$. 

Thus, in estimating (\ref{1.19g}) we need only consider $A=\{ 1,\ldots,m\}$ and 
those cases where each of the $3m$ $\De^{ h}$'s are assigned either to unique
$u^{\ze}$ factors, or if two are assigned to the same $u^{\ze}$ factor, it is not of the form $u^{\ze}(
x_{ i}-y_{ i}), u^{\ze}(
x_{ i}-z_{ i})$ or $u^{\ze}(
z_{ i}-y_{ i})$. 

For ease of exposition, in the following  calculations we first replace the right hand side of (\ref{pr1}) by $\{ \De_{ x}^{ h} f( x)\}g( x)+f( x)\{  \De_{ x}^{ h}g( x)\}$, and return at the end of the proof to explain why this doesn't affect the final result. We use the notation
\begin{equation}
E_{\ast}\(\(\int \lc (\De_{x}^{h}L^{ x}_{ \la_{\ze}})^{
3} - 6u_{h,-h}^{\ze}(0)   L^{ x}_{\la_{\ze}}\De_{x}^{h}L^{ x}_{ \la_{\ze}} 
- 6\(u_{h,-h}^{\ze}\)^{2}  L^{ x}_{\la_{\ze}}\rc\,dx\)^{ m}
\)\label{pr.1}
\end{equation}
to denote the expression obtained with this replacement.
We can thus write
\bea && E_{\ast}\(\( \int \lc (\De_{x}^{h}L^{ x}_{ \la_{\ze}})^{
3} - 6u_{h,-h}^{\ze}(0)   L^{ x}_{\la_{\ze}}\De_{x}^{h}L^{ x}_{ \la_{\ze}} 
- 6\(u_{h,-h}^{\ze}\)^{2}  L^{ x}_{\la_{\ze}}\rc\,dx\)^{ m}
\)\nn\\ &&=6^{ m}\sum_{ \pi,a}\int \mathcal{T}_{h}( x;\,\pi,a)\,dx\label{1.20g}
\eea with
\begin{equation}
\mathcal{T}_{h}( x;\,\pi,a) =\prod_{ j=1}^{ 3m}\(\De^{ h}_{ x_{ \pi( j)}}\)^{a_{ 1}(j)}
\(\De^{ h}_{ x_{ \pi( j-1)}}\)^{a_{ 2}(j)}\,u^{\ze}( x_{\pi(j)}- x_{\pi(j-1)})\label{1.21g}
\end{equation} where the sum runs over all  maps $\pi\,:\,[1,\ldots, 3m]\mapsto
[1,\ldots, m]$ with $|\pi^{ -1}(i )|=3$ for each $i$, and all `assignments' 
$a=(a_{ 1},a_{ 2})\,:\,[1,\ldots, 2m]\mapsto \{ 0,1\}\times \{ 0,1\}$ with the
property that for each $i$ there will be exactly three operators of the form $\De^{
h}_{ x_{i}}$ in (\ref{1.21g}), and if $a( j)=( 1,1)$ for any $j$, then 
$ x_{\pi(j)}\neq x_{\pi(j-1)}$. The factor $6^{ m}=(3!)^{m}$ in (\ref{1.20g})  comes from 
  the fact that $|\pi^{ -1}(i )|=3$ for
each $i$.

  Let  $m=2n$. Assume first that  $a=e$ where now $e( 2j)=( 1,1)$, $e( 2j-1)=( 0,0)$ for all $j$.

 \subsection{$a =e$ with $\pi$ compatible with a pairing}\label{ss-3.1t}
 
When $a =e$ we have
 \begin{equation}\qquad
\mathcal{T}_{h}( x;\,\pi,e) =\prod_{ j=1}^{3n}u^{\ze}( x_{\pi(2j-1)}- x_{\pi(2j-2)})\,\De^{ h}\De^{- h}\,u^{\ze}( x_{\pi(2j)}- x_{\pi(2j-1)}).\label{91.1}
\end{equation}

Let $\mathcal{P}=\{(l_{2i-1},l_{2i})\,,\,1\leq i\leq n\}$ be a pairing of the integers $[1,2n]$. Let $\pi$   as in (\ref{91.1}) be     such that for each $1\leq j\leq 3n$, 
$\{\pi(2j-1), \pi(2j)\}=\{l_{2i-1},l_{2i}\}$ for some, necessarily unique, $ 1\leq i\leq n$. In this case we say that $\pi$ is compatible with the pairing $\mathcal{P}$ and write this  as $ \pi \sim \mathcal{P}$. (Note that  when we write $\{\pi(2j-1), \pi(2j)\}=\{l_{2i-1},l_{2i}\}$ we mean as two  sets, so, according to what $\pi$ is, we may have  $\pi(2j-1)=l_{2i-1}$, $\pi(2j )=l_{2i }$ or $\pi(2j-1)=l_{2i}$, $\pi(2j )=l_{2i-1 }$.   )   
In this case we have 
 \begin{equation} \qquad
\mathcal{T}_{h}( x;\,\pi,e) =\prod_{ i=1}^{ n}\(\De^{ h}\De^{- h}\,u^{\ze}( x_{l_{2i}}- x_{l_{2i-1}})\)^{3}\prod_{ j=1}^{ 3n}u^{\ze}( x_{\pi(2j-1)}- x_{\pi(2j-2)})\,.\label{91.2}
\end{equation}

We now show that 
\begin{equation}
\int \mathcal{T}_{h}( x;\,\pi,e)\prod_{j=1}^{2n}\,dx_{j} =\int \mathcal{T}_{1,h}( x;\,\pi,a)\prod_{j=1}^{2n}\,dx_{j}+O(h^{4n+1})\label{91.3}
\end{equation}
where
\begin{eqnarray}\qquad
\mathcal{T}_{1,h}( x;\,\pi,e)&=&\prod_{ i=1}^{ n}\(1_{\{|x_{l_{2i}}-x_{l_{2i-1}}|\leq h\}}\)\(\De^{ h}\De^{- h}\,u^{\ze}( x_{l_{2i}}- x_{l_{2i-1}})\)^{3}
\nn\\
&& \hspace{.8 in}  \times \prod_{ j=1}^{ 3n}u^{\ze}( x_{\pi(2j-1)}- x_{\pi(2j-2)}).\label{91.4}
\end{eqnarray}
To prove (\ref{91.3}) we write
\begin{equation}
1=\prod_{ i=1}^{ n}\( 1_{\{|x_{l_{2i}}-x_{l_{2i-1}}|\leq h\}}+1_{\{|x_{l_{2i}}-x_{l_{2i-1}}|\geq h\}}\),\label{}
\end{equation}
insert this inside the integral on the left hand side of (\ref{91.3}) and expand the product. It then suffices to show that
\begin{eqnarray} 
&& \int \prod_{i\in A} 1_{\{|x_{l_{2i}}-x_{l_{2i-1}}|\leq h\}}\prod_{i\in A^{c}} 1_{\{|x_{l_{2i}}-x_{l_{2i-1}}|\geq h\}} \(w^{\ze}( x_{l_{2i}}- x_{l_{2i-1}})\)^{3}
\nn\\
&&   \hspace{.5 in} \times \prod_{ j=1}^{ 3n}u^{\ze}( x_{\pi(2j-1)}- x_{\pi(2j-2)})\prod_{j=1}^{2n}\,dx_{j}=O(h^{4n+1})\label{91.4}
\end{eqnarray}
whenever $|A^{c}|\geq 1$. To see this we first  choose $j_{k},\,k=1,\ldots,n$ so that 
\[ \{x_{\pi(2j_{k}-1)}- x_{\pi(2j_{k}-2)},\,k=1,\ldots,n \}\cup \{x_{l_{2i}}- x_{l_{2i-1}},\,i=1,\ldots,n \}\] generate $\{x_{1},\ldots, x_{2n}\}$.
After changing variables, (\ref{91.4}) follows  easily from (\ref{li.13}), (\ref{1.30gb}) and  the fact that $u^{\ze}$ is bounded and integrable. 

We then study  
\begin{equation} \hspace{.4 in}
\int \mathcal{T}_{1,h}( x;\,\pi,e)\prod_{j=1}^{2n}\,dx_{j}. \label{f91.36}
\end{equation}
 Recall that for  each $1\le j\le 3n$, 
$\{\pi(2j-1), \pi(2j)\}=\{l_{2i-1},l_{2i}\} $, for some $1\le i\le n$.  We identify these relationships by setting $i=\si (j) $ when  $\{\pi(2j-1), \pi(2j)\}=\{l_{2i-1},l_{2i}\} $.  In the present situation   this means that $\si:\,[1,3n]\mapsto [1,n]$  with $|\si^{-1}(i)|=3$ for each $1\leq i\leq n$. (One for each occurrence of $\{l_{2i-1},l_{2i}\} $). We write
\begin{eqnarray}
&&\prod_{ j=1}^{3 n}\, u^{\ze}(x_{\pi(2j-1)}-x_{\pi(2j-2)})
\label{f91.37}\\
&&\qquad=\prod_{ j=1}^{3n}\,\( u^{\ze}(x_{l_{2\si (j)-1}}-x_{l_{2\si (j-1)-1}})+\De^{h_{j}}u^{\ze}(x_{l_{2\si (j)-1}}-x_{l_{2\si (j-1)-1}})\) ,   \nn
\end{eqnarray}
where $h_{j}=(x_{\pi(2j-1)}-x_{l_{2\si (j)-1}})+(x_{l_{2\si (j-1)-1}}-x_{\pi(2j-2)})$. 
Note that because of the presence of  the term $\prod_{i=1}^{n}\(1_{\{|x_{l_{2i}}-x_{l_{2i-1}}|\leq h\}}\)$ in the integral in (\ref{f91.36})  we need only be concerned with values of $|h_{j}|\leq 2h$, $1\le j\le 3n$.

We expand the product on the right hand side of (\ref{f91.37}) and obtain a sum of many terms. Using (\ref{1.3x})  and the fact that $|h_{j}|\leq 2h$, $1\le j\le 3n$ we can see as above that 
\begin{eqnarray}
&&\int \mathcal{T}_{1,h}( x;\,\pi,e)\prod_{j=1}^{2n}\,dx_{j}
\label{91.6a}\\
&&=\int \prod_{ i=1}^{ n}\(1_{\{|x_{l_{2i}}-x_{l_{2i-1}}|\leq h\}}\)  \prod_{ i=1}^{ n}\(\De^{ h}\De^{- h}\,u^{\ze}( x_{l_{2i}}- x_{l_{2i-1}})\)^{3} \nonumber\\
&&\hspace{1 in}\prod_{ j=1}^{ 3n}\,  u^{\ze}(x_{l_{2\si (j)-1}}-x_{l_{2\si (j-1)-1}})\prod_{j=1}^{2n}\,dx_{j}+O(h^{4n+1})\nonumber
\end{eqnarray}
where $x_{-1}=0$. Once again we can now see that
\begin{eqnarray}
&&\int \mathcal{T}_{1,h}( x;\,\pi,e)\prod_{j=1}^{2n}\,dx_{j}
\label{91.6}\\
&&=\int  \prod_{ i=1}^{ n}\(\De^{ h}\De^{- h}\,u^{\ze}( x_{l_{2i}}- x_{l_{2i-1}})\)^{3} \nonumber\\
&&\hspace{1 in}\prod_{ j=1}^{ 3n}\,  u^{\ze}(x_{l_{2\si (j)-1}}-x_{l_{2\si (j-1)-1}})\prod_{j=1}^{2n}\,dx_{j}+O(h^{4n+1}).\nonumber
\end{eqnarray}

Using translation invariance and then (\ref{1.30g}) we have
\begin{eqnarray}
&&\int   \prod_{ i=1}^{ n}\(\De^{ h}\De^{- h}\,u^{\ze}( x_{l_{2i}}- x_{l_{2i-1}})\)^{3} \prod_{ j=1}^{ 3n}\,  u^{\ze}(x_{l_{2\si (j)-1}}-x_{l_{2\si (j-1)-1}})\prod_{j=1}^{2n}\,dx_{j}
\nn\\
&& =\int   \prod_{ i=1}^{ n}\(\De^{ h}\De^{- h}\,u^{\ze}( x_{l_{2i}})\)^{3} \prod_{ j=1}^{ 3n}\,  u^{\ze}(x_{l_{2\si (j)-1}}-x_{l_{2\si (j-1)-1}})\prod_{k=1}^{2n}\,dx_{l_{k}}  \nonumber\\
&& =(4+O( h))^{n}h^{ 4n}\int     \prod_{ j=1}^{ 3n}\,  u^{\ze}(x_{l_{2\si (j)-1}}-x_{l_{2\si (j-1)-1}})\prod_{k=1}^{n}\,dx_{l_{2k-1}}.  \label{91.7}
\end{eqnarray}
Rewriting this and summarizing, we have shown that 
\begin{eqnarray}
&&\int \mathcal{T}_{h}( x;\,\pi,e)\prod_{j=1}^{2n}\,dx_{j}
\label{91.8}\\
&& =4^{n}h^{ 4n}\int     \prod_{ j=1}^{ 3n}\,  u^{\ze}(y_{\si (j)}-y_{\si (j-1)})\prod_{k=1}^{n}\,dy_{k} +O(h^{4n+1}) \nonumber
\end{eqnarray}
with $y_{0}=0$.
 
Let $\mathcal{M}$ denote the set of maps $\si$   from $[1,\ldots,3n]$ to $[1,\ldots,n]$ such that $|\si^{ -1}( i)|=3$ for all $i$.
For each pairing $\mathcal{P}$, any  
$\pi\sim \mathcal{P}$ gives rise as above to a map $\si\in \mathcal{M}$. Also, any of the $2^{ 3n}$ $\pi$'s obtained by permuting the
$2$ elements in each of the $3n$ pairs, give rise to the same $\si$. In addition, for any  $\si'\in \mathcal{M}$, we can   reorder  the $3n$ pairs of $\pi$ to obtain a new   $\pi'\sim \mathcal{P}$ which gives rise to $\si'$. Thus we have shown that 
\begin{eqnarray}  &&\sum_{ \pi\sim  \mathcal{P}}\int\mathcal{T}_{h}(
x;\,\pi,e)\prod_{ j=1}^{ 2n}\,dx_{j}\label{1.35g}\\ &&=  \(  32 h^{ 4}\)^{ n} \sum_{ \si\in  \mathcal{M}}\int     \prod_{ j=1}^{ 3n}\,  u^{\ze}(y_{\si (j)}-y_{\si (j-1)})\prod_{k=1}^{n}\,dy_{k}+O(h^{4n+1})\nn\\ &&=  \( {16\over 3}h^{4}\)^{ n} E\lc\(\int (L^{ x}_{ \la_{\ze}})^{ 3}\,dx\)^{
n}\rc+O(h^{4n+1})\nn
\end{eqnarray}
where the last line follows from Kac's moment formula, compare (\ref{1.20g}).
To complete this subsection, let $\mathcal{G} $ denote the set of $\pi$ which are compatible with some   pairing 
$\mathcal{P}$. Since  there are ${( 2n)! \over
2^{ n}n!}$  pairings of $[1,\ldots,2n]$, we have shown that
\begin{eqnarray}  &&\sum_{ \pi \in \mathcal{G} }\int\mathcal{T}_{h}(
x;\,\pi,e)\prod_{ j=1}^{ 2n}\,dx_{j}\label{1.35g}\\ &&= {( 2n)! \over
2^{ n}n!}\( {16\over 3}h^{ 4}\)^{ n} E\lc\(\int (L^{ x}_{ \la_{\ze}})^{ 3}\,dx\)^{
n}\rc+O(h^{4n+1}).\nn
\end{eqnarray}

 \subsection{$a=e $ but   $\pi$ not compatible with a pairing}\label{ss-3.2t}
 
 In this subsection we show that  
\begin{equation}
\sum_{ \pi \not \in \mathcal{G} }\Big |\int\mathcal{T}_{h}(
x;\,\pi,e)\prod_{ j=1}^{ 2n}\,dx_{j}\Big |=O(h^{4n+1}).\label{91.10}
\end{equation}

We return to (\ref{91.1}) to obtain 
 \bea
 &&
\Big |\int\mathcal{T}_{h}(
x;\,\pi,e)\prod_{ j=1}^{ 2n}\,dx_{j}\Big |\label{91.12}\\
&&\leq \int \prod_{ j=1}^{3n}u^{\ze}( x_{\pi(2j-1)}- x_{\pi(2j-2)})\,w^{\ze}( x_{\pi(2j)}- x_{\pi(2j-1)})\prod_{ j=1}^{ 2n}\,dx_{j}
.\nn
\eea

Let us show that when  $\pi$ not compatible with a pairing we can find $n+1$ linearly independent vectors from among the $3n$ vectors
\begin{equation}
x_{\pi(2j)}- x_{\pi(2j-1)}, \hspace{.3 in}1\leq j\leq 3n.\label{v.1}
\end{equation}
We will say that $x$ and $y$ are both `contained' in $x-y$. Since $|\pi^{-1}(i)|=3$ for each 
$1\leq i\leq 2n$, we can find $j_{1}$ such that $x_{1}$ is contained in $x_{\pi(2j_{1})}- x_{\pi(2j_{1}-1)}$. In addition,   $x_{\pi(2j_{1})}- x_{\pi(2j_{1}-1)}$ will contain another $x_{i_{1}}$. 
We then pick an integer from  $[2,\ldots, 2n]-\{i_{1}\}$, say $i_{2}$ and then find $x_{\pi(2j_{2})}- x_{\pi(2j_{2}-1)}$ which contains $x_{i_{2}}$. $x_{\pi(2j_{2})}- x_{\pi(2j_{2}-1)}$ will contain another 
$x_{i_{3}}$ where we may have $i_{3}=1$ or $i_{3}=i_{1}$. In any event it is clear that in this manner we can obtain a sequence of vectors
\begin{equation}
x_{\pi(2j_{i})}- x_{\pi(2j_{i}-1)}, \hspace{.3 in}1\leq i\leq n\label{v.2}
\end{equation}
which are linearly independent, since for each $i$,  $x_{\pi(2j_{i})}- x_{\pi(2j_{i}-1)}$ contains some $x_{k}$ not contained in any of the preceding $x_{\pi(2j_{l})}- x_{\pi(2j_{l}-1)}, 1\leq l<i $. 
Then, if the $n$ vectors in (\ref{v.2}) contain all $x_{i}, 1\leq i\leq 2n$, it follows that  the $n$ vectors in (\ref{v.2}) must contain disjoint pairs of $x_{i}$'s. As a consequence they cannot generate any vector of the form $x_{i}-x_{i'}$ which is not among the $n$ vectors in (\ref{v.2}). Since by our assumption   $\pi$ is not compatible with a pairing, there are vectors of the form $x_{\pi(2j)}- x_{\pi(2j-1)}$
which are different from the $n$ vectors in (\ref{v.2}). This proves our claim that we can find $n+1$ linearly independent vectors from among the $3n$ vectors of (\ref{v.1}) in case the $n$ vectors in (\ref{v.2}) contain all $x_{i}, 1\leq i\leq 2n$. But if they do not contain all $x_{i}, 1\leq i\leq 2n$, say they do not contain $x_{k}$. There is some vector in (\ref{v.1}) which contains $x_{k}$, and it is clearly linearly independent of  the vectors in (\ref{v.2}).

Thus we have a sequence 
\begin{equation}
x_{\pi(2j_{i})}- x_{\pi(2j_{i}-1)}, \hspace{.3 in}1\leq i\leq n+1\label{v.3}
\end{equation}
  of linearly independent vectors. Let $J=\{j_{i},\,1\leq i\leq n+1\}$. We use (\ref{1.3x}) to bound (\ref{91.12}) by 
   \bea
 &&
\Big |\int\mathcal{T}_{h}(
x;\,\pi,e)\prod_{ j=1}^{ 2n}\,dx_{j}\Big |\label{v.4}\\
&&\leq Ch^{2n-1}\int \prod_{ j=1}^{3n}u^{\ze}( x_{\pi(2j-1)}- x_{\pi(2j-2)})\,\prod_{j\in J^{c}} u^{\ze}( x_{\pi(2j)}- x_{\pi(2j-1)})  \nonumber\\
&&\hspace{1 in}\,\prod_{j\in J}w^{\ze}( x_{\pi(2j)}- x_{\pi(2j-1)})\prod_{ j=1}^{ 2n}\,dx_{j}
.\nn
\eea
We can complete the set of $n+1$ vectors in (\ref{v.3}) to a basis of $x_{i}, 1\leq i\leq 2n$   by choosing $n-1$  vectors from among the vectors appearing as arguments of $u^{\ze}$ in the second line of (\ref{v.4}). We then bound the remaining $u^{\ze}$ factors by a constant, change variables and use (\ref{li.13})  with $q=1$ to see that the integral on the right hand side of (\ref{v.4}) is bounded by $Ch^{2(n+1)}$. Combining this with (\ref{v.4}) proves (\ref{91.10}).

  \subsection{$a\neq e $}\label{ss-3.3t}
  
We now claim that
\begin{equation}
\sum_{\pi }\sum_{  a \neq e}\Big |\int \mathcal{T}_{h}( x;\,\pi ,a ) \,\prod_{ j=1}^{ 2n}\,dx_{j}\Big |=O(h^{4n+1}). \label{91.20}
\end{equation} 

If $\mathcal{T}_{h}( x;\,\pi ,a )$ contains $k<3n$ factors of the form $w^{\ze}$, then we will have 
$2(3n-k)$ factors of the form $v^{\ze}$. We use (\ref{1.3x}) then to bound 
\begin{equation}
\Big |\int \mathcal{T}_{h}( x;\,\pi ,a ) \,\prod_{ j=1}^{ 2n}\,dx_{j}\Big |\leq C h^{2(3n-k)}
\int \mathcal{I}_{h}( x;\,\pi ,a ) \,\prod_{ j=1}^{ 2n}\,dx_{j}. \label{91.20s}
\end{equation} 
where $\mathcal{I}_{h}( x;\,\pi ,a ) $ is   similar to $\mathcal{T}_{h}( x;\,\pi ,a )$ except that we have bounded the integrand by its absolute value and replaced each $v^{\ze}$ by $u^{\ze}$. We now show how to get a good bound on the integral of $\mathcal{I}_{h}( x;\,\pi ,a ) $.

Choose $j_{1}=1,j_{2},\ldots,j_{2n}$ such that 
\begin{eqnarray}
&&\mbox{span }\{x_{\pi( j_{1})} , x_{\pi(j_{2})}- x_{\pi(j_{2}-1)},\ldots, x_{\pi(j_{2n})}- x_{\pi(j_{2n}-1)}\}
\nn\\
&&=\mbox{span }\{x_{1},\ldots, x_{2n}\}   \label{sp.3}
\end{eqnarray}
It is easy to see that this can be done.
We now show that we can choose a permutation   $\si_{1},\si_{2},\ldots, \si_{2n}$ of 
$[1,2n]$ such that   for any $1\leq k\leq 2n$ 
 \begin{eqnarray} 
 &&x_{\pi(j_{k})},x_{\pi(j_{k}-1)}\in \{x_{\si_{1}} , x_{\si_{2}},\ldots, x_{\si_{k}}\}.   \label{sp.2}
\end{eqnarray} 

 We take $\si_{1}=\pi( j_{1})=\pi( 1)$ and   choose the $\si_{2},\ldots, \si_{2n}$ by induction so that (\ref{sp.2}) holds for each $1\leq k\leq 2n$. This clearly holds for $k=1$, since by definition $x_{\pi(0)}=0$. Assume we have chosen $\si_{1}, \si_{2},\ldots, \si_{l}$  so that (\ref{sp.2}) holds for all  $k\leq l$.   Then among the remaining $\{  x_{\pi(j_{l+1})}- x_{\pi(j_{l+1}-1)},\ldots, x_{\pi(j_{2n})}- x_{\pi(j_{2n}-1)}\}$ there will be at least one $i$ such  
that either $x_{\pi(j_{i})}$ or $x_{\pi(j_{i}-1)}$ is equal to one of $x_{\si_{1}} , x_{\si_{2}},\ldots, x_{\si_{l}}$. This is because each element of $\{  x_{\pi(j_{l+1})}- x_{\pi(j_{l+1}-1)},\ldots, x_{\pi(j_{2n})}- x_{\pi(j_{2n}-1)}\}$ is a difference of $x$'s so that by themselves we could never
have 
\begin{eqnarray}
&&\mbox{span }\{  x_{\pi(j_{l+1})}- x_{\pi(j_{l+1}-1)},\ldots, x_{\pi(j_{2n})}- x_{\pi(j_{2n}-1)}\}
\nn\\
&&=\mbox{span }\( \{x_{1},\ldots, x_{2n}\}-\{x_{\si_{1}} , x_{\si_{2}},\ldots, x_{\si_{l}}\}\).   \label{sp.3}
\end{eqnarray}
We then take such an $i$ and if $x_{\pi(j_{i+1})}$  is equal to one of $x_{\si_{1}} , x_{\si_{2}},\ldots, x_{\si_{l}}$ set $\si_{l+1}=\pi(j_{i+1}-1)$, while if   $x_{\pi(j_{i+1}-1)}$ is equal to one of $x_{\si_{1}} , x_{\si_{2}},\ldots, x_{\si_{l}}$ set $\si_{l+1}=\pi(j_{i+1})$. Then (\ref{sp.2}) holds for   $k= l+1$ and completes our induction.

We will prove (\ref{91.10}) by first bounding the $dx_{ \si_{2n}}$ integral in (\ref{91.12}) involving all factors containing 
 $x_{ \si_{2n}}$. We then bound the $dx_{ \si_{2n-1}}$ integral involving all \underline{remaining} factors containing 
 $x_{\si_{2n-1}}$. We then iterate this procedure bounding in turn the $dx_{\si_{2n}}, dx_{\si_{2n-1}},\ldots,dx_{\si_{1}}$ integrals. (\ref{sp.2}) guarantees that at each stage we are integrating a non-empty product of bounded integrable functions. Note that by (\ref{1.3y}) and  (\ref{li.13}) with $q=1$
\bea
\sup_{a_{i}}\int \prod_{i=1}^{p}w^{\ze}(y+a_{i})\,dy  \leq   Ch^{p-1}\sup_{a_{1}} 
 \int  w^{\ze} (y+a_{1})\,dy =O(h^{p+1}) \label{com.13x}
 \eea
for all $p\geq 1$. 

Let $p_{i}$ denote the number of \underline{remaining} $w^{\ze}$ factors containing 
 $x_{\si_{i}}$ after we have bounded in turn the $dx_{\si_{2n}}, dx_{\si_{2n-1}},\ldots,dx_{\si_{i+1}}$ integrals, that is,
 $p_{i}$ denotes the number of  $w^{\ze}$ factors containing 
 $x_{\si_{i}}$ but not any of $x_{\si_{2n}}, x_{\si_{2n-1}},\ldots,x_{\si_{i+1}}$. Since there are a total of $k<3n$ factors of the form $w^{\ze}$ in (\ref{91.20s}), we have that $\sum_{i=1}^{2n}p_{i}=k$. Let  $k_{0}=|\{i\,|\,p_{i}\neq 0\}|$. 
 
 If we apply our bounding procedure using (\ref{com.13x}) together with the fact that 
 $u^{\ze}$ is bounded and integrable we see that  
 \begin{equation}
\int \mathcal{I}_{h}( x;\,\pi ,a )  \,\prod_{ j=1}^{ 2n}\,dx_{j}=O(h^{\sum_{i=1}^{2n}(p_{i}+1_{\{p_{i}\neq 0\}})}) 
=O(h^{k+k_{0}}).\label{com.14x}
\end{equation}
 It is easy to see that in (\ref{91.20s})   each $x_{j}$ appears in at most $3$ factors of the form $w^{\ze}$. Thus each $p_{i}\leq 3$. Since  $\sum_{i=1}^{2n}p_{i}=k$, we must have   $k_{0}\geq k/3$.   Combining (\ref{com.14x}) with  (\ref{91.20s}) we have 
\begin{equation}
\Big |\int \mathcal{T}_{h}( x;\,\pi ,a ) \,\prod_{ j=1}^{ 2n}\,dx_{j}\Big |\leq C h^{2(3n-k)+4k/3}
=Ch^{4n+2(n-k/3)} \label{91.20t}
\end{equation} 
 which proves (\ref{91.20}) since $k<3n$.

Combining (\ref{91.20}), (note the factor $6^{m}=36^{n}$) with the results of Subsections \ref{ss-3.1t}-\ref{ss-3.2t} we have thus shown that
\begin{eqnarray} 
&&E_{\ast}\(\( \int \lc (\De_{x}^{h}L^{ x}_{ \la_{\ze}})^{
3} - 6u_{h,-h}^{\ze}(0)   L^{ x}_{\la_{\ze}}\De_{x}^{h}L^{ x}_{ \la_{\ze}} 
- 6\(u_{h,-h}^{\ze}\)^{2}  L^{ x}_{\la_{\ze}}\rc\,dx\)^{ 2n}
\)\nn\\ &&= {( 2n)! \over 2^{ n}n!}\( 192 h^{ 4}\)^{
n} E\lc\(\int (L^{ x}_{ \la_{\ze}})^{ 3}\,dx\)^{ n}\rc+O(h^{ 4n+1} ).\label{1.39g}
\end{eqnarray}

We now explain why we obtain the same expression on the right hand side when we have $E$ instead of $E_{\ast}$.  In the paragraph following (\ref{pr.1})  we make use of precise cancellations to handle bound variables,  which a priori might be affected by our modification of (\ref{pr1}). 
Consider how (\ref{1.20g}) would look if we had used the 
product formula (\ref{pr1}). Note that any estimates we used will still apply  since these estimates involve integrating or bounding by the supremum, neither of which are affected by replacing any of the $x$'s are replaced by 
$x+h$. (\ref{91.6}) will be affected, but note from (\ref{pr1}) that the only  terms of the form $u^{\ze}(x-y)$ that may have $x$ replaced by $x\pm h$ are those   to which $\De^{h}_{x}$ is not applied.  Similarly $y$ may be replaced by $y\pm h$ only if $\De^{h}_{y}$ is not applied to a term of the form $u^{\ze}(x-y)$ .   Consequently,     we still have all terms of the form $\De^{h}\De^{-h}u^{\ze}$.   Thus we obtain (\ref{91.8}), except that some of the remaining $u^{\ze}(x-y)$ may be replaced by  $u^{\ze}(x-y\pm h)$. Using (\ref{1.3x}) then leads to (\ref{91.8}).

  Thus, it only remains to show that for each
$n$
\begin{eqnarray}  &&E\(\(\int \lc (\De_{x}^{h}L^{ x}_{ \la_{\ze}})^{
3} - 6u_{h,-h}^{\ze}(0)   L^{ x}_{\la_{\ze}}\De_{x}^{h}L^{ x}_{ \la_{\ze}} 
- 6\(u_{h,-h}^{\ze}\)^{2}  L^{ x}_{\la_{\ze}}\rc\,dx\)^{ 2n+1}
\)\nn\\ &&\hspace{ 3in}=O(h^{ 2( 2n+1) +1} ).\label{1.40g}
\end{eqnarray} This follows from the fact that we cannot form any $\pi$'s   in
$\mathcal{G}$.
\qed

\section{Proof of Lemma \ref{lem-6.3a}}\label{sec-6.1}

We use $E^{y,z}(\cdot)$ to denote expectation with respect to the independent Brownian motions $B_{t}$ starting at $y$ and $\wt B_{t}$ starting at $z$.

\bl\label{lem-6.3a}
Let $\la_{\ze}, \la_{\ze'}$ be independent  exponential random variables with mean $1/\ze, 1/\ze'$.
For each integer $m\ge 0$,
\begin{eqnarray} &&
\lim_{ h\rar 0}E^{y,z}\(\({ \int \lc  (\De_{x}^{h}L^{ x}_{\la_{\ze}})^{
2}-2\De^{h}\De^{-h}u^{\ze}(0)  L^{ x}_{\la_{\ze}}       \rc\De_{x}^{h}\wt L^{ x}_{\la_{\ze'}}\,dx \over h^{ 2}}\)^{m}\)\label{16.30}\\ 
&&\hspace{ .5in}  =\left\{\begin{array}{ll}
  \displaystyle  {( 2n)!\over 2^{ n}n!}\(   \displaystyle   64 \)^{ n} E^{y,z}\lc\(\int (L^{ x}_{ \la_{\ze}})^{ 2}\wt L^{ x}_{ \la_{\ze'}}\,dx\)^{
n}\rc &\mbox{ if }m=2n\\\\
0&\mbox{ otherwise}
\end{array}
\right.
\nn
\end{eqnarray}
uniformly in $y,z$.
  \el

  {\bf  Proof of Lemma \ref{lem-6.3a}}: 
For any
integer
$m$ we have
\bea && \qquad E^{y,z}\(\(  \int \lc  (\De_{x}^{h}L^{ x}_{\la_{\ze}})^{
2}-2\De^{h}\De^{-h}u^{\ze}(0)  L^{ x}_{\la_{\ze}}       \rc\De_{x}^{h}\wt L^{ x}_{\la_{\ze'}}\,dx\)^{ m}
\)\label{16.31}\\ &&=E^{y,z}\(\prod_{ i=1}^{ m}\(  \int \int \lc  (\De_{x_{i}}^{h}L^{ x_{i}}_{\la_{\ze}})^{
2}-2\De^{h}\De^{-h}u^{\ze}(0)  L^{ x_{i}}_{\la_{\ze}}       \rc\De_{x_{i}}^{h}\wt L^{ x_{i}}_{\la_{\ze'}}\,dx\,dx_{ i}\)
\)\nn\\ &&=\sum_{A\subseteq \{ 1,\ldots,m\}} ( -1)^{m- |A|}(2\De^{h}\De^{-h}u^{\ze}(0))^{|A^{c}|}\nn\\ &&
\hspace{ .7in}E^{y,z}\(
\(\prod_{ i\in A}  \int (\De_{x_{ i}}^{h}L^{ x_{ i}}_{ \la_{\ze}})^{
2}\De_{x_{ i}}^{h}\wt L^{ x_{ i}}_{ \la_{\ze}}\,dx_{
i}\)\(\prod_{ k\in A^{ c}}\int L^{ x_{ k}}_{
\la_{\ze}}\De_{x_{k}}^{h}\wt L^{ x_{ k}}_{ \la_{\ze}}\,dx_{k}\) \).\nn
\eea

We initially calculate 
\begin{equation}(2\De^{h}\De^{-h}u^{\ze}(0))^{|A^{c}|} E^{y,z}\(\prod_{ i\in A} \De_{x_{ i}}^{h}L^{ x_{ i}}_{ \la_{\ze}}\De_{y_{ i}}^{h}L^{ y_{ i}}_{ \la_{\ze}}\De_{z_{ i}}^{h}\wt L^{ z_{ i}}_{ \la_{\ze}}
\prod_{ k\in A^{ c}}L^{ u_{ k}}_{\la_{\ze}}\De_{v_{ k}}^{h}\wt L^{ v_{ k}}_{ \la_{\ze}}\)\label{16.32}
\end{equation}
 and eventually we set $y_{ j}=x_{ j}=z_{j}$ and $u_{j}=v_{ j}$ for all $j$.
  Using (\ref{1.2w}) we have 
\begin{eqnarray} && E^{y,z}\(\prod_{ i\in A} \De_{x_{ i}}^{h}L^{ x_{ i}}_{ \la_{\ze}}\De_{y_{ i}}^{h}L^{ y_{ i}}_{ \la_{\ze}}\De_{z_{ i}}^{h}\wt L^{ z_{ i}}_{ \la_{\ze}}
\prod_{ k\in A^{ c}}L^{ u_{ k}}_{\la_{\ze}}\De_{v_{ k}}^{h}\wt L^{ v_{ k}}_{ \la_{\ze}}\)\label{16.33}\\ && =\(\prod_{ i\in A}\De_{ x_{ i}}^{ h}\De_{ y_{ i}}^{
h}\De_{ z_{ i}}^{ h}\prod_{ k\in A^{ c}} \De_{v_{ k}}^{h} \)E^{y,z}\(\prod_{ i\in A}L^{ x_{ i}}_{
\la_{\ze}} L^{ y_{ i}}_{ \la_{\ze}}\wt L^{ z_{ i}}_{ \la_{\ze}}\prod_{ k\in A^{ c}} L^{ u_{ k}}_{\la_{\ze}}\wt L^{ v_{ k}}_{ \la_{\ze}}\)
\nonumber\\ && =\(\prod_{ i\in A}\De_{ x_{ i}}^{ h}\De_{ y_{ i}}^{
h}\De_{ z_{ i}}^{ h}\prod_{ k\in A^{ c}} \De_{v_{ k}}^{h} \)\nonumber\\
&&\hspace{1 in}E^{y}\(\prod_{ i\in A}L^{ x_{ i}}_{
\la_{\ze}} L^{ y_{ i}}_{ \la_{\ze}}\prod_{ k\in A^{ c}} L^{ u_{ k}}_{\la_{\ze}}\)E^{z}\(\prod_{ i\in A}\wt L^{ z_{ i}}_{ \la_{\ze}}\prod_{ k\in A^{ c}} \wt L^{ v_{ k}}_{ \la_{\ze}}\)
\nonumber\\ && =\(\prod_{ i\in A}\De_{ x_{ i}}^{ h}\De_{ y_{ i}}^{
h} \)\,\,
\sum_{ \si }\prod_{ j=1}^{m+|A|}u^{\ze}( \si( j)-\si( j-1)) \nonumber\\ && \hspace{.5 in}
\(\prod_{ i\in A}\De_{ z_{ i}}^{ h}\prod_{ k\in A^{ c}} \De_{v_{ k}}^{h} \)\,\,
\sum_{ \si'}\prod_{ j=1}^{m}u^{\ze}( \si'( j)-\si'( j-1)) \nonumber
\end{eqnarray} 
where the first sum runs over all
bijections\[\si:\,[1,\ldots,m+|A|]\mapsto 
\{ x_{ i}, y_{ i},   i\in A\}\cup \{ u_{ i},  i\in A^{ c} \}.\] 
with $\si(0)=y$ and the second sum runs over all
bijections\[\si':\,[1,\ldots,m]\mapsto 
\{  z_{i},  i\in A\}\cup \{ v_{ i},  i\in A^{ c} \}\]
with $\si'(0)=z$. 

We then use the product
to expand the right hand side of (\ref{16.33}) into a sum of many terms, over all
$\si$ and all ways to allocate each $\De_{ x_{ i}}^{ h},\De_{ y_{ i}}^{ h}$ or $\De_{ z_{ i}}^{ h}$ to
a single $v$ factor. 

Consider first the case where $A=\{ 1,\ldots,m\}$. For a given term in the above
expansion, we will say that
$x_{ i}$ is bound if $x_{ i},y_{ i} $ are adjacent, (in other words either $(x_{ i},y_{ i})=(\si(j),\si(j+1))$ or $(y_{ i},x_{ i})=(\si(j),\si(j+1))$ for some $j$), and both $\De_{ x_{ i}}^{ h}$ and
$\De_{ y_{ i}}^{ h}$ are attached to the factor $u^{\ze}( x_{ i}-y_{ i})$.  
Setting $ x_{ i}=y_{ i}$   turns the factor $\De_{ x_{ i}}^{ h}\De_{ y_{ i}}^{ h}u^{\ze}( x_{ i}-y_{ i})$ into $\De^{h}\De^{-h}u^{\ze}(0)$, and as in the proof of Lemma \ref{lem-2weak} we can assume that our use of (\ref{pr1}) for $\De_{ x_{ i}}^{ h}$ and
$\De_{ y_{ i}}^{ h}$ does not introduce a $\pm h$ in the arguments of other factors.
For every such $\si$ there is precisely one other $\si$ which agrees with $\si$ except that it 
permutes 
$x_{ i},y_{ i} $, we obtain a factor of $2\De^{h}\De^{-h}u^{\ze}(0)$. This is precisely what we
would have obtained if instead of $\De_{ x_{ i}}^{ h}\De_{ y_{ i}}^{ h}L^{ x_{
i}}_{\la_{\ze}} L^{ y_{ i}}_{ \la_{\ze}}$ in (\ref{16.33}) we had $2\De^{h}\De^{-h}u^{\ze}(0)L^{ u_{ i}}_{\la_{\ze}} $.  Consider then any term 
in the expansion of  (\ref{16.33}) with $A=\{ 1,\ldots,m\}$ and $J=\{i\,|\,x_{ i} \mbox{ is bound}\}$.
By  (\ref{1.19g}) there will be an identical contribution from the last line of  (\ref{16.33}) for any other  $A$ with $A^{c}\subseteq J$.
Since 
$\sum_{A^{c}\subseteq J} ( -1)^{ |A^{c}|}=0$, we see that in the
expansion of (\ref{16.31}) there will not be any contributions from bound $x$'s.
Furthermore, this completely exhausts the contribution to (\ref{16.31}) of all  
$A\neq \{ 1,\ldots,m\}$. 

Thus in estimating (\ref{16.33}) we need only consider $A=\{ 1,\ldots,m\}$ and 
those cases where  if two $\De^{ h}$'s are assigned to the same $u^{\ze}$ factor, it is not of the form $u^{\ze}(
x_{ i}-y_{ i})$. 

Again as in the proof of Lemma \ref{lem-2weak}, in the following calculation we first replace the right hand side of (\ref{pr1}) by $\{ \De_{ x}^{ h} f( x)\}g( x)+f( x)\{  \De_{ x}^{ h}g( x)\}$, and return at the end of the proof to explain why this doesn't affect the final result. We use the notation
\begin{equation}
E^{y,z}_{\ast}\(\( \int  \lc  (\De_{x}^{h}L^{ x}_{\la_{\ze}})^{
2}-2\De^{h}\De^{-h}u^{\ze}(0)  L^{ x}_{\la_{\ze}}       \rc\De_{x}^{h}\wt L^{ x}_{\la_{\ze'}}\,dx\)^{ m}
\)\label{pr.16}
\end{equation}
to denote the expression obtained with this replacement.

We can thus write
\bea &&\hspace{ .1in} E^{y,z}_{\ast}\(\( \int  \lc  (\De_{x}^{h}L^{ x}_{\la_{\ze}})^{
2}-2\De^{h}\De^{-h}u^{\ze}(0)  L^{ x}_{\la_{\ze}}       \rc\De_{x}^{h}\wt L^{ x}_{\la_{\ze'}}\,dx\)^{ m}
\)\label{16.34}\\ &&=2^{ m}\sum_{ \pi ,a }\,\sum_{  \pi', a'}\int \mathcal{T}_{h}( x;\,\pi,\pi',a,a')\,dx\nn
\eea 
with
\bea
&&
\mathcal{T}_{h}( x;\, \pi,\pi',a,a') \nonumber\\
&&=\prod_{ j=1}^{ 2m}\(\De^{ h}_{ x_{ \pi( j)}}\)^{a_{ 1}(j)}
\(\De^{ h}_{ x_{ \pi( j-1)}}\)^{a_{ 2}(j)}\,u^{\ze}( x_{\pi(j)}- x_{\pi(j-1)})\label{16.21g}\\
& &\times\prod_{ j=1}^{ m}\(\De^{ h}_{ x_{ \pi( j)}}\)^{a'_{ 1}(j)}
\(\De^{ h}_{ x_{ \pi( j-1)}}\)^{a'_{ 2}(j)}\,u^{\ze'}( x_{\pi'(j)}- x_{\pi'(j-1)})\nn
\eea
where the first  sum runs over all  maps $\pi\,:\,[1,\ldots, 2m]\mapsto
[1,\ldots, m]$ with $|\pi^{ -1}(i )|=2$ for each $i$, and all `assignments' 
$a=(a_{ 1},a_{ 2})\,:\,[1,\ldots, 2m]\mapsto \{ 0,1\}\times \{ 0,1\}$ with the
property that for each $i$ there will be exactly two factors of the form $\De^{
h}_{ x_{i}}$ in the second line of  (\ref{16.21g}), and if $a( j)=( 1,1)$ for any $j$, then 
$ x_{\pi(j)}\neq x_{\pi(j-1)}$. The factor $2^{ m}$ in (\ref{16.34})  comes from 
  the fact that $|\pi^{ -1}(i )|=2$ for
each $i$. Similarly, the second sum runs over all  permutations $\pi'\,:\,[1,\ldots, m]\mapsto
[1,\ldots, m]$, and all `assignments' 
$a'=(a'_{ 1},a'_{ 2})\,:\,[1,\ldots, 2m]\mapsto \{ 0,1\}\times \{ 0,1\}$  with the
property that for each $i$ there will be exactly one factor of the form $\De^{
h}_{ x_{i}}$ in  in the last line of  (\ref{16.21g}). Here we have set $x_{ \pi(0)}=y, x_{ \pi'(0)}=z.$

From this point on the proof is very similar to that of Lemma \ref{lem-2weak}. Let  $m=2n$. Assume first that  $a=e$ where now $e( 2j)=( 1,1)$, $e( 2j-1)=( 0,0)$ for all $j$, and similarly for $a'$.

 \subsection{$a =a'=e$ with $\pi, \pi'$ compatible with a pairing}\label{ss-6.1t}
 
When $a  =a'=e$ we have
 \bea
 &&
\mathcal{T}_{h}( x;\, \pi,\pi',e,e) =\prod_{ j=1}^{2n}u^{\ze}( x_{\pi(2j-1)}- x_{\pi(2j-2)})\,\De^{ h}\De^{- h}\,u^{\ze}( x_{\pi(2j)}- x_{\pi(2j-1)})\nonumber\\
&&  \times \prod_{ j=1}^{n}u^{\ze'}( x_{\pi'(2j-1)}- x_{\pi'(2j-2)})\,\De^{ h}\De^{- h}\,u^{\ze'}( x_{\pi'(2j)}- x_{\pi'(2j-1)})  .\label{691.1}
\eea

Let $\mathcal{P}=\{(l_{2i-1},l_{2i})\,,\,1\leq i\leq n\}$ be a pairing of the integers $[1,2n]$. Let $\pi$   as in (\ref{691.1}) be     such that for each $1\leq j\leq 2n$, 
$\{\pi(2j-1), \pi(2j)\}=\{l_{2i-1},l_{2i}\}$ for some, necessarily unique, $ 1\leq i\leq n$. In this case we say that $\pi$ is compatible with the pairing $\mathcal{P}$ and write this  as $ \pi \sim \mathcal{P}$.   Similarly we say that  $\pi'$ is compatible with the pairing $\mathcal{P}$ if for  each $1\leq j\leq n$, 
$\{\pi(2j-1), \pi(2j)\}=\{l_{2i-1},l_{2i}\}$ for some, necessarily unique, $ 1\leq i\leq n$, and write this  as $ \pi' \sim \mathcal{P}$.  
If  $\pi,  \pi' \sim \mathcal{P}$ we have 
 \bea
 &&
\mathcal{T}_{h}( x;\, \pi,\pi',e,e)\nonumber\\
&& =\prod_{ i=1}^{ n}\(\De^{ h}\De^{- h}\,u^{\ze}( x_{l_{2i}}- x_{l_{2i-1}})\)^{2}\prod_{ j=1}^{ 2n}u^{\ze}( x_{\pi(2j-1)}- x_{\pi(2j-2)})\,\nonumber\\
&&\hspace{.1 in}  \times \prod_{ i=1}^{ n} \De^{ h}\De^{- h}\,u^{\ze'}( x_{l_{2i}}- x_{l_{2i-1}}) \prod_{ j=1}^{ n}u^{\ze'}( x_{\pi'(2j-1)}- x_{\pi'(2j-2)})     .\label{691.2}
\eea

Set $\si (j)=i $ when  $\{\pi(2j-1), \pi(2j)\}=\{l_{2i-1},l_{2i}\} $, so that $\si:\,[1,2n]\mapsto [1,n]$  with $|\si^{-1}(i)|=2$ for each $1\leq i\leq n$. Similarly, set $\si' (j)=i $ when  $\{\pi'(2j-1), \pi'(2j)\}=\{l_{2i-1},l_{2i}\} $, so that that $\si'$ is a permutation of $ [1,n] $.

As in Sub-section \ref{ss-3.1t} we can  show that
\begin{eqnarray}
&&\int \mathcal{T}_{h}( x;\, \pi,\pi',e,e)\prod_{j=1}^{2n}\,dx_{j}
\label{691.8}\\
&& =4^{n}h^{ 4n}\int     \prod_{ j=1}^{ 2n}\,  u^{\ze}(y_{\si (j)}-y_{\si (j-1)})
\prod_{ j=1}^{ n}\,  u^{\ze'}(y_{\si' (j)}-y_{\si' (j-1)})\prod_{k=1}^{n}\,dy_{k} +O(h^{4n+1}) \nonumber
\end{eqnarray}
with   $y_{ \si(0)}=y, y_{ \si'(0)}=z$ and error term uniform in $y,z$.
 
Let $\mathcal{M}_{d}$ denote the set of maps $\si$   from $[1,\ldots,dn]$ to $[1,\ldots,n]$ such that $|\si^{ -1}( i)|=d$ for all $i$.
For each pairing $\mathcal{P}$,  each  map
$\pi\sim \mathcal{P}$ gives rise as above to a map $\si\in \mathcal{M}_{2}$. Also, any of the $2^{ 2n}$ $\pi$'s obtained by permuting the
$2$ elements in each of the $2n$ pairs, give rise to the same $\si$. In addition, for any  $\wh\si\in \mathcal{M}_{2}$, we can   reorder  the $2n$ pairs of $\pi$ to obtain a new   $\wh \pi\sim \mathcal{P}$ which gives rise to $\wh \si$. A similar analysis applies to our $\pi'$.

Thus we have shown that 
\begin{eqnarray}  &&\sum_{ \pi,\pi'\sim  \mathcal{P}}\int\mathcal{T}_{h}( x;\, \pi,\pi',e,e)\prod_{ j=1}^{ 2n}\,dx_{j}\label{61.35g}\\ &&=  \(  32 h^{ 4}\)^{ n} \sum_{ \si\in  \mathcal{M}_{2},\,\si'\in \mathcal{M}_{1}}\int     \prod_{ j=1}^{ 2n}\,  u^{\ze}(y_{\si (j)}-y_{\si (j-1)})\prod_{k=1}^{n}\,dy_{k}+O(h^{4n+1})\nn\\ &&=  \(  16 h^{4}\)^{ n} E^{y,z}\lc\(\int (L^{ x}_{ \la_{\ze}})^{ 2}L^{ x}_{ \la_{\ze'}}\,dx\)^{
n}\rc+O(h^{4n+1})\nn
\end{eqnarray}
where the last line follows from Kac's moment formula.
To complete this subsection, let $\mathcal{G} $ denote the set of $\pi,\pi'$ which are compatible with some   pairing 
$\mathcal{P}$. Since  there are ${( 2n)! \over
2^{ n}n!}$  pairings of $[1,\ldots,2n]$, we have shown that
\begin{eqnarray}  &&\sum_{ \pi,\pi' \in \mathcal{G} }\int\mathcal{T}_{h}( x;\, \pi,\pi',e,e)\prod_{ j=1}^{ 2n}\,dx_{j}\label{61.35g}\\ &&= {( 2n)! \over
2^{ n}n!}\(  16 h^{ 4}\)^{ n} E^{y,z}\lc\(\int (L^{ x}_{ \la_{\ze}})^{ 2}L^{ x}_{ \la_{\ze'}}\,dx\)^{
n}\rc+O(h^{4n+1}).\nn
\end{eqnarray}

The rest of the proof follows as in the proof of Lemma \ref{lem-2weak}.
  \qed

  \section{Proof of Lemma \ref{lem-3.1j}}\label{sec-1.1}

  \bl\label{lem-3.1j} For each integer $n\ge 0$,
\begin{eqnarray} &&
\lim_{ h\rar 0}E\(\({ \int   \De_{x}^{h}L^{ x}_{t_{1}}  \,\(\De_{x}^{h}  L^{ x}_{t_{2} }\circ\th_{t_{1}}\) \, \De_{x}^{h}\wt L^{ x}_{t_{3}} \,dx \over h^{ 2}}\)^{2n}\)\label{6.3lev}\\ 
&&\hspace{ .5in}  =  {( 2n)!\over 2^{ n}n!}\(   \displaystyle   64  \)^{ n} E\lc\(\int L^{ x}_{t_{1}} 
\(  L^{ x}_{t_{2} }\circ\th_{t_{1}}\)     \wt L^{ x}_{ t_{3}}\,dx\)^{
n}\rc,
\nn
\end{eqnarray}
locally uniformly in $t_{1},t_{2},t_{3}$ on $t_{1}>0$.
\el
  
{\bf Proof of Lemma \ref{lem-3.1j} }	The proof of this Lemma is easier than that of Lemmas  
\ref{lem-2weak} and \ref{lem-6.3a} since there are no subtraction terms. However, there are complications due to the fact that we now work with non-random times $t_{1},t_{2},t_{3}$ rather than exponential times. 

We begin by writing
\begin{eqnarray}
&&E\(\( \int   \De_{x}^{h}L^{ x}_{t_{1}}  \,\(\De_{x}^{h}  L^{ x}_{t_{2} }\circ\th_{t_{1}}\) \, \De_{x}^{h}\wt L^{ x}_{t_{3}} \,dx  \)^{2n}\)
\label{semi.1}\\
&&= E\(\prod_{i=1}^{2n}\( \int   \De_{x_{i}}^{h}L^{ x_{i}}_{t_{1}}  \,\(\De_{x_{i}}^{h}  L^{ x_{i}}_{t_{2} }\circ\th_{t_{1}}\) \, \De_{x_{i}}^{h}\wt L^{ x_{i}}_{t_{3}} \,dx_{i}  \) \)  \nonumber
\\
&&=\int\, E\(\prod_{i=1}^{2n}\(   \De_{x_{i}}^{h}L^{ x_{i}}_{t_{1}}  \,\(\De_{x_{i}}^{h}  L^{ x_{i}}_{t_{2} }\circ\th_{t_{1}}\) \, \De_{x_{i}}^{h}\wt L^{ x_{i}}_{t_{3}}   \) \)\prod_{i=1}^{2n} \,dx_{i} \nonumber
\end{eqnarray}
We first evaluate
\begin{eqnarray}
&&E\(\prod_{i=1}^{2n}\(     \De_{x_{i}}^{h}L^{ x_{i}}_{t_{1}}  \,\(\De_{y_{i}}^{h}  L^{ y_{i}}_{t_{2} }\circ\th_{t_{1}}\) \, \De_{z_{i}}^{h}\wt L^{ z_{i}}_{t_{3}}   \) \)
\label{semi.2}\\
&& =\prod_{i=1}^{2n}\( \De_{x_{i}}^{h} \De_{y_{i}}^{h} \De_{z_{i}}^{h}\) E\(\prod_{i=1}^{2n}\(    L^{ x_{i}}_{t_{1}}  \,\(   L^{ y_{i}}_{t_{2} }\circ\th_{t_{1}}\) \,  \wt L^{ z_{i}}_{t_{3}}   \) \) \nonumber\\
&& =\prod_{i=1}^{2n}\( \De_{x_{i}}^{h} \De_{y_{i}}^{h}  \) E\(\prod_{i=1}^{2n}     L^{ x_{i}}_{t_{1}}  \,\(   L^{ y_{i}}_{t_{2} }\circ\th_{t_{1}}\)     \) \nonumber\\
&& \hspace{1.5 in} \prod_{i=1}^{2n}\(\De_{z_{i}}^{h}\) E\(\prod_{i=1}^{2n}      \wt L^{ z_{i}}_{t_{3}}     \) \nonumber
\end{eqnarray}
and then set all $y_{i}=z_{i}=x_{i}$.

By Kac's moment formula
\begin{eqnarray}
&& E\(\prod_{i=1}^{2n}     L^{ x_{i}}_{t_{1}}  \,\(   L^{ y_{i}}_{t_{2} }\circ\th_{t_{1}}\)     \)
\label{semi.3}\\
&&= \sum_{\pi_{1},\pi_{2}}\int_{\{\sum_{j=1}^{2n} r_{1,j}\leq t_{1}\} } \prod_{j=1}^{2n}p_{r_{1,j}}(x_{\pi_{1}(j)}-x_{\pi_{1}(j-1)}) \nonumber   \\
&&\hspace{1 in} \int_{\{\sum_{j=1}^{2n} r_{2,j}\leq t_{2}\} }
p_{r_{2,1}+(t_{1}-\sum_{j=1}^{2n} r_{1,j})}(y_{\pi_{2}(1)}-x_{\pi_{1}(2n)})\nonumber\\
&&\hspace{1.5 in} \prod_{j=2}^{2n}p_{r_{2,j}}(y_{ \pi_{2}(j)}-y_{ \pi_{2}(j-1)})\prod_{j=1}^{2n}\,dr_{1,j}\,dr_{2,j} \nonumber  
\end{eqnarray}
and
\begin{eqnarray}
&&E\(\prod_{i=1}^{2n}      \wt L^{ z_{i}}_{t_{3}}    \)
\label{semi.3a}\\
&& = \sum_{\pi_{3}}\int_{\{\sum_{j=1}^{2n} r_{3,j}\leq t_{3}\} } \prod_{j=1}^{2n}p_{r_{3,j}}(z_{ \pi_{3}(j)}-z_{ \pi_{3}(j-1)})\,dr_{3,j} \nonumber
\end{eqnarray}
where the sums run over all permutations  $\pi_{j}$  of $\{1,\ldots, 2n\}$ and we also define $\pi_{j}(0)   =0$
and $x_{0}  =z_{0} =0$.

We then use the product
rule (\ref{pr1}) as before, 
to expand the right hand side of (\ref{semi.2}) into a sum of many terms, and then setting $y_{i}=z_{i}=x_{i}$ we obtain:
\bea && E\(\prod_{i=1}^{2n}\(     \De_{x_{i}}^{h}L^{ x_{i}}_{t_{1}}  \,\(\De_{x_{i}}^{h}  L^{ x_{i}}_{t_{2} }\circ\th_{t_{1}}\) \, \De_{x_{i}}^{h}\wt L^{ x_{i}}_{t_{3}}   \) \)
 \label{f1.20gi}\\ &&\qquad= \sum_{ \pi ,a }\int \mathcal{T}^{\sharp}_{h}( x;\,\pi ,a )\,dx\nn
\eea 
where   $x=(x_{1},\ldots, x_{2n}), \pi=(\pi_{1},\pi_{2}, \pi_{3}), a=(a_{1},a_{2}, a_{d})$ and 
\begin{eqnarray}
&&
\mathcal{T}^{\sharp}_{h}( x;\,\pi ,a ) \label{f1.21gi} \\
&& = \prod_{d=1}^{3}\int_{\RR_{d} }  \prod_{ j=1}^{ 2n}\(\(\De^{ h}_{ x_{ \pi_{d}( j)}}\)^{a_{d, 1}(j)}
\(\De^{ h}_{ x_{ \pi_{d}( j-1)}}\)^{a_{d, 2}(j)}\,p^{\sharp}_{\bar r_{d,j}}(x_{\pi_{d} (j)}-x_{\pi_{d} (j-1)})\)\nn \\
&&\hspace{4 in}\prod_{ j=1}^{ m}\,dr_{d,j}.\nn
\end{eqnarray}
 In (\ref{f1.20gi}) the sum runs over all  triples of permutations $(\pi_{1},\pi_{2}, \pi_{3})$  and   all   
$a_{d}  =(a_{d, 1},a_{ d,2})\,:\,[1,\ldots, 2n]\mapsto \{ 0,1\} \times \{ 0,1\} $, with the
restriction that for each $d,i$ there is exactly one  factor  of the form $\De^{
h}_{ x_{\pi_{d} (i)}}$.  (Here we define $(\De_{x_{i}}^{h})^{0}=1 $ and $(\De_{0}^{h}) =1 $. We have also set $\pi_{1} (0)=\pi_{3} (0)=0$ and $\pi_{2} (0)=\pi_{1} (2n)$.)  In (\ref{f1.21gi}) we set  $\RR_{d}=\{\sum_{j=1}^{2n} r_{d,j}\leq t_{d}\}$, $p_{r}^{\sharp}(x)$ may be either 
$p_{r} (x), p_{r} (x+h)$ or $p_{r} (x-h)$, and $\bar r_{d,j}=  r_{d,j}$ unless $d=2,j=1$ in which case $\bar r_{2,1}=r_{2,1}+(t_{1}-\sum_{j=1}^{2n} r_{1,j})$.
It is important to recognize that in  the right hand side of (\ref{f1.21gi}) each difference operator is applied to only one of the  terms $p_{\cdot}(\cdot)$.

  Instead of    (\ref{f1.21gi}) we first analyze 
  \begin{eqnarray}
  &&
\mathcal{T} _{h}( x;\,\pi ,a ) \label{f1.21gm} \\
&& = \prod_{d=1}^{3}\int_{\RR_{d} }  \prod_{ j=1}^{ 2n}\(\(\De^{ h}_{ x_{ \pi_{d}( j)}}\)^{a_{d, 1}(j)}
\(\De^{ h}_{ x_{ \pi_{d}( j-1)}}\)^{a_{d, 2}(j)}\,p_{\bar r_{d,j}}(x_{\pi_{d} (j)}-x_{\pi_{d} (j-1)})\)\nn \\
&&\hspace{4 in}\prod_{ j=1}^{ m}\,dr_{d,j}.\nn
\end{eqnarray}
  This differs from     (\ref{f1.21gi}) in that we have replaced all $p_{r}^{\sharp}(x)$ by
$p_{r} (x)$. At before it   will be seen that this has no effect on the asymptotics.      
      
As before, we first consider the case that $a_{1}=a_{2}=a_{3}=e$
where    $e=(e(1),\ldots,e(2n))$ and
 $e(2j)=(1,1),\,e(2j-1)=(0,0)$,   $j=1\ldots n$. 
Let $\mathcal{P}=\{(l_{2i-1},l_{2i})\,,\,1\leq i\leq n\}$ be a pairing of the integers $[1,2n]$. Let $\pi_{1},\pi_{2}, \pi_{3}$   be  permutations  of   $[1,2n]$   such that for each $1\leq d\leq 3, 1\leq j\leq n$, 
$\{\pi_{d}(2j-1), \pi_{d}(2j)\}=\{l_{2i-1},l_{2i}\}$ for some, necessarily unique, $ 1\leq i\leq n$. In this case we say that $\pi $ is compatible with the pairing $\mathcal{P}$. (  Note   that $\{\pi_{d}(2j-1), \pi_{d}(2j)\} $  is  not necessarily the same for each $d$.) Then  by (\ref{f1.21gm})  
\begin{eqnarray}
&& 
\mathcal{T}_{h}( x;\,\pi,e) \label{f1.21gia}\\
&&=\prod_{d=1}^{3}\int_{\{\sum_{j=1}^{n} r_{d,j}+s_{d,j}\leq t_{d}\} }  \prod_{ j=1}^{ n}\(\De^{ h}\De^{ -h}
\,p_{r_{d,j}}(x_{\pi_{d}(2j)}-x_{\pi_{d}(2j-1)})\)\,dr_{d,j}\nn\\ &&\hspace{
1.1in}\times\prod_{ j=1}^{ n}\, p_{\bar s_{d,j}}(x_{\pi_{d}(2j-1)}-x_{\pi_{d}(2j-2)}) \,\, \,ds_{d,j}.\nn
\end{eqnarray}
where again $\bar s_{d,j}= s_{d,j}$ unless $(d,j)=(2,1)$ in which case we have 
$\bar s_{2,1}= s_{2,1}+(t_{1}-\sum_{j=1}^{n} r_{1,j}+s_{1,j})$.

Set $\si_{d} (j)=i $ when  $\{\pi_{d}(2j-1), \pi_{d}(2j)\}=\{l_{2i-1},l_{2i}\} $. Using the approach of Sub-section \ref{ss-3.1t} together with the estimates of Lemma \ref{lem-vpropt} in place of the estimates of Lemma \ref{lem-vprop} (in estimating error terms we take absolute values of all integrands and extend the time integration of each term to $[0,T]$ with $T=2\max (t_{1},t_{2},t_{3})$) we can show that 
\begin{equation} 
\int  
 \mathcal{T}_{h}( x;\,\pi ,e )\prod_{ j=1}^{ 2n}\,dx_{j}=\int  
\wt \mathcal{T}_{h}( x;\,\pi ,e )\prod_{ j=1}^{ 2n}\,dx_{j}+O(h^{4n+1/2 }) \label{9.43}
\end{equation}
where 
\begin{eqnarray}
\lefteqn{
\wt \mathcal{T}_{h}( x;\,\pi ,e)  
 =\prod_{d=1}^{3}\int_{\{\sum_{j=1}^{n} r_{d,j}+s_{d,j}\leq t_{d}\} }  \prod_{ i=1}^{ n}\De^{ h}\De^{ -h}
\,p_{r_{d,i}}(x_{l_{2i}}-x_{l_{2i-1}})\,  dr_{d,i}}\nn\\ &&\hspace{
1in}\times\prod_{ j=1}^{ n}\,  p_{\bar s_{d,j}}(x_{l_{2\si_{d} (j)-1}}-x_{l_{2\si_{d} (j-1)-1}})\,\, \,ds_{d,j}\label{f1.21giaww}. 
\end{eqnarray}
The fact that the error term in (\ref{9.43}) is $O(h^{4n+1/2 })$ and not $O(h^{4n+1  })$
is due to the fact that we use (\ref{9.30gb}) instead of (\ref{1.30gb}).

Let   $\wt A_{h} (\pi ,e )$ denote  the integral on the right hand side of (\ref{9.43}) so that
\begin{eqnarray}
  \lefteqn{
\wt A_{h} (\pi ,e )\nn}\\
&&\label{f9.44d}=\int \prod_{d=1}^{3}\int_{\{\sum_{j=1}^{n} r_{d,j}+s_{d,j}\leq t_{d}\} }  \prod_{ i=1}^{ n}\De^{ h}\De^{ -h}
\,p_{r_{d,i}}(x_{l_{2i}}-x_{l_{2i-1}})\,  dr_{d,i}\nn\\ &&\hspace{
1in}\times\prod_{ j=1}^{ n}\,  p_{\bar s_{d,j}}(x_{l_{2\si_{d} (j)-1}}-x_{l_{2\si_{d} (j-1)-1}})\,\, \,ds_{d,j}
\prod_{i=1}^{2n}\,dx_{i}.\nn
\eea
  We make the change of   variables $x_{l_{2i}}\to  x_{l_{2i}}+x_{l_{2i-1}}  $, $i=1,\ldots,n$ and  write this as  
\bea
&& \int \prod_{d=1}^{3}\int_{\{\sum_{j=1}^{n} r_{d,j}+s_{d,j}\leq t_{d}\} }  \prod_{ i=1}^{ n}\De^{ h}\De^{ -h}
\,p_{r_{d,i}}(x_{l_{2i}})\,  dr_{d,i}\nn\\ &&\hspace{
1in}\times\prod_{ j=1}^{ n}\,  p_{\bar s_{d,j}}(x_{l_{2\si_{d} (j)-1}}-x_{l_{2\si_{d} (j-1)-1}})\,\, \,ds_{d,j}
\prod_{i=1}^{2n}\,dx_{i}..\nn
\end{eqnarray}

We  now rearrange the integrals with respect to $x_{2},x_{4},\ldots,x_{2n}$ and get
\begin{eqnarray} 
&& \wt A_{h} (\pi ,e )\label{simp.1} \\ 
&&=
\int \(\int_{\{\sum_{j=1}^{n} r_{d,j}+s_{d,j}\leq t_{d}, \forall d\} } \right.  \prod_{i=1}^{n} \(\int \( \prod_{d=1}^{3} \De^{ h}\De^{ -h}
\,p_{r_{d,i}}(x ) \)  \,dx\)   \nn\\
&&\hspace{
.2in}  \times\prod_{d=1}^{3} \prod_{ j=1}^{ n}\,   p_{\bar s_{d,i}}(x_{l_{2\si_{d} (j)-1}}-x_{l_{2\si_{d} (j-1)-1}})\,\left.  \,ds_{d,i}\,dr_{d,i}  \) \prod_{i=1}^{n}\,dx_{l_{2i-1}}.\nn
\end{eqnarray}
Let
\begin{eqnarray}
\lefteqn{ F(\si ;s)
\nn}\\
&&\,\,: = \int \prod_{d=1}^{3}\prod_{ j=1}^{ n}\,   p_{\bar s_{d,i}}(x_{l_{2\si_{d} (j)-1}}-x_{l_{2\si_{d} (j-1)-1}})  \prod_{i=1}^{n}\,dx_{l_{2i-1}} \label{f9.2d}\\
&& \,\,\, =\int \prod_{d=1}^{3}\prod_{ j=1}^{ n}\,   p_{\bar s_{d,i}}(y_{ \si_{d} (j)  }-y_{ \si_{d} (j-1) })\, \prod_{i=1}^{n}\,dy_{i},   \nonumber
\end{eqnarray} 
  where we set $y_{i}=x_{l_{2i-1}}$.  
We can now write
\begin{eqnarray} 
&&
 \wt A_{h} (\pi ,e )\nn\\ 
&&=\label{9.44b}
  \int  \int_{\{\sum_{i=1}^{n}r_{d,i}+s_{d,i}\leq t_{d}, \forall d\}}      
  F(\si ;s)\prod_{d=1}^{3} \prod_{i=1}^{n} \,ds_{d,i}  \\ 
  &&
 \hspace{1 in} \prod_{i=1}^{n} \(\int \prod_{d=1}^{3}\(  \De^{ h}\De^{ -h}
\,p_{r_{d,i}}(x ) \)  \,dx\)\prod_{d=1}^{3}\,dr_{d,i} .\nn
\end{eqnarray}
Using $2-e^{ih\la } -e^{-ih\la }=2-2\cos (\la h)=4\sin^{2} (\la h/2)$ we can write
\begin{eqnarray}
&&G_{h}(r)=:\int \prod_{d=1}^{3}\(  \De^{ h}\De^{ -h}
\,p_{r_{d,i}}(x ) \)  \,dx
 \label{fou.1}\\
&& =\int \(\prod_{d=1}^{3}\( {1 \over 2\pi}\int e^{ix\la_{d,i}} 
\(2-e^{ih\la_{d,i}} -e^{-ih\la_{d,i}} \) e^{-r_{d,i}\la^{2}_{d,i}/2}\,d\la_{d,i}  \) \) \,dx  \nonumber\\
&& =\({4 \over 2\pi }\)^{3}\int  \(  \int e^{ix\sum_{d=1}^{3}\la_{d,i}} \prod_{d=1}^{3}
\sin^{2} (\la_{d,i} h/2) e^{-r_{d,i}\la^{2}_{d,i}/2}\,d\la_{d,i}   \) \,dx  \nonumber\\
&& =\({4 \over 2\pi }\)^{3} \int  \( \int e^{ix\sum_{d=2}^{3}\la_{d,i}} \(  \int    
e^{ix \la_{d,1}}\sin^{2} (\la_{1,i} h/2) e^{-r_{1,i}\la^{2}_{1,i}/2}\,d\la_{d,1}\)\,dx\right.   \nonumber \\
&& \hspace{2.5 in}\left. \prod_{d=2}^{3}
\sin^{2} (\la_{d,i} h/2) e^{-r_{d,i}\la^{2}_{d,i}/2}\,d\la_{d,i}   \)   \nn\\
&& =\({4 \over 2\pi }\)^{3}   \(  \int    
\sin^{2} (\la_{1,i} h/2) e^{-r_{1,i}\la^{2}_{1,i}/2}\right.   \nonumber \\
&& \hspace{2.5 in}\left. \prod_{d=2}^{3}
\sin^{2} (\la_{d,i} h/2) e^{-r_{d,i}\la^{2}_{d,i}/2}\,d\la_{d,i}   \)  \nn
\end{eqnarray}
with $\la_{1,i}=:-\sum_{d=2}^{3}\la_{d,i}$ in the last equality. For the last equality we used Fourier inversion.

Since $G_{h}, F\geq 0$, we have the following upper and lower bounds for $ \wt A_{h} (\pi ,e )$
\begin{eqnarray} 
\lefteqn{  \wt A_{h} (\pi ,e )\label{f9.25ga}\nn}\\ &&\leq 
\(\int_{[0,\ff]^{3}}\(\int \prod_{d=1}^{3}\(\De^{ h}\De^{ -h}\,p_{r_{d} }(x)\)
 \,dx\)\,    \,dr_{d}    \)^{n}\label{9.44b}\\ &&\hspace{
.5in}\times   \int_{\{\sum_{i=1}^{n} s_{d,i}\leq t_{d}, \forall d\}}    
  F(\si ;s ) \prod_{d=1}^{3}\prod_{i=1}^{n} \,ds_{d,i}   \nn
\end{eqnarray} 
and 
\begin{eqnarray} 
\lefteqn{  \wt A_{h} (\pi ,e )\nn}\\ &&\geq 
 \(\int_{[0,h]^{3}}\(\int \prod_{d=1}^{3}\(\De^{ h}\De^{ -h}\,p_{r_{d} }(x)\)
  \,dx\)\,    \,dr_{d}   \)^{n}\label{f9.25ga}\\ &&\hspace{
.5in}\times  \int_{\{\sum_{i=1}^{n} s_{d,i}\leq t_{d}-nh, \forall d\}}     
  F(\si ;s ) \prod_{d=1}^{3}\prod_{i=1}^{n} \,ds_{d,i}   \nn .
\end{eqnarray}
   We show that the two sides of the inequalities are asymptoticallly equivalent as $h\to 0$.
The following Lemma is proven below. 
\bl\label{lem-spread}
\begin{eqnarray}
\lefteqn{\int_{\{\sum_{i=1}^{n} s_{d,i}\leq t_{d}, \forall d\}}     
  F(\si ;s )  \prod_{d=1}^{3}\prod_{i=1}^{n} \,ds_{d,i}   
\nn}\\
&& -     \int_{\{\sum_{i=1}^{n} s_{d,i}\leq t_{d}-nh, \forall d\}}     
  F(\si ;s ) \prod_{d=1}^{3} \prod_{i=1}^{n} \,ds_{d,i}   \label{f9.4d} \leq C_{T}h.\nonumber
\end{eqnarray}
\el
Referring to (\ref{9.43}), using (\ref{big.1})-(\ref{big.2}) 
 we see that 
 \begin{eqnarray} 
 && 
\int  
 \mathcal{T}_{h}( x;\,\pi ,e )\prod_{ j=1}^{ 2n}\,dx_{j}=\wt A_{h} (\pi ,e ) +O(h^{4n+1/2})
  \label{f9.5}\\
 && =
(8h^{4} )^{n} \int_{\{\sum_{i=1}^{n} s_{d,i}\leq t_{d}, \forall d\}}     
  F(\si ;s ) \prod_{d=1}^{3}\prod_{i=1}^{n} \,ds_{d,i}+O(h^{4n+1/2}) \nn\\
&& =(8h^{4} )^{n}\nonumber\\
&&\hspace{.3 in}\int  \(\prod_{d=1}^{3}\int_{\{\sum_{i=1}^{n} s_{d,i}\leq t_{d}\}}\prod_{ i=1}^{ n}\, p_{\bar s_{d,i}}(y_{\si_{d} (i)}-y_{\si_{d} (i-1)})\,\prod_{i=1}^{n} \,ds_{d,i}\) \prod_{i=1}^{n}\,dy_{i}\nn\\
&&\hspace{3 in}+O(h^{4n+1/2}).   \nonumber
 \end{eqnarray}
 
  Recall that, in the paragraph containing (\ref{f1.21gia}),  for a given pairing $\mathcal{P}=\{(l_{2i-1},l_{2i})\,,\,1\leq i\leq n\}$  of the integers $[1,2n]$,  we   define  what it means for a collection of permutations $\pi=(\pi_{1},\pi_{2}, \pi_{3})$ of $[1,2n]$ to be compatible with $\mathcal{P}$. We write this  as $(\pi_{1},\pi_{2}, \pi_{3}) \sim \mathcal{P}$.  Obviously, there are many such pairs.  We can interchange the two elements of the pair $
\pi_{d}(2j-1), \pi_{d}(2j)$  without changing (\ref{f9.5}).   There are   $2^{3n}$ ways to do this. Furthermore, by permuting the pairs $\{\pi_{d}(2j-1), \pi_{d}(2j)\}$ we give rise in (\ref{f9.5}) to all possible permutations  $\si_{d}$ of $[1,n]$. We thus obtain 
 \begin{eqnarray} 
\lefteqn{ \sum_{(\pi_{1},\pi_{2}, \pi_{3})\sim \mathcal{P}}\int  
 \mathcal{T}_{h}( x;\,\pi ,e )\prod_{ j=1}^{ 2n}\,dx_{j} \nn}\\
&&\,\, =(2^{3}8h^{4} )^{n}\sum_{ \si }\int  \(\prod_{d=1}^{3}\int_{\{\sum_{i=1}^{n} s_{d,i}\leq t_{d}\}}\right.
  \nn\\
&&\hspace{1 in}\left. \prod_{ i=1}^{ n}\, p_{\bar s_{d,i}}(y_{\si_{d} (i)}-y_{\si_{d} (i-1)})\,\prod_{i=1}^{n} \,ds_{d,i}\) \prod_{i=1}^{n}\,dy_{i}+O(h^{4n+1/2})  \nn\\ 
&&=(64h^{4} )^{n}E\lc\(\int   L^{ x}_{ t_{1}}(L^{ x}_{ t_{2}}\circ\th_{t_{1}}) \wt L^{ x}_{ t_{3}} \,dx\)^{ n}\rc+O(h^{4n+1/2}). \label{f9.45r}
 \end{eqnarray}
Here the  sum  in the second line runs over all   permutations  $\si=(\si_{1},\cdots, \si_{q})$ of $\{1,\ldots, n\}$ and we set $\si_{d}(0)= 0$. The      fourth   line follows from Kac's moment formula.  

Since   there are ${( 2n)!\over 2^{ n}n!}$ pairings of the $2n$
elements $\{1,\ldots, m=2n\}$ we see that
 \begin{eqnarray} \quad
&&\sum_{  \mathcal{P}}\sum_{(\pi_{1},\pi_{2}, \pi_{3})\sim \mathcal{P}}\int \mathcal{T}_{h}( x;\,\pi ,e) \,\prod_{ j=1}^{ 2n}\,dx_{j}\label{f9.45s}\\
&&   =  {( 2n)!\over 2^{ n}n!} (64h^{4} )^{n}E\lc\(\int   L^{ x}_{ t_{1}}(L^{ x}_{ t_{2}}\circ\th_{t_{1}}) \wt L^{ x}_{ t_{3}} \,dx\)^{ n}\rc+O(h^{4n+1/2})\nn
 \end{eqnarray}
 where the first sum runs over all pairings $\mathcal{P}$ of $\{1,\ldots, 2n\}$.
 
 Given the estimates of Lemma \ref{lem-vpropt} we can show as in Section 4 that the contributions to (\ref{f1.20gi}) from  $(a_{1}, a_{2}, a_{3})= (e,e,e)$ for  $\pi$ not compatible with a pairing is $O(h^{4n+1/2})$. The arguments of Section 4 will give a similar bound for $(a_{1}, a_{2}, a_{3})\neq (e,e,e)$    with one possible exception. This will happen if $\pi_{2}(1)=\pi_{1}(2n)$ so that the argument of the term $p_{\bar s_{2,1}}$ is zero, and two $\De$ operators are applied to this $p$. In that case, since $\De^{h}\De^{-h}p_{\bar s_{2,1}}(0)=2\De^{h}p_{\bar s_{2,1}}(0)$, we seem to have lost one  $\De$ operator, all of which are used in Section 4 to obtain the required error estimate. The remedy will be found in the special nature of $\bar s_{2,1}$ as we now explain.

 Instead of extending  the time integration of each term to $[0,T]$, we first consider the region where
  $\sum_{j=1}^{n} r_{1,j}+s_{1,j}\leq t_{1}/2$ so that $\bar s_{2,1}\geq t_{1}/2$. Then
     \begin{equation}
 |\De^{h}p_{\bar s_{2,1}}(0)|={c| e^{-h^{2}/2\bar s_{2,1}} -1| \over \bar s^{1/2}_{2,1}}\leq {c h^{2}  \over \bar s_{2,1}^{3/2}}\leq c(t_{1})h^{2}.\label{cn.1}
 \end{equation}
 We then extend   the time integration of each term to $[0,T]$ and proceed as before. On the other hand, if  $\sum_{j=1}^{n} r_{1,j}+s_{1,j}\geq t_{1}/2$, then for some $j$ we have either 
 $r_{1,j}\geq \de=:t_{1}/(4n)$ or $s_{1,j}\geq \de$. Say it is the latter. We then use the $s_{1,j}$
integration for the bound, see (\ref{cn.1}),
\begin{eqnarray}
&& \int_{0}^{T}\int_{0}^{T}|\De^{h}p_{  s_{2,1}+s_{1,j}}(0)|\,ds_{2,1}\,ds_{1,j}
\label{cn.3}\\
&&\leq ch^{2}   \int_{0}^{T}\int_{0}^{T}{ 1  \over (s_{2,1}+s_{1,j})^{3/2}}\,ds_{2,1}\,ds_{1,j}\leq ch^{2}\nonumber
\end{eqnarray} 
 and bound the other  term involving $  s_{1,j} $, be it $ p_{s_{1,j}}(x),  |\De^{h}p_{s_{1,j}}(x)|$ or $ |\De^{h}\De^{-h}p_{s_{1,j}}(x)|$, by its supremum over $s_{1,j}\geq \de $, using Lemma \ref{lem-vpropd}.
  \qed
  
{\bf    Proof of Lemma \ref{lem-spread}: }

Let $A\subseteq [0,t_{1}]^{n}\times [0,t_{2}]^{n}\times [0,t_{3}]^{n}$ and set 
\begin{equation}
I(A)= \int_{A}      
  F(\si ;s)\prod_{d=1}^{3} \prod_{i=1}^{n} \,ds_{d,i}.\label{F.1}
\end{equation}
To prove Lemma \ref{lem-spread} it suffices to show that
\begin{equation}
I(A)\leq C_{T}|A|^{1/2}.\label{F.2}
\end{equation}

We have 
\begin{eqnarray}
&&I(A)
\label{F.3}\\
&&=\int_{A} \(\int \prod_{d=1}^{3}\prod_{ j=1}^{ n}\,   p_{s_{d,i}}(y_{ \si_{d} (j)  }-y_{ \si_{d} (j-1) })\, \prod_{i=1}^{n}\,dy_{i}\)        \prod_{d=1}^{3} \prod_{i=1}^{n} \,ds_{d,i}\nonumber\\
&&= \int \(\int_{A} \prod_{d=1}^{3}\prod_{ j=1}^{ n}\,   p_{s_{d,i}}(y_{ \si_{d} (j)  }-y_{ \si_{d} (j-1) })\,  \,ds_{d,i} \)          \prod_{i=1}^{n}\,dy_{i}.\nonumber
\end{eqnarray}
Then by the Cauchy-Schwarz inequality
\begin{eqnarray}
&&I(A)
\label{F.4}\\
&&\leq  |A|^{1/2} \int \(\int_{[0,T]^{3n}} \prod_{d=1}^{3}\prod_{ j=1}^{ n}\,   p^{2}_{s_{d,i}}(y_{ \si_{d} (j)  }-y_{ \si_{d} (j-1) })\,  \,ds_{d,i} \)^{1/2}          \prod_{i=1}^{n}\,dy_{i}\nonumber\\
&&\leq  |A|^{1/2} e^{3nT}\int \( \prod_{d=1}^{3}\prod_{ j=1}^{ n}\,  f(y_{ \si_{d} (j)  }-y_{ \si_{d} (j-1) })  \)^{1/2}          \prod_{i=1}^{n}\,dy_{i}\nonumber
\end{eqnarray}
where
\begin{equation}
f(y)=\int_{0}^{\ff} e^{-s  }  p^{2}_{s }(y)\,  \,ds. \label{F.5}
\end{equation}
Since $f(y)$ is the 1-potential density of planar Brownian motion evaluated at $(\sqrt{2}\,\,y,0)$, we know that $f(y)$ has a logaritmic singularity at $y=0$ and has exponential falloff at $\ff $ so that the last integral in (\ref{F.4}) is finite.
\qed

\section{Proof of Lemma \ref{lem-var}}\label{sec-var}

\bl\label{lem-var}Fix $T<\ff$. For all $s,t\leq T$
\bea
&&
\E\bigg[ \(
\int \lc  (\De_{x}^{h}L^{ x}_{t})^{2}-4hL^{ x}_{t}      \rc\De_{x}^{h}\wt L^{ x}_{s} \,dx\)^{2}\bigg]=
32 h^{4}E\(\int  ( L^{ x}_{t})^{2}  \wt L^{ x}_{s} \,dx\)\nonumber\\
&&\hspace{2.6 in}+O\((s\wedge t)^{\ep}h^{4+\ep}\).\label{va.0}
\eea
\el

{\bf  Proof of Lemma \ref{lem-var}: }In order to prove (\ref{va.0}) we must make use of the subtraction on the left hand side  to eliminate all terms which are not $O\( h^{4 }\)$, then isolate the main contribution which is the first term on the right hand side, and  estimate all error terms. As we will see the terms which are not $O\( h^{4 }\)$ come from `bound' variables.  Because we are not using exponential times, the subtractions do not exactly eliminate all bound variables, which makes the analysis more complicated than in previous sections. 

We first write
\begin{equation}
\E\bigg[ \(
\int \lc  (\De_{x}^{h}L^{ x}_{t})^{2}-4hL^{ x}_{t}      \rc\De_{x}^{h}\wt L^{ x}_{s} \,dx\)^{2}\bigg]=I_{1}-8hI_{2}+16h^{2}I_{3}\label{rv.1}
\end{equation}
where
\begin{eqnarray}
&&I_{1}=\E\bigg[ \(
\int \(\De_{x}^{h}L^{ x}_{t} \)^{2}     \De_{x}^{h}\wt L^{ x}_{s} \,dx\)^{2}\bigg] \nonumber\\
&&
=\E\bigg[  
\int \(\De_{x}^{h}L^{ x}_{t} \)^{2}     \De_{x}^{h}\wt L^{ x}_{s} \,dx\int \(\De_{y}^{h}L^{ y}_{t} \)^{2}     \De_{y}^{h}\wt L^{ y}_{s} \,dy\bigg] 
\label{km.10a}\\
&&=\int\int E\( \(\De_{x}^{h}L^{ x}_{t} \)^{2} \(\De_{y}^{h}L^{ y}_{t} \)^{2}    \)E\( \De_{x}^{h}\wt L^{ x}_{s}  \De_{y}^{h}\wt L^{ y}_{s}  \) \,dx\,dy    \nn
\end{eqnarray}
\begin{eqnarray}
&&I_{2}=\E\bigg[  
\int  L^{ x}_{t}   \De_{x}^{h}\wt L^{ x}_{s} \,dx\int \(\De_{y}^{h}L^{ y}_{t} \)^{2}    \De_{y}^{h}\wt L^{ y}_{s} \,dy\bigg] 
\label{km.10b}\\
&&=\int\int E\(  L^{ x}_{t} \(\De_{y}^{h}L^{ y}_{t} \)^{2}   \)E\( \De_{x}^{h}\wt L^{ x}_{s}  \De_{y}^{h}\wt L^{ y}_{s}  \) \,dx\,dy    \nn
\end{eqnarray}
and
\begin{eqnarray}
&&I_{3}=\E\bigg[ \(
\int  L^{ x}_{t}       \De_{x}^{h}\wt L^{ x}_{s} \,dx\)^{2}\bigg] \nonumber\\
&&
=\E\bigg[  
\int  L^{ x}_{t}     \De_{x}^{h}\wt L^{ x}_{s} \,dx\int  L^{ y}_{t}      \De_{y}^{h}\wt L^{ y}_{s} \,dy\bigg] 
\label{km.10c}\\
&&=\int\int E\(  L^{ x}_{t}   L^{ y}_{t}     \)E\( \De_{x}^{h}\wt L^{ x}_{s}  \De_{y}^{h}\wt L^{ y}_{s}  \) \,dx\,dy.    \nn
\end{eqnarray}

By Kac's moment formula and (\ref{pr1}) we have 
    \begin{equation}
G_{s}(x,y) =:E\( \De_{x}^{h}\wt L^{ x}_{s}  \De_{y}^{h}\wt L^{ y}_{s}  \) =\int_{\{s_{1}+s_{2}\leq s\}} F_{s}(x,y) \,ds_{1}\,ds_{2}
\label{km.11b}
  \end{equation}
  where  \bea
  &&
F_{s}(x,y)\label{km.20b}\\
&&= 
 \De^{h}p_{s_{1}}(x)\,\De^{h}p_{s_{2}}(y-x-h) + 
p_{s_{1}}(x)\,\De^{h}\De^{-h}p_{s_{2}}(y-x)\nn\\
&&+ 
 \De^{h}p_{s_{1}}(y)\,\De^{h}p_{s_{2}}(x-y-h) + 
p_{s_{1}}(y)\,\De^{h}\De^{-h}p_{s_{2}}(x-y).\nn
  \eea
 For any $\ep>0$
  \begin{eqnarray}
  &&|G_{s}(x,y)| \leq cs^{\ep/2 }v^{1-\ep}_{s}(x)v_{s}(y-x-h)+cs^{\ep/2 }u^{1-\ep}_{s}(x)w_{s}(y-x)
\nn\\
  &&\hspace{.5 in}+cs^{\ep/2 }v^{1-\ep}_{s}(y)v_{s}(x-y-h)+cs^{\ep/2 }u^{1-\ep}_{s}(y)w_{s}(x-y).    \label{gbd.1}
  \end{eqnarray} 
To see this we note the bounds
  \begin{equation}
  \int_{0}^{s}p_{r}(x)\,dr\leq \int_{0}^{s}p_{r}(0)\,dr\leq cs^{1/2}\label{int.1}
  \end{equation}
  and
    \begin{equation}
  \int_{0}^{s}|\De^{h}p_{r}(x)|\,dr\leq 2\int_{0}^{s}p_{r}(0)\,dr\leq cs^{1/2}.\label{int.1v}
  \end{equation}
and interpolate  to obtain
\begin{equation}
u_{s}(x)\leq cs^{\ep/2 }u^{1-\ep}_{s}(x),\,\hspace{.4 in}v_{s}(x)\leq cs^{\ep/2 }v^{1-\ep}_{s}(x). \label{int.2}
\end{equation}

It follows from  (\ref{gbd.1}) and  Lemma \ref{lem-vpropt} that for  any $\ep>0$
  \begin{equation}
 \int |G_{s}(x,y)|\,dx\,dy\leq c s^{\ep/2 }h^{2-\ep}.\label{gbd.2}
  \end{equation}
  
  Clearly
  \begin{equation}
E\(  L^{ x}_{t}  L^{ y}_{t}  \)=\int_{\{t_{1}+t_{2}\leq t\}} \(A_{t_{1},t_{2} }(x,y)+A_{t_{1},t_{2} }(y,x)\)\,dt_{1}\,dt_{2}
\label{km.11c}
  \end{equation}
  where
  \begin{equation}
 A_{t_{1},t_{2} }(x,y)=
 p_{t_{1}}(x)p_{t_{2}}(y-x). \label{km.20a}
  \end{equation} 

By Kac's moment formula and (\ref{pr1}), compare (\ref{1.21g}),
  \begin{eqnarray}
&&E\(  L^{ x}_{t}  \(\De_{y}^{h}L^{ y}_{t} \)^{2}    \)
\label{km.11b}\\
&&=2 \sum_{\pi',a' }  \int_{\{\sum_{i=1}^{3}t_{i}\leq t\}}\prod_{i=1}^{3}
\(\De_{ \pi'(i)}^{h}\)^{a'_{1}(i)}
\(\De_{ \pi'(i-1)}^{h}\)^{a'_{2}(i)}p^{\sharp}_{t_{i}}( \pi'(i) - \pi'(i-1) )\,dt_{i} \nonumber
\end{eqnarray}
where  the sum runs over all  maps $\pi'\,:\,[1,2, 3]\mapsto
\{x,y\}$ with $|\pi'^{ -1}(x )|=1,\,|\pi'^{ -1}(y )|=2$, and all `assignments' 
$a'=(a'_{ 1},a'_{ 2})\,:\,[1,2, 3]\mapsto \{ 0,1\}\times \{ 0,1\}$ with the
property that  there will be exactly two factors of the form $\De^{
h}_{ y}$ in (\ref{km.11b}) and none of the form $\De^{
h}_{ x}$. The factor $2$    comes from 
  the fact that $|\pi'^{ -1}(x )|=1,\,|\pi'^{ -1}(y )|=2$. Recall that $p^{\sharp}_{t }(x)$ can be 
  $p_{t }(x), p_{t }(x+h)$ or $p_{t }(x-h)$, but we always have $\De^{h}\De^{-h}p^{\sharp}_{t }(x)=\De^{h}\De^{-h}p_{t }(x)$. Also, we always take the $p_{t }(\cdot)$ for a bound variable to be the $g$ in (\ref{pr1}).
  
  Bound variables can come only from $\pi'_{1}=(x,y,y)$ and $\pi'_{2}=( y,y,x)$. Setting
  \begin{equation}
  f_{t}(h)=p_{t }(0)-p_{t }(h)\label{eff}
  \end{equation}
  we can write the contributions of $\pi'_{1}$ and $\pi'_{2}$ arising from a bound variable as
  \begin{eqnarray}
&&\wt D_{\pi'_{1},t}(x,y)=p_{t_{1}}(x) \,  p_{t_{2}}(y-x)\,
\(\De^{h}\De^{-h}p_{t_{3}}(0)\)
\label{km.11ba}\\
&&\hspace{.7 in}=2p_{t_{1}}(x) \,  p_{t_{2}}(y-x)\,
f_{t_{3}}(h)  \nonumber\\
&&\hspace{.7 in}=2A_{t_{1},t_{2} }(x,y)\,
f_{t_{3}}(h)  \nonumber
\end{eqnarray}
and
  \begin{eqnarray}
&&\wt D_{\pi'_{2},t}(x,y)=p_{t_{1}}(y)\,    
\(\De^{h}\De^{-h}p_{t_{2}}(0)\)\, p_{t_{3}}(x-y) 
\label{km.11ba}\\
&& \hspace{.6 in} =2 p_{t_{1}}(y)\,    
f_{t_{2}}(h)\, p_{t_{3}}(x-y)  \nonumber\\
&& \hspace{.6 in} = 2A_{t_{1},t_{3} }(y,x)\,    
f_{t_{2}}(h).   \nonumber
\end{eqnarray}
 
The non-bound contributions for $\pi'_{1}$ and $\pi'_{2}$ are 
\begin{equation}
\wt B_{\pi'_{1},t}(x,y)=p^{\sharp}_{t_{1}}(x) \,  \De^{h}p^{\sharp}_{t_{2}}(y-x)\,
\De^{h}p_{t_{3}}(0)\label{}
\end{equation}
and
 \begin{eqnarray}
 &&\wt B_{\pi'_{2},t}(x,y)=  p^{\sharp}_{t_{1}}(y)\,   \, \De^{-h}p^{\sharp}_{t_{2}}(0) \, \De^{-h}p^{\sharp}_{t_{3}}(x-y)
 \label{km.11ba}\\
&& \hspace{.6 in} + \De^{h}p^{\sharp}_{t_{1}}(y)\,   \,p^{\sharp}_{t_{2}}(0) \,  \De^{-h}p^{\sharp}_{t_{3}}(x-y) \nn\\
&&  \hspace{.6 in}+ \De^{h}p^{\sharp}_{t_{1}}(y)\,   \, \De^{h}p^{\sharp}_{t_{2}}(0) \, p^{\sharp}_{t_{3}}(x-y). \nn 
 \end{eqnarray} 
 and in addition there is a term from $\pi'_{3}=(y,x,y)$ which is 
  \begin{eqnarray}
 &&\wt B_{\pi'_{2},t}(x,y)= \De^{h}  p^{\sharp}_{t_{1}}(y)\,p^{\sharp}_{t_{2}}(x-y)  \,\De^{h} p^{\sharp}_{t_{3}}(y-x)  
 \label{km.11ba}\\
&& \hspace{.6 in} + p^{\sharp}_{t_{1}}(y)\, \De^{-h} p^{\sharp}_{t_{2}}(x-y)  \,\De^{h} p^{\sharp}_{t_{3}}(y-x).
         \nn 
 \end{eqnarray}

  We observe that by (\ref{int.2}) and  Lemma \ref{lem-vpropt}, for any $1\leq j \leq 3$
  \begin{equation}
\sup_{x,y} \int_{\{\sum_{i=1}^{3}t_{i}\leq t\}} |\wt  B_{\pi'_{j}, t}(x,y)|\prod_{i=1}^{3}\,dt_{i} \leq ct^{\ep/2}h^{2-\ep}.\label{sb.1}
  \end{equation}
Hence in view of   (\ref{gbd.2})   we see that  for  any $\ep>0$ and $1\leq j \leq 3$
  \begin{equation}
 h\int \( \int_{\{\sum_{i=1}^{3}t_{i}\leq t\}} | \wt  B_{\pi'_{j}, t}(x,y)|\prod_{i=1}^{3}\,dt_{i}\)|G_{s}(x,y)|\,dx\,dy=O\((s\wedge t)^{\ep}h^{4+\ep}\).\label{gbd.22}
  \end{equation}

  Similarly
\begin{eqnarray}
&&E\( \(\De_{x}^{h}L^{ x}_{t} \)^{2} \(\De_{y}^{h}L^{ y}_{t} \)^{2}    \)
\label{km.11a}\\
&&= 4 \sum_{\pi,a }  \int_{\{\sum_{i=1}^{4}t_{i}\leq t\}}\prod_{i=1}^{4}
\(\De_{ \pi(i)}^{h}\)^{a_{1}(i)}
\(\De_{ \pi(i-1)}^{h}\)^{a_{2}(i)}p^{\sharp}_{t_{i}}( \pi(i) - \pi(i-1) )\,dt_{i} \nonumber
\end{eqnarray}
where  the sum runs over all  maps $\pi\,:\,[1,\ldots, 4]\mapsto
\{x,y\}$ with $|\pi^{ -1}(x )|=|\pi^{ -1}(y )|=2$, and all `assignments' 
$a=(a_{ 1},a_{ 2})\,:\,[1,\ldots, 4]\mapsto \{ 0,1\}\times \{ 0,1\}$ with the
property that  there will be exactly two factors of the form $\De^{
h}_{ x}$ in (\ref{km.11a}) and similarly for $\De^{
h}_{ y}$. The factor $4=2^{ 2}$    comes from 
  the fact that $|\pi^{ -1}(x )|=|\pi^{ -1}(y )|=2$.
  
  Writing $\pi$ as a sequence $(\pi(1),\pi(2),\pi(3),\pi(4))$, we first consider 
  $\pi_{1}=(x,x,y,y)$ and $\pi_{2}=( y,y,x,x)$. These are the only $\pi$'s which have two bound variables. We can write the contribution of $\pi_{1}$ arising from two bound variables as
  \bea
  &&
 D_{\pi_{1},t}(x,y)=p_{t_{1}}(x)\,\(\De^{h}\De^{-h}p_{t_{2}}(0)\) p_{t_{3}}(y-x)\,
\(\De^{h}\De^{-h}p_{t_{4}}(0)\)\nn\\
 &&\hspace{.6 in}=4 p_{t_{1}}(x)\,f_{t_{2}}(h) p_{t_{3}}(y-x)\, f_{t_{4}}(h)
\nn\\
 &&\hspace{.6 in}=4 f_{t_{2}}(h) \, f_{t_{4}}(h) A_{t_{1},t_{3} }(x,y)
\label{km.11ab}
  \eea
  and similarly
  \begin{equation}
 D_{\pi_{2},t}(x,y)=4f_{t_{2}}(h) \, f_{t_{4}}(h) A_{t_{1},t_{3} }(y,x). \label{dp2}
  \end{equation}
The contribution of $\pi_{1}$ arising from one  bound variable is
    \bea
  &&
 B_{\pi_{1},t}(x,y)=p^{\sharp}_{t_{1}}(x)\,\(\De^{-h}p^{\sharp}_{t_{2}}(0)\)\,\De^{-h} p^{\sharp}_{t_{3}}(y-x)\,
\(\De^{h}\De^{-h}p_{t_{4}}(0)\)\nonumber\\
&& \hspace{.6 in}+ \De^{h}p^{\sharp}_{t_{1}}(x)\, p^{\sharp}_{t_{2}}(0) \,\De^{-h} p^{\sharp}_{t_{3}}(y-x)\,
\(\De^{h}\De^{-h}p_{t_{4}}(0)\)         \nn\\
&& \hspace{.6 in}+ \De^{h}p^{\sharp}_{t_{1}}(x)\, \De^{h}p^{\sharp}_{t_{2}}(0) \, p^{\sharp}_{t_{3}}(y-x)\,
\(\De^{h}\De^{-h}p_{t_{4}}(0)\)         \label{km.11ac}\\
&& \hspace{.6 in}+  p^{\sharp}_{t_{1}}(x)\,\(\De^{h}\De^{-h}p_{t_{2}}(0)\) \, \De^{h}p^{\sharp}_{t_{3}}(y-x)\,
\De^{h}p^{\sharp}_{t_{4}}(0)         \nn
  \eea
    and similar terms for $\pi_{2}$.
  
  This is also a contribution of $\pi_{3}=(x,y,y,x)$ arising from one  bound variable
    \bea
  &&
 B_{\pi_{3},t}(x,y)= \De^{h}p^{\sharp}_{t_{1}}(x)\,   \, p^{\sharp}_{t_{2}}(y-x)\,
\(\De^{h}\De^{-h}p_{t_{3}}(0)\)\, \De^{h}p^{\sharp}_{t_{4}}(x-y)         \nn\\
&& \hspace{.6 in}+ p^{\sharp}_{t_{1}}(x)\,   \, \De^{-h}p^{\sharp}_{t_{2}}(y-x)\,
\(\De^{h}\De^{-h}p_{t_{3}}(0)\)\, \De^{h}p^{\sharp}_{t_{4}}(x-y)         \nn\\
&& \hspace{.6 in}+ \De^{h}p^{\sharp}_{t_{1}}(x)\,   \, \De^{-h}p^{\sharp}_{t_{2}}(y-x)\,
\(\De^{h}\De^{-h}p_{t_{3}}(0)\)\, p^{\sharp}_{t_{4}}(x-y)         \nn
  \eea
  and similar terms for $\pi_{4}=(y,x,x,y)$.

As before, we observe that by (\ref{int.2}) and  Lemma \ref{lem-vpropt}, for any $1\leq j \leq 4$
  \begin{equation}
\sup_{x,y} \int_{\{\sum_{i=1}^{4}t_{i}\leq t\}} | B_{\pi_{j}, t}(x,y)|\prod_{i=1}^{4}\,dt_{i} \leq ct^{\ep/2}h^{3-\ep}.\label{sb.1}
  \end{equation}
Hence in view of   (\ref{gbd.2})   we see that  for  any $\ep>0$ and $1\leq j \leq 4$
  \begin{equation}
 \int \( \int_{\{\sum_{i=1}^{4}t_{i}\leq t\}} | B_{\pi_{j}, t}(x,y)|\prod_{i=1}^{4}\,dt_{i}\)|G_{s}(x,y)|\,dx\,dy=O\((s\wedge t)^{\ep}h^{4+\ep}\).\label{gbd.2c}
  \end{equation}

Taking note of the factor $4$ in (\ref{km.11a}) and the factor  $2$ in (\ref{km.11b}) we now  show that  
  \bea
  &&
4\int\(\int_{\{\sum_{i=1}^{4}t_{i}\leq t\}} \( D_{\pi_{1},t}(x,y)+ D_{\pi_{2},t}(x,y)\)\prod_{i=1}^{4}\,dt_{i}\)G_{s}(x,y)\,dx\,dy \nn\\
 &&-16h \int \(\int_{\{\sum_{i=1}^{3}t_{i}\leq t\}}   \(\wt D_{\pi'_{1},t}(x,y)+\wt D_{\pi'_{2},t}(x,y)\) \prod_{i=1}^{3}\,dt_{i}\,dx\,dy\)G_{s}(x,y)\nn\\
 &&+ 16h^{2}\int \(\int_{\{t_{1}+t_{2}\leq t\}} \(A_{t_{1},t_{2}}(x,y)+A_{t_{1},t_{2}}(y,x)\)\,dt_{1}\,dt_{2}\)G_{s}(x,y)\,dx\,dy\nonumber\\
 &&=O\((s\wedge t)^{\ep}h^{4+\ep}\).\label{bd.2a}
  \eea

  We begin by rewritting (\ref{bd.2a}). By symmetry it suffices to show that
    \bea
  &&
8\int\(\int_{\{\sum_{i=1}^{4}t_{i}\leq t\}}  D_{\pi_{1},t}(x,y) \prod_{i=1}^{4}\,dt_{i}\)     G_{s}(x,y)\,dx\,dy\label{bd.3a}\\
&&-32h\int \(\int_{\{\sum_{i=1}^{3}t_{i}\leq t\}}  \wt D_{\pi'_{1},t}(x,y) \prod_{i=1}^{3}\,dt_{i}\)     G_{s}(x,y)\,dx\,dy\nn\\
 &&+ 32h^{2}\int\(\int_{\{t_{1}+t_{2}\leq t\}}  A_{t_{1},t_{2}}(x,y) \,dt_{1}\,dt_{2}\)     G_{s}(x,y)\,dx\,dy=O\((s\wedge t)^{\ep}h^{4+\ep}\).\nn
  \eea
    Using the above expressions for $D_{\pi_{1},t}(x,y), \wt D_{\pi'_{1},t}(x,y)$ and relabeling the $t_{i}'s $ this is equivalent to showing that 
      \bea
  &&
32\int\(\int_{\{\sum_{i=1}^{4}t_{i}\leq t\}}  A_{t_{1},t_{2}}(x,y)  f_{t_{3}}(h) f_{t_{4}}(h)\prod_{i=1}^{4}\,dt_{i}\)     G_{s}(x,y)\,dx\,dy\label{bd.3ab}\\
&&-64h\int \(\int_{\{\sum_{i=1}^{3}t_{i}\leq t\}} A_{t_{1},t_{2}}(x,y)  f_{t_{3}}(h) \prod_{i=1}^{3}\,dt_{i}\)     G_{s}(x,y)\,dx\,dy\nn\\
 &&+ 32h^{2}\int\(\int_{\{t_{1}+t_{2}\leq t\}}  A_{t_{1},t_{2}}(x,y) \,dt_{1}\,dt_{2}\)     G_{s}(x,y)\,dx\,dy=O\((s\wedge t)^{\ep}h^{4+\ep}\).\nn
  \eea
This comes down to making precise  the intuitive notion that 
  $ f_{r}(h)$ is $h$ times a delta-function in $r$, (in which case the left hand side would vanish).
  
  To this end we note 
 \begin{equation}
 \int_{0}^{\ff} f_{r}(h)\,dr= \int_{0}^{\ff} (p_{r}(0)-p_{r}(h))\,dr=h\label{bd.6}
 \end{equation}
 and   for any $\de>0$
\begin{equation}
 \int_{\de}^{\ff} f_{r}(h)\,dr=\int_{\de}^{\ff}  {1-e^{-h^{2}/2r} \over \sqrt{2\pi r}} \,dr\leq \int_{\de}^{\ff}  {  h^{2}/2r\over \sqrt{2\pi r}} \,dr=O(h^{2}/\sqrt{\de}).\label{kacv.7aa}
\end{equation} 

We also note that
\begin{eqnarray}
&&\int_{\{t-2h^{\ep'}\leq  t_{1}+t_{2}\leq t \}}p_{t_{1}}(x)\, p_{t_{2}}(y-x)\,
 \,dt_{1}\,dt_{2}
\label{bd.5a}\\
&&\leq   c\int_{\{t-2h^{\ep'}\leq  t_{1}+t_{2}\leq t \}}{1 \over \sqrt{t_{1}}}\, {1 \over \sqrt{t_{2}}}\,
 \,dt_{1}\,dt_{2}\leq Ct^{2/3}h^{\ep'/4}. \nonumber
\end{eqnarray}
  
  We then write
  \begin{eqnarray}
  &&\int_{\{\sum_{i=1}^{4}t_{i}\leq t\}}  A_{t_{1},t_{2}}(x,y)  f_{t_{3}}(h) f_{t_{4}}(h)\prod_{i=1}^{4}\,dt_{i}
  \label{db.10}\\
  &&=\(\int_{\{t_{1}+t_{2}\leq t-2h^{\ep'}\}}  A_{t_{1},t_{2}}(x,y)  \,dt_{1}\,dt_{2} \)\(\int_{0}^{h^{\ep'}}f_{r}(h)\,dr\)^{2}  \nonumber\\
  &&+\int_{ C(t,h)}  A_{t_{1},t_{2}}(x,y)  f_{t_{3}}(h) f_{t_{4}}(h)\prod_{i=1}^{4}\,dt_{i}   \nonumber
  \end{eqnarray}
  where 
 \bea
 &&
 C(t,h)=\{\sum_{i=1}^{4}t_{i}\leq t\}-\{ t_{1}+t_{2}\leq t-2h^{\ep'}\}\times \{t_{3},t_{4}\leq h^{\ep'}\}\label{bd.5}\\
 &&\hspace{.6 in}\subseteq \([0,t]^{4}\cap \{t_{3},t_{4}\leq h^{\ep'}\}^{c}\)\cup 
 \{t-2h^{\ep'}\leq  t_{1}+t_{2}\leq t \}.\nn
 \eea 
  Using (\ref{bd.6})-(\ref{bd.5a}) we see that for $\ep'$ small
  \begin{eqnarray}
  &&\(\int_{\{t_{1}+t_{2}\leq t-2h^{\ep}\}}  A_{t_{1},t_{2}}(x,y)  \,dt_{1}\,dt_{2} \)\(\int_{0}^{h^{\ep'}}f_{r}(h)\,dr\)^{2}
  \label{db.11}\\
  &&=h^{2} \int_{\{t_{1}+t_{2}\leq t\}}  A_{t_{1},t_{2}}(x,y) \,dt_{1}\,dt_{2}+O(t^{2/3}h^{2+\ep'/4})    \nonumber
  \end{eqnarray}
  and
\begin{equation}
\int_{ C(t,h)}  A_{t_{1},t_{2}}(x,y)  f_{t_{3}}(h) f_{t_{4}}(h)\prod_{i=1}^{4}\,dt_{i} =O(t^{2/3}h^{2+\ep'/4}).\label{db.12}
\end{equation}
A similar analysis applies to  the second term in (\ref{bd.3ab}). Then taking $\ep'=8\ep$ and using  (\ref{gbd.2}) completes the proof of 
 (\ref{bd.3ab}). 
 
 We have now dealt with all terms coming from $I_{2}, I_{3}$ and it only remains to consider the contribution of non-bound variables to $I_{1}$. We will show that this is
  \begin{eqnarray}
  && 32 h^{4}E\(\int  ( L^{ x}_{t})^{2}  \wt L^{ x}_{s} \,dx\)+ O\((s\wedge t)^{\ep}h^{4+\ep}\).   \label{bd.2c}
  \end{eqnarray}
  
 The proof of (\ref{bd.2c}) follows closely the proof of Lemma \ref{lem-3.1j}. The main contribution comes from $\pi=(x,y,x,y)$ or $ (y,x,y, x)$ and $a=e$. Taking $\pi=(x,y,x,y)$   and $a=e$ we have  
  \begin{eqnarray}
  && 4  \int\( \int_{\{\sum_{i=1}^{4}t_{i}\leq t\}}p_{t_{1}}( x ) \De^{h}\De^{-h}p_{t_{2}}(y- x )p_{t_{3}}(y- x )\De^{h}\De^{-h}p_{t_{4}}(y- x )\right.\nonumber\\
  &&  \left. \hspace{3 in}\prod_{i=1}^{4}\,dt_{i}\)     G_{s}(x,y)\,dx\,dy
\nn   \label{bd.14}
  \end{eqnarray}
 
Since as before 
  \begin{eqnarray}
  && |\int_{\{\sum_{i=1}^{4}t_{i}\leq t\}}p_{t_{1}}( x ) \De^{h}\De^{-h}p_{t_{2}}(y- x )p_{t_{3}}(y- x )\De^{h}\De^{-h}p_{t_{4}}(y- x )\prod_{i=1}^{4}\,dt_{i} |    
\nn\\
&&\hspace{1 in}\leq  ct^{\ep/2} u^{1-\ep}_{t}(x)u_{t}(y-x)w^{2}_{t}(y-x), \label{bd.15}
  \end{eqnarray}
we see  that up to  terms that are  $O\((s\wedge t)^{\ep}h^{4+\ep}\)$ we can replace $G_{s}(x,y)$ in (\ref{bd.14})
by
\begin{equation}
\int_{\{s_{1}+s_{2}\leq s\}}\(p_{s_{1}}(x)\,\De^{h}\De^{-h}p_{s_{2}}(y-x)+p_{s_{1}}(y)\,\De^{h}\De^{-h}p_{s_{2}}(x-y)\)\,ds_{1}\,ds_{2}.\label{bd.16}
\end{equation}

Thus consider
  \begin{eqnarray}
  && 4  \int\( \int_{\{\sum_{i=1}^{4}t_{i}\leq t\}}p_{t_{1}}( x ) \De^{h}\De^{-h}p_{t_{2}}(y- x )p_{t_{3}}(y- x )\De^{h}\De^{-h}p_{t_{4}}(y- x )\right.\nonumber\\
  &&  \left. \hspace{,5 in}\prod_{i=1}^{4}\,dt_{i}\)\( \int_{\{s_{1}+s_{2}\leq s\}} p_{s_{1}}(x)\,\De^{h}\De^{-h}p_{s_{2}}(y-x)   \,ds_{1}\,ds_{2} \)\,dx\,dy
\nn   \label{bd.17}
  \end{eqnarray}
  It now follows as in the proof of Lemma \ref{lem-3.1j} that up 
to the error terms allowed in (\ref{bd.2c}) this is equal to 
\begin{equation}
16h^{4}\int\( \int_{\{ t_{1}+ t_{2}\leq t\}}p_{t_{1}}( x ) p_{t_{2}}(0) \,dt_{1} \,dt_{2}\)
\( \int_{\{ s_{1} \leq s\}}p_{s_{1}}( x )  \,ds_{1} \)\,dx.\label{bd.18}
\end{equation}
The second term (\ref{bd.16}) gives the same  contribution since up to another error term we can replace $p_{s_{1}}(y)$ by $p_{s_{1}}(x)$. There is a similar contribution from  $ \pi= (y,x,y, x)$. Thus altogether we have 
\begin{equation}
64\int\( \int_{\{ t_{1}+ t_{2}\leq t\}}p_{t_{1}}( x ) p_{t_{2}}(0) \,dt_{1} \,dt_{2}\)
\( \int_{\{ s_{1} \leq s\}}p_{s_{1}}( x )  \,ds_{1} \)\,dx.\label{bd.18}
\end{equation}
Since by Kac's moment formula
\begin{eqnarray}
&&E\(\int  ( L^{ x}_{t})^{2}  \wt L^{ x}_{s} \,dx\)
\label{bd.19}\\
&& =2 \int\( \int_{\{ t_{1}+ t_{2}\leq t\}}p_{t_{1}}( x ) p_{t_{2}}(0) \,dt_{1} \,dt_{2}\)
\( \int_{\{ s_{1} \leq s\}}p_{s_{1}}( x )  \,ds_{1} \)\,dx \nonumber
\end{eqnarray}
we obtain the main contribution to (\ref{bd.2c}). The fact that all remaining $\pi,a$   give error terms is now easy and left to the reader.\qed

\section{Proof of Lemmas \ref{lem-vprop}--\ref{lem-big}}\label{sec-Prooflemvprop}

{\bf  Proof of Lemma \ref{lem-vprop}}
Since   
\begin{eqnarray} \lefteqn{
\De_{ x}^{ h}\De_{ y}^{ h} u^{\al}(x-y)\label{1.8w}}\\ && =
\{u^{\al}(x-y)-u^{\al}(x-y-h)\}\nonumber-
\{u^{\al}(x-y+h)-u^{\al}(x-y)\}\nonumber
\end{eqnarray}  
we have 
\begin{eqnarray} &&
\De_{ x}^{ h}\De_{ y}^{ h} u^{\al}(x-y)\Bigg\vert_{ y=x} =
\{u^{\al}(0)-u^{\al}(-h)\}-
\{u^{\al}(h)-u^{\al}(0)\}\nn\\ &&\hspace{ 1in}=2(u^{\al}(0)-u^{\al}(h) )=
2\({1-e^{-\sqrt{2\al}\,h} \over \sqrt{2\al}}\),\label{1.8a}
\end{eqnarray}
which gives (\ref{1.8}).

  To obtain  (\ref{1.3x}) we note that
\begin{equation}
\De_{x}^{ h}\,u^{\al}(x)=\({e^{-\sqrt{2\al}|x+h|} -e^{-\sqrt{2\al}|x|} \over \sqrt{2\al}}\).\label{pot.3ow}
\end{equation}
Therefore
\bea
|\De_{x}^{ h}\,u^{\al}(x)|&\le&{ e^{-\sqrt{2\al}|x |}\over \sqrt{2\al}}\left| e^{ \sqrt{2\al}(|x|-|x+h|)}-1\right| \label{pot.3owa}\\
&\le&  e^{-\sqrt{2\al}|x |}\( ||x|-|x+h||+O( ||x|-|x+h||^{2})\) \nn
\eea
which gives  (\ref{1.3x}), (since we allow $C$ to depend on $\al$.)

To obtain  (\ref{1.3y}) we simply note that 
\begin{equation}
|\De^{ h}\De^{ -h}\,u^{\al}( x)|=|2u^{\al}( x)-u^{\al}( x+h)-u^{\al}( x-h)|\leq 2v^{\al}( x)\label{ff.4}
\end{equation}
where we used the fact that $u^{\al}( x)$ is an even function. The first part of (\ref{1.3y}) then follows from (\ref{1.3x}).
When $|x|\geq h$ we have  
\begin{eqnarray} 
\De^{ h}\De^{ -h}\,u^{\al}( x)&=& 2u^{\al}( x)-u^{\al}( x+h)-u^{\al}( x-h)\label{1.26gd}\\ 
&=& u^{\al}( x)\(2-e^{-\sqrt{2\al}\,h}-e^{\sqrt{2\al}\,h}\)\nn .
\end{eqnarray}
The statement in  (\ref{1.30gb}) follows trivially from (\ref{1.3y}).

For (\ref{1.30g}) we note that for $|x|\leq h$
\begin{eqnarray} 
\De^{ h}\De^{ -h}\,u^{\al}( x)&=& 2u^{\al}( x)-u^{\al}( x+h)-u^{\al}( x-h)\label{1.26g}\\ &=& (1-u^{\al}(
x+h))+(1-u^{\al}( x-h))-2(1-u^{\al}( x))\nn\\ &=& | x+h|+ | x-h|-2 | x|+O( h^{ 2}).\nn
\end{eqnarray} 
 When $0\leq x\leq h$ we therefore have 
\begin{equation}
\De^{ h}\De^{ -h}\,u^{\al}( x)=x+h+h-x-2x+O( h^{ 2})=( 2+O( h))(h-x).\label{1.27g}
\end{equation} 
 Consequently
\bea
\int_{0}^{ h} \(\De^{ h}\De^{ -h}\,u^{\al}(x)\)^{q} \,dx \nn
&=& ( 2^{q}+O( h))\int_{0}^{
h}(h-x)^{q}\,dx\\
&=&(  2^{q}/(q+1)+O( h))h^{ q+1}.\label{1.28g}
\eea
 Similarly, when $-h\leq x\leq 0$ it follows from  (\ref{1.26g}) that
\begin{equation}
\De^{ h}\De^{ -h}\,u^{\al}( x)=h-x+x+h+2x+O( h^{ 2})=( 2+O( h))(h+x),\label{1.29g}
\end{equation} 
 Consequently
\bea \int_{-h}^{ 0}\(\De^{ h}\De^{ -h}\,u^{\al}(x)\)^{q} \,dx&=& ( 2^{q}+O( h))\int_{-h}^{
0}(h+x)^{q}\,dx\nn\\
&=&(2^{q}/(q+1)+O( h))h^{ q+1}.\label{1.30gx}
\eea
Using  (\ref{1.28g}), (\ref{1.30gx}) and   (\ref{1.30gb}) we get (\ref{1.30g}).

To obtain (\ref{li.13}) we write 
\begin{eqnarray}
&&\int |\De^{ h}\De^{- h}\,u^{\al}(y) |^{q}\,dy 
\label{li.13a}\\
&&\qquad= \int_{|y|\leq h} |\De^{ h}\De^{- h}\,u^{\al}(y) |^{q}\,dy   + \int_{|y|\geq h} |\De^{ h}\De^{- h}\,u^{\al}(y) |^{q}\,dy \nonumber\\
&&\qquad\leq Ch^{q}\int_{|y|\leq h} 1\,dy   + Ch^{2q}\int_{|y|\geq h}  u^{\al}(y) \,dy=O( h^{
q+1}) ,\nonumber
\end{eqnarray}
where for the last line we use (\ref{1.3y}).
 \qed
 
 {\bf  Proof of Lemma \ref{lem-vpropt}}
 It follows from the fact that $p_{r}(x)\leq p_{r}(y)$ for all $r$ if $|y|\leq |x|$,  (\ref{pot.1w}), and (\ref{1.3x}) that
 \begin{eqnarray}
  \int_{0}^{T} |\De ^{ h}\,p_{t}(x)|\,dt
& \leq &e^{T/2}\int_{0}^{\ff}e^{-t/2} |\De ^{ h}\,p_{t}(x)|\,dt  \nonumber\\
 & =& e^{T/2}\bigg| \De ^{ h}\(\int_{0}^{\ff}e^{-t/2}\,p_{t}(x)\,dt\)\bigg | \label{9.7}\\
 & =&e^{T/2}|\De ^{ h}\,u^{1/2}(x)|\leq C_{T} h\, e^{-|x|}.  \nonumber
 \end{eqnarray}
This gives (\ref{9.3x}).

 For (\ref{9.3w}), we note that 
 \bea
\bigg|{d^{2} \over dx^{2}}p_{t}(x )\bigg|&=&\bigg|{x^{2}/t-1 \over t\sqrt{2\pi t}}e^{-x^{2}/2t}\bigg| \label{9.9}\\
&\le& {C\over   t^{3/2}}\({x^{2}   \over 2t}+1\)e^{-x^{2}/2t}\le \frac{C}{t^{3/2}}e^{-x^{2}/4t},\nn
\eea
since $\sup_{s>0}se^{-s }<\ff$.
We use     this and  Taylor's theorem to see that  for some   $0\leq h_{t}',h_{t}''\leq h$,
\begin{eqnarray}
 |\De^{ h}\De^{ -h} p_{t}(x )|
\label{9.11} & =&|2 p_{t}(x )-p_{t}(x+h )-p_{t}(x-h )| \label{2.26k}\\
&  =&{h^{2} \over 2}\bigg|{d^{2} \over dx^{2}}p_{t}(x+h_{t}')+{d^{2} \over dx^{2}}p_{t}(x-h_{t}'')\bigg|  \nonumber\\
& \leq &{Ch^{2} \over t^{3/2}}
\(e^{-(x+h_{t}')^{2}/4t}+e^{-(x+h_{t}'')^{2}/4t} \)  \nonumber.
\eea
Therefore, when   $|x|\geq 2h$,    
\begin{equation}
 |\De^{ h}\De^{ -h} p_{t}(x )|   \leq  {Ch^{2} \over t^{3/2}}
 e^{-x^{2}/16t}.
   \end{equation}
 Consequently, when   $|x|\geq 2h$, 
\bea
\int_{0}^{T} |\De^{ h}\De^{ -h} p_{t}(x )|\,dt&\leq& Ch^{2} \int_{0}^{T} { e^{-x^{2}/16t}\over t^{3/2}}\,dt\nn\\
&\leq & Ch^{2} e^{-x^{2}/32T} \int_{0}^{\ff} { e^{-x^{2}/32t}\over t^{3/2}}\,dt
\label{ggg}\\
 &=& Ch^{2}{e^{-x^{2}/32T} \over |x| }\int_{0}^{\ff} { e^{-1/32t}\over t^{3/2}}\,dt\leq C_{T}h^{2}{e^{-x^{2}/32T} \over |x| }\nn,
\eea
which proves (\ref{9.3w}).

Using  (\ref{9.3x}) and (\ref{9.3w}) we see that  
 \begin{eqnarray}
\lefteqn{\int  w_{T}^{q}(x)\,dx
\label{2.28}}\\
 &&=\int_{|x|\leq 2h}  \(\int_{0}^{T} |\De^{ h}\De^{ -h} p_{t}(x )|\,dt\)^{q}\,dx\nn\\
 &&\qquad+ \int_{|x|\geq 2h}  \(\int_{0}^{T} |\De^{ h}\De^{ -h} p_{t}(x )|\,dt\)^{q}\,dx  \nonumber\\
 &&\leq 4\int_{|x|\leq 2h}  \(\int_{0}^{T} |\De^{ h}  p_{t}(x )|\,dt\)^{q}\,dx\nn\\
 &&\qquad+ \int_{|x|\geq 2h}  \(\int_{0}^{T} |\De^{ h}\De^{ -h} p_{t}(x )|\,dt\)^{q}\,dx  \nonumber\\
 &&\leq C_{T}\int_{|x|\leq 2h}  h^{q}\,dx+ C_{T}h^{2q}\int_{|x|\geq 2h}  {1 \over |x|^{q}}\,dx\leq C_{T}h^{q+1}  \nn,
 \end{eqnarray}
which gives us (\ref{9.30g}). 

For (\ref{9.30gb}) we   note that when $h\le 1/4$, $\sqrt h\ge 2h$. Therefore, it folows from   (\ref{9.3w})   that
 \bea
 &&
\int_{|x|\geq \sqrt{h}}w_{T}^{q}(x)\,dx \label{9.14 }\\
&&\leq  C_{T}h^{2q}\int_{|x|\geq \sqrt{h}}  {1 \over |x|^{q}}\,dx \leq C_{T}h^{3q/2+1/2}.\nn
 \eea 
 
Finally, to obtain (\ref{9.13t}) we use  (\ref{9.3x}) and (\ref{9.3w}) to see that
\begin{eqnarray}
&&\int w_{T} (x) \,dx
\label{9.1}\\&&
 =\int_{|x|\leq h} \int_{0}^{T} |\De^{ h}\De^{ -h} p_{t}(x )|\,dt \,dx+\int_{|x|\geq h}  \int_{0}^{1} |\De^{ h}\De^{ -h} p_{t}(x )|\,dt \,dx   \nonumber\\
&&\leq 2\int_{|x|\leq h} \int_{0}^{T} |\De^{ h}  p_{t}(x )|\,dt \,dx+\int_{|x|\geq h}  \int_{0}^{1} |\De^{ h}\De^{ -h} p_{t}(x )|\,dt \,dx   \nonumber\\
&&\leq C_{T}\int_{|x|\leq h} h\,dx+C_{T}h^{2}\int_{|x|\geq h} {e^{-x^{2}/8}\over |x|} \,dx\leq C_{T}h^{ 2}\log h.   \nonumber
\end{eqnarray}
\qed

  \begin{remark}
 {\rm   Using Remark 2.1 and (\ref{9.7}) it is easy to check that we obtain the analog of (\ref{9.3x} ) for all $|h|\leq 1$  if on the right hand side   we replace $h$ by $|h|$.}
  \end{remark}

  {\bf  Proof of Lemma \ref{lem-vpropd}} The proof of (\ref{d9.300}) is immediate. (\ref{d9.3w})
  follows from (\ref{9.11}), and a similar application of the mean value theorem gives (\ref{d9.3x}).
  \qed

{\bf  Proof of Lemma \ref{lem-big}}  Using $2-e^{ihp } -e^{-ihp }=2-2\cos (hp)=4\sin^{2} (p h/2)$ we can write
\bea
&&
\int_{0}^{\ff} \De^{ h}\De^{ -h} p_{t}(x )\,dt \label{big.3}\\
&&={1 \over 2\pi}\int_{0}^{\ff} \int e^{ipx}(2-e^{ihp } -e^{-ihp })e^{-tp^{2}/2} \,dp\,dt\nn\\
&&={4 \over 2\pi}\int_{0}^{\ff} \int e^{ipx}\sin^{2} (p h/2)e^{-tp^{2}/2} \,dp\,dt\nn\\
&&={8 \over 2\pi}  \int e^{ipx}{\sin^{2} (p h/2) \over p^{2}}  \,dp. \nn
\eea
Similarly
\begin{equation}
\int_{0}^{h} \De^{ h}\De^{ -h} p_{t}(x )\,dt={8 \over 2\pi}  \int e^{ipx}{\sin^{2} (p h/2) \over p^{2}}\(1-e^{-hp^{2}/2}\)  \,dp\label{big.4}
\end{equation}
and
\begin{equation}\qquad
\De^{ h}\De^{ -h} u^{1/2}(x )=\int_{0}^{\ff}e^{-t/2} \De^{ h}\De^{ -h} p_{t}(x )\,dt={8 \over 2\pi}  \int e^{ipx}{\sin^{2} (p h/2) \over 1+p^{2}}  \,dp.\label{big.5}
\end{equation}

 Using (\ref{big.3}) and the Fourier inversion formula we see that
 \begin{eqnarray}
 &&\int  \(\int_{0}^{\ff} \De^{ h}\De^{ -h} p_{t}(x )\,dt\)^{q}\,dx
 \label{big.7}\\
 &&=\({8 \over 2\pi}\)^{q}\int  \( \int e^{ipx}{\sin^{2} (p h/2) \over p^{2}}  \,dp\)^{q}\,dx   \nonumber\\
 &&=\({8 \over 2\pi}\)^{q}\int  \( \int e^{ix\sum_{j=1}^{q}p_{j}}\prod_{j=1}^{q}{\sin^{2} (p_{j} h/2) \over p_{j}^{2}}  \,dp_{j}\)\,dx   \nonumber\\
 &&=\({8 \over 2\pi}\)^{q} \int \(\int e^{ix\sum_{j=2}^{q}p_{j}}\(\int e^{ix p_{1}}{\sin^{2} (p_{1} h/2) \over p_{1}^{2}}  \,dp_{1}     \)\,dx\)\nn\\
 &&\hspace{3 in}\prod_{j=2}^{q}{\sin^{2} (p_{j} h/2) \over p_{j}^{2}}  \,dp_{j}  \nonumber\\
 &&= {8^{q} \over (2\pi)^{q-1}}  \int  {\sin^{2} (p_{1} h/2) \over p_{1}^{2}}\prod_{j=2}^{q}{\sin^{2} (p_{j} h/2) \over p_{j}^{2}}  \,dp_{j}   \nonumber
 \end{eqnarray}
 where now $p_{1}=\sum_{j=2}^{q}p_{j}$. Scaling in $h$ we then obtain
 \bea
 &&
\int  \(\int_{0}^{\ff} \De^{ h}\De^{ -h} p_{t}(x )\,dt\)^{q}\,dx \label{big.8}\\
&&=
 {8^{q} h^{q+1}\over (2\pi)^{q-1}}  \int  {\sin^{2} (p_{1} /2) \over p_{1}^{2}}\prod_{j=2}^{q}{\sin^{2} (p_{j} /2) \over p_{j}^{2}}  \,dp_{j}.\nn
\eea
Similarly we see that
 \bea
 &&\qquad
\int  \(\int_{0}^{h} \De^{ h}\De^{ -h} p_{t}(x )\,dt\)^{q}\,dx \label{big.9}\\
&&=
 {8^{q} h^{q+1}\over (2\pi)^{q-1}}  \int  {\sin^{2} (p_{1} /2) \over p_{1}^{2}}
 \(1-e^{-p_{1}^{2}/2h}\)\prod_{j=2}^{q}{\sin^{2} (p_{j} /2) \over p_{j}^{2}}\(1-e^{-p_{j}^{2}/2h}\)  \,dp_{j}\nn
\eea
and
 \bea
 &&
\int  \(  \De^{ h}\De^{ -h} u^{1/2}(x ) \)^{q}\,dx \label{big.10}\\
&&=
 {8^{q} h^{q+1}\over (2\pi)^{q-1}}  \int  {\sin^{2} (p_{1} /2) \over h^{2}+p_{1}^{2}}\prod_{j=2}^{q}{\sin^{2} (p_{j} /2) \over h^{2}+p_{j}^{2}}  \,dp_{j}.\nn
\eea
Using the fact that ${\sin^{2} (p  /2) \over p^{2}}$ is bounded and 
\begin{equation}
\int e^{-p^{2}/2h}\,dp=Ch^{1/2}\label{big.11}
\end{equation}
 our Lemma follows from comparing (\ref{big.8})-(\ref{big.10}) with (\ref{1.30g}).\qed

 {\bf  Acknowledgment.} I would like to thank David Nualart for pointing out an error in the first draft of this paper.

\def\noopsort#1{} \def\printfirst#1#2{#1}
\def\singleletter#1{#1}
      \def\switchargs#1#2{#2#1}
\def\bibsameauth{\leavevmode\vrule height .1ex
      depth 0pt width 2.3em\relax\,}
\makeatletter
\renewcommand{\@biblabel}[1]{\hfill#1.}\makeatother
\newcommand{\bysame}{\leavevmode\hbox to3em{\hrulefill}\,}

\bigskip
\noindent
\begin{tabular}{lll} 
      & Jay Rosen \\
      & Department of Mathematics \\
     &College of Staten Island, CUNY \\
     &Staten Island, NY 10314 \\ &jrosen30@optimum.net  
\end{tabular}

\end{document}